\documentclass[11pt]{amsart}

\usepackage{DKstyle}
\usepackage{bm}
\usepackage{float}
\usepackage{longtable}
\usepackage[usenames,dvipsnames]{color}

\usepackage[active]{srcltx}
\usepackage{hyperref}

\definecolor{blue}{rgb}{0,0,1}
\definecolor{red}{rgb}{1,0,0}
\definecolor{green}{rgb}{0,.6,.2}
\definecolor{purple}{rgb}{1,0,1}
\long\def\red#1\endred{\textcolor{red}{#1}}
\long\def\blue#1\endblue{\textcolor{blue}{#1}}
\long\def\purple#1\endpurple{\textcolor{purple}{ #1}}
\long\def\green#1\endgreen{\textcolor{green}{#1}}

\newcommand{\matr}[4]{\left( \begin{matrix} #1 & #2 \\ #3 & #4\end{matrix} \right) }
\newcommand{\smatr}[4]{\left( \begin{smallmatrix} #1 & #2 \\ #3 & #4\end{smallmatrix} \right) }

\newcommand{\AB}{\operatorname{AB}}
\newcommand{\TT}{\operatorname{\hspace{1pt} T \hspace{-1pt} T}}
\newcommand{\Nr}{\operatorname{Nr}}
\newcommand{\Gh}{\operatorname{Gh}}

\newcommand{\intl}{\operatorname{int}}

\def\scrA{{\mathcal A}}
\def\scrB{{\mathcal B}}
\def\scrD{{\mathcal D}}
\def\scrF{{\mathcal F}}
\def\scrL{{\mathcal L}}
\def\scrI{{\mathcal I}}
\def\scrO{{\mathcal O}}
\def\scrP{{\mathcal P}}

\def\scrT{{\mathcal T}}

\def\scrW{{\mathcal W}}

\def\S{\operatorname{S{}}}

\def\covol{\operatorname{covol}}

\def\vece{{\text{\boldmath$e$}}}

\def\vecm{{\text{\boldmath$m$}}}

\def\vecv{{\text{\boldmath$v$}}}

\def\vecw{{\text{\boldmath$w$}}}
\def\vecx{{\text{\boldmath$x$}}}

\def\vecy{{\text{\boldmath$y$}}}
\def\vecz{{\text{\boldmath$z$}}}

\def\vectheta{{\text{\boldmath$\theta$}}}

\newcommand{\bn}{\mathbf{0}}
\newcommand{\GaG}{\Gamma\bs \mathrm{G}}

\newcommand{\tF}{\widetilde{F}}

\newcommand{\ttheta}{\tilde{\theta}}

\newcommand{\ts}{\tilde{s}}

\newcommand{\hvecv}{\widehat{\vecv}}
\newcommand{\Mod}[1]{\ (\mathrm{mod}\ #1)}
\newcommand{\vertiii}[1]{{\left\vert\kern-0.25ex\left\vert\kern-0.25ex\left\vert #1
    \right\vert\kern-0.25ex\right\vert\kern-0.25ex\right\vert}}
 \makeatletter
\newcommand*{\rom}[1]{\expandafter\@slowromancap\romannumeral #1@}
\makeatother
\theoremstyle{plain}

\newcommand{\col}{\: : \:}
\renewcommand{\bs}{\backslash}

\newcommand{\tff}{\widetilde{\ff}}
\newcommand{\ve}{\varepsilon}

\usepackage{caption}
\usepackage[labelfont=rm]{subcaption}

\usepackage{comment}


\dedicatory{}
\commby{}

\begin{document}
\title[large gaps between directions in planar quasicrystals]
{Asymptotic estimates of large gaps between directions in certain planar quasicrystals}

\author[G.\ Hammarhjelm]{Gustav Hammarhjelm}

\author[A.\ Str\"ombergsson]{Andreas Str\"ombergsson}
\address{Department of Mathematics, Uppsala University, Box 480, SE-75106, Uppsala, Sweden}
\email{astrombe@math.uu.se}

\author[S.\ Yu]{Shucheng Yu}
\address{School of Mathematical Sciences, University of Science and Technology of China (USTC), 230026, Hefei, China}
\email{yusc@ustc.edu.cn}

\date{May 2, 2025}


\thanks{All three authors were supported by the Knut and Alice Wallenberg Foundation. 
S.Y.\ was also supported by the 
National Key R\&D Program of China No. 2024YFA1015100}

\maketitle

{\centering Dedicated to Gustav Hammarhjelm (1992--2022).\par}

\begin{abstract}
For quasicrystals of cut-and-project type in $\R^d$,
it was proved by Marklof and Str\"ombergsson 
\cite{MarklofStrombergsson2015}
that the limit local statistical properties of the directions
to the points in the set
are described by certain $\SL_d(\R)$-invariant point processes.
In the present paper we make 
a detailed study of the tail asymptotics of
the limiting  
gap statistics of the directions,
for certain specific classes of 
planar quasicrystals.

\end{abstract}

\tableofcontents

\section{Introduction}

\subsection{Gaps between directions in planar point sets} 
Given a locally finite point set $\scrP$ in the plane $\R^2$,
we consider its set of  
\textit{directions},
that is, the set of points $\hvecv:=\|\vecv\|^{-1}\vecv$ on the unit circle $\S_1^1$ \label{p:unitcirc}
as $\vecv$ runs through $\scrP\smallsetminus\{\bn\}$.
For each $R>0$, let $\Delta_R$ be the finite multi-set of directions to points in
$\scrP$ within distance $R$ from the origin,
i.e.,
\begin{align}\label{deltaRDEF}
\Delta_R:=\{\hvecv\col\vecv\in\scrP,\:0<\|\vecv\|\leq R\}.
\end{align}
Throughout the paper we will assume that $\scrP$ has an
asymptotic density $c_\scrP>0$,
meaning that for any bounded set $\scrD\subset\R^2$ with boundary of 
measure zero,
\begin{align}\label{asymptdensityDEF}
\lim_{R\to\infty}\frac{\#(\scrP\cap R\scrD)}{R^2}=c_{\scrP}\Area(\scrD).
\end{align}
It then follows that the multi-set $\Delta_R$ 
has cardinality $\sim c_\scrP\hspace{1pt} \pi R^2$ as $R\to\infty$,
and furthermore that this multi-set 
becomes asymptotically equidistributed along the unit circle,
in the sense that for any arc $I\subset\S_1^1$,
\begin{align*}
\lim_{R\to\infty}\frac{\#(\Delta_R\cap I)}{\#\Delta_R}=\frac{|I|}{2\pi},
\end{align*}
where $|I|$ denotes the length of $I$ (in particular $|\S_1^1|=2\pi$).

In this situation,
it is natural to consider finer questions about the local statistics
of the points in $\Delta_R$.
A particular statistics which has been much studied 
is the distribution of normalized gaps between the points in $\Delta_R$,
as $R\to\infty$.
For example, when $\scrP$ is the set of primitive lattice points in $\Z^2$,
the limiting distribution of normalized gaps was
explicitly determined by Boca, Cobeli and Zaharescu
\cite{BocaCobeliZaharescu2000}.
More generally, if $\scrP$ is an arbitrary lattice (possibly translated) in $\R^d$,
for any $d\geq2$,
the limiting distribution of general local statistics of directions
was proved to exist by Marklof and Str\"ombergsson \cite{MarklofStrombergsson2010};
see also \cite{ElBazMarklofVinogradov2015}
and \cite{KimMarklof2025}
regarding convergence of related moments. 
In particular it was noted in \cite{MarklofStrombergsson2010}
that the explicit limiting distribution
computed in \cite{BocaCobeliZaharescu2000}
remains valid when replacing $\Z^2$ by an arbitrary lattice in $\R^2$;
furthermore, 
when $\scrP$ is an arbitrary 'irrational' translate of a lattice in $\R^2$,
the limiting distribution of gaps between directions
coincides with the gap distribution for the fractional parts of $\sqrt n$
calculated by Elkies and McMullen 
\cite{nEcM2001};
see \cite[{Sec.\ 1.2}]{MarklofStrombergsson2010}
and \cite{Marklofsquareroots2024}.

Athreya and Chaika \cite[{Prop.\ 3.10}]{jAjC2012}
have proved that the limiting distribution of gaps between directions
exists in the case when 
$\scrP$ is the set of holonomy vectors of either the saddle connections or periodic cylinders
on a generic translation surface;
see also 
\cite{jAjCsL2015},
\cite{cUgW2016},
\cite{lKjW2024},
\cite{aS2020a},
\cite{jBtMaMcUhW2021}
and
\cite{OsmanSoutherlandWang2024}
for later studies with more precise results in specific cases.
Also in hyperbolic $n$-space,
for any point set which is the orbit
of a lattice within the group of orientation preserving isometries,
the analogous limiting gap distribution
and also more general local statistics of directions,
has been proved to exist
by Marklof and Vinogradov \cite{jMiV2014};
for related work see 
\cite{fBaPaZ2013},
\cite{dKaK2013},
\cite{mRaS2014}.

\vspace{2pt}

The present paper concerns the case of $\scrP$ being a regular cut-and-project set
(also referred to as a Euclidean model set).
For this case,
the limit distribution of normalized gaps between directions
was proved to exist, and given a complicated but explicit description,
by Marklof and Str\"ombergsson
in \cite{MarklofStrombergsson2015};
see Theorem \ref{gaplimitdistrexistsTHM} below.
See also R\"uhr, Smilansky and Weiss \cite{RuhrSmilanskyWeiss2024}
for closely related investigations,
and El-Baz \cite{dE2017}
for an extension to so called adelic model sets.
The results of  
\cite{MarklofStrombergsson2015}
answered some questions which had been raised
by Baake, G\"otze, Huck and Jakobi in
\cite{mBfGcHtJ2014},
where a numerical investigation was carried out of
the normalized gap distribution
between directions 
for several vertex sets
coming from aperiodic tilings,
some of these being of cut-and-project type and others not.
Building on the results of 
\cite{MarklofStrombergsson2015},
Hammarhjelm \cite{Hammarhjelm2022}
explicitly determined the limit of the \textit{minimal} normalized gap between the
directions to the points in $\scrP$,
for $\scrP$ belonging to either of two families of planar quasicrystals,
including both the Ammann-Beenker point set and the
vertex sets of some rhombic Penrose tilings.
Also in \cite{Hammarhjelm2022},
the asymptotic density of visible points was
explicitly determined for several families of planar quasicrystals,
including the two families just mentioned.

\vspace{2pt}

\begin{figure}
\begin{center}
\framebox{
\begin{minipage}{0.4\textwidth}
\unitlength0.1\textwidth
\begin{picture}(10,9.5)(0,0.8)
\put(0.5,1){\includegraphics[width=0.9\textwidth]{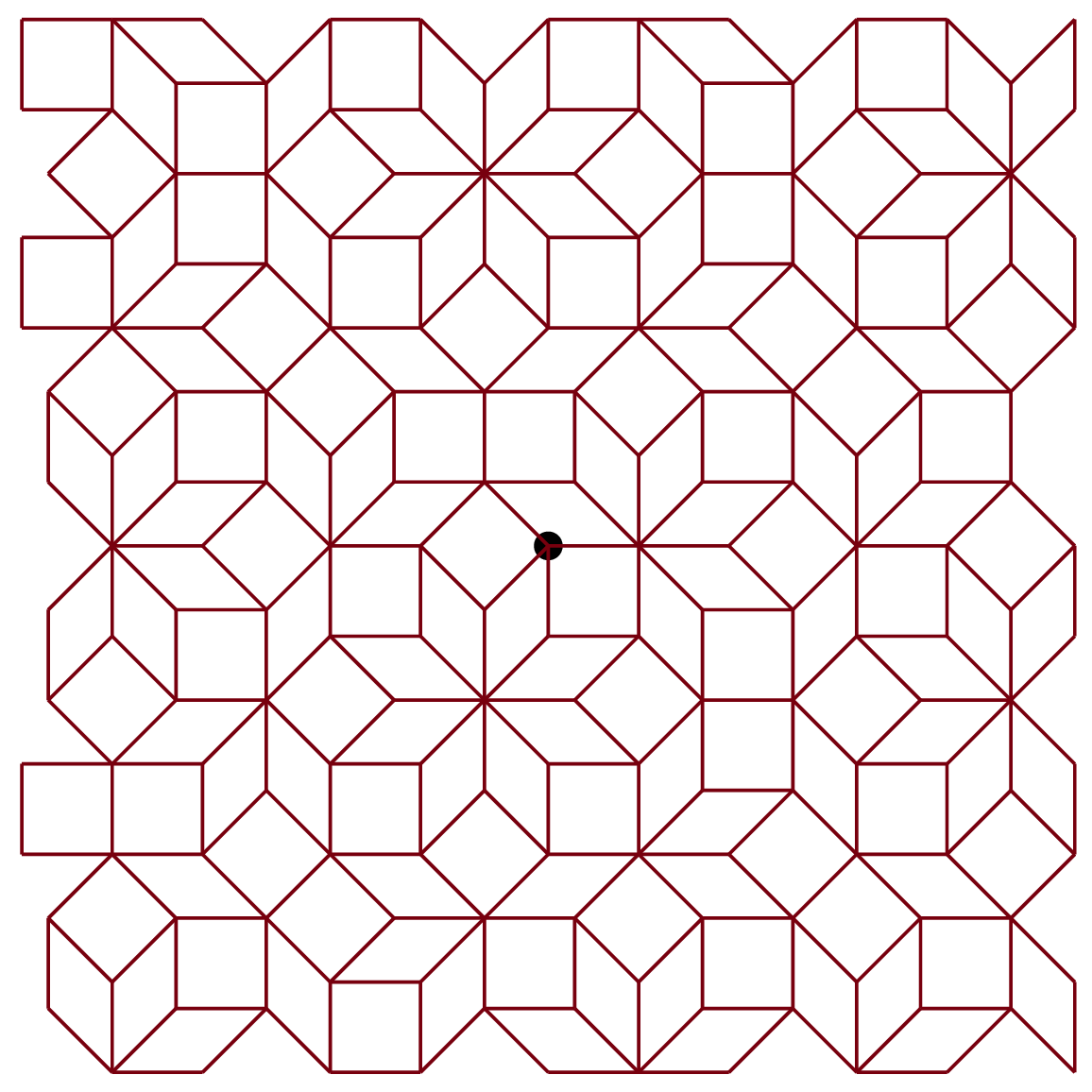}}
\end{picture}
\end{minipage}
}\hspace{6pt}\framebox{
\begin{minipage}{0.4\textwidth}
\unitlength0.1\textwidth
\begin{picture}(10,9.5)(0,0.8)
\put(0.5,1){\includegraphics[width=0.9\textwidth]{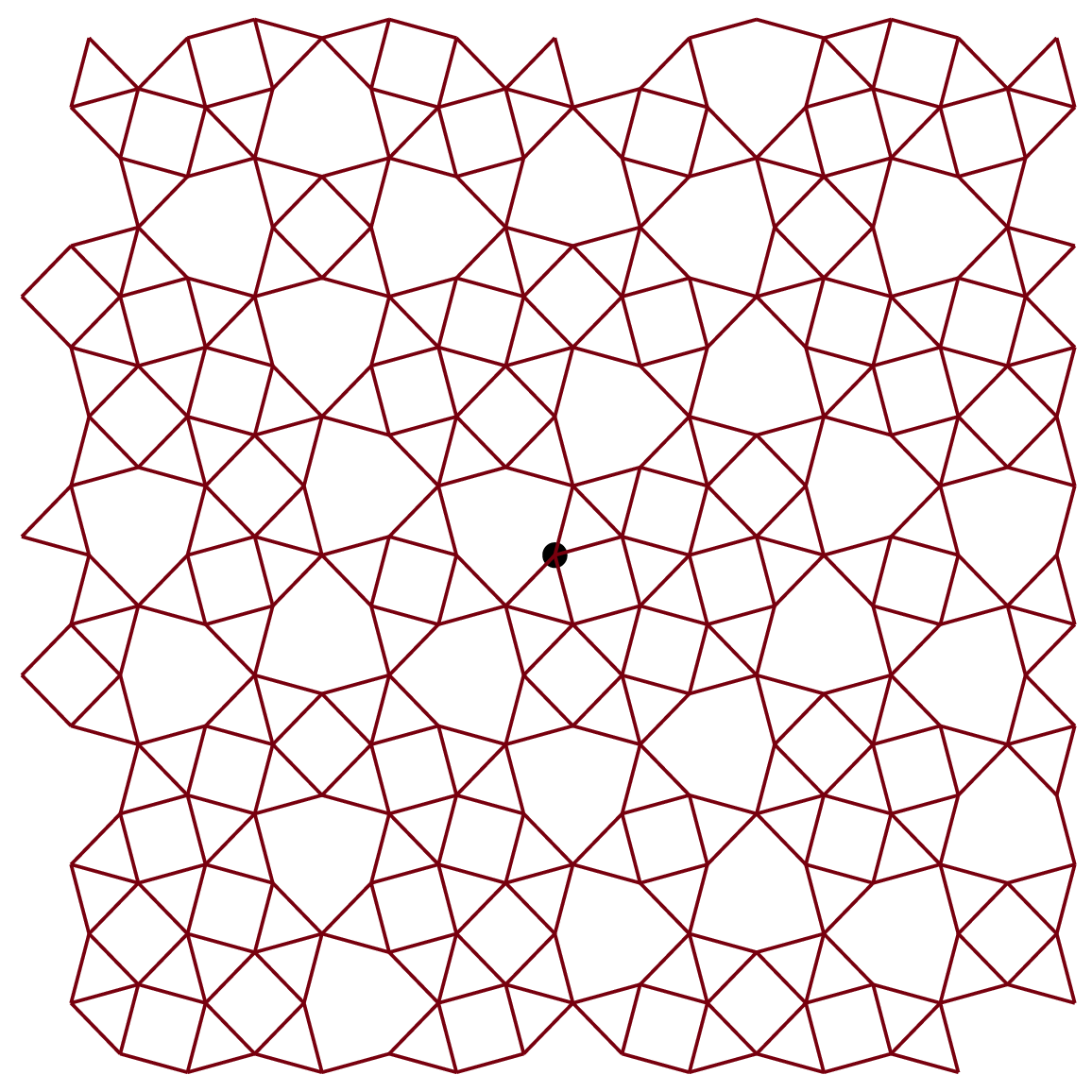}}
\end{picture}
\end{minipage}
}
\\[5pt]
\framebox{
\begin{minipage}{0.4\textwidth}
\unitlength0.1\textwidth
\begin{picture}(10,9.5)(0,0.8)
\put(0.5,1){\includegraphics[width=0.9\textwidth]{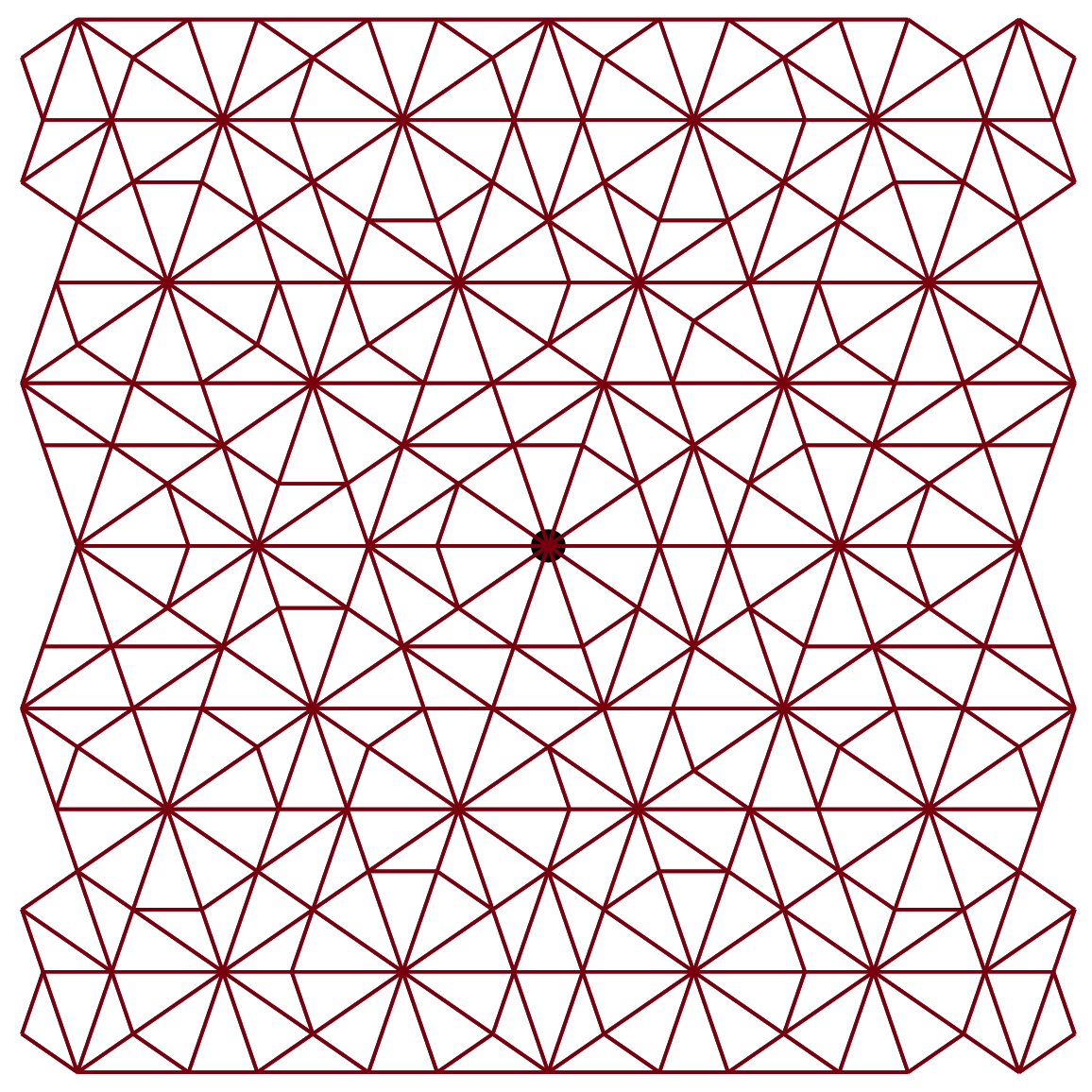}}
\end{picture}
\end{minipage}
}\hspace{6pt}\framebox{
\begin{minipage}{0.4\textwidth}
\unitlength0.1\textwidth
\begin{picture}(10,9.5)(0,0.8)
\put(0.5,1){\includegraphics[width=0.9\textwidth]{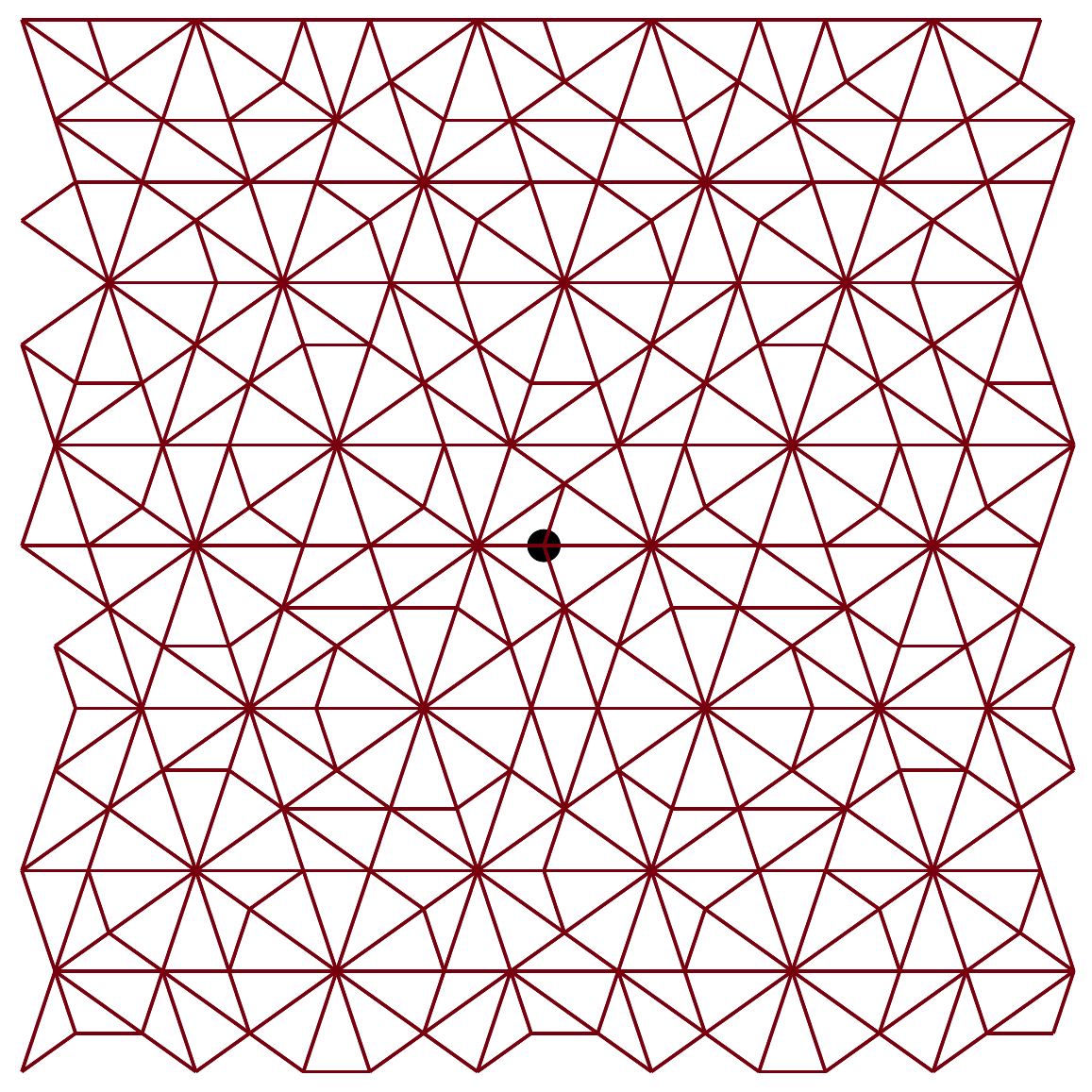}}
\end{picture}
\end{minipage}
}
\end{center}
\caption{Patches of an Ammann-Beenker tiling with vertex set 
$\scrP_{\AB,(0.8,0.1)}$ 
(top left panel),
a G\"ahler's shield tiling with vertex set $\scrP_{\Gh,(1.7,0.6)}$ 
(top right panel),
and T\"ubingen triangle tilings with vertex sets
$\scrP_{\TT,(0.001,0)}$ (bottom left panel)
and $\scrP_{\TT,(1.3,0.4)}$ (bottom right panel).
The notation is explained in Section \ref{examplessec}.
In each patch, the origin is marked by a small black disc.}
\label{ABandGpatches}
\end{figure}

The main purpose of the present paper is to continue the explicit study
begun in \cite{Hammarhjelm2022}
of the limit distribution of normalized gaps between directions 
in the case of $\scrP$ belonging to certain 
families of planar cut-and-project sets,
which include some well-known classes of quasicrystals; see Figure \ref{ABandGpatches}.
Our focus will be on asymptotics for
\textit{large} gaps between directions.
The present work grew out of the initial study
carried out by Hammarhjelm
\cite{Hammarhjelm2021}. 

\vspace{4pt}

We remark that the methods developed in the present paper 
can be expected to be useful also for
questions related to the \textit{Lorentz gas}
on a cut-and-project  
scatterer configuration
-- namely, 
for the task of obtaining asymptotic estimates
for the transition probabilities
in the 
transport (Markov) process which arises
when considering such a Lorentz gas
in the limit of low scatterer density
\cite{MarklofStrombergsson2014},
\cite{MarklofStrombergsson2024}.
For the case of a lattice scatterer configuration,
such asymptotic estimates were obtained in 
\cite{MarklofStrombergsson2011},
and found an important application in
\cite{MarklofToth2016}.

\newpage
\subsection{The limit distribution of gaps for cut-and-project sets} 
To recall the precise definition of cut-and-project sets considered in
\cite{MarklofStrombergsson2015}, 
let $d,m\geq1$, set $n=d+m$,
and let 
$\scrL$ be a lattice of full rank in $\R^n=\R^d\times\R^m$. 
We refer to $\R^d$ and $\R^m$ as the \textit{physical space} and 
\textit{internal space,} respectively.
In the present paper the dimension of the physical space will always be $d=2$.
We write $\pi$ and $\pi_{\intl}$ for the orthogonal projection of
$\R^n$ 
onto the first $d$ coordinates and last $m$ coordinates\label{p:projections}.
Let $\scrA$ be the closure of $\pi_{\intl}(\scrL)$;\label{scrADEF}
this is a closed abelian subgroup of $\R^m$.
Let $\scrA^\circ$ be the identity component of $\scrA$;\label{scrAcircDEF}
this is a linear subspace of $\R^m$,
and set $m'=\dim\scrA=\dim\scrA^\circ$.
Let $\scrW$, the \textit{window}, 
be a bounded subset of $\scrA$ with nonempty interior
(with respect to the topology of $\scrA$)\label{p:window}. 
We will always assume that $\scrW$ and $\scrL$ are such that
\begin{align}\label{LWproperty}
\forall\, \vecy_1,\vecy_2\in\scrL:\quad
\text{if $\pi_{\intl}(\vecy_1)\in\scrW$ and $\pi_{\intl}(\vecy_2)\in\scrW$
and $\pi(\vecy_1)=\pi(\vecy_2)$, then $\vecy_1=\vecy_2$.}
\end{align}
We now define the \textit{cut-and-project set}
associated to $\scrW$ and $\scrL$ to be
\begin{align}\label{scrPWLdef}
\scrP=\scrP(\scrW,\scrL)=\bigl\{\pi(\vecy)\col\vecy\in\scrL,\:\pi_{\intl}(\vecy)\in\scrW\bigr\}\subset\R^d.
\end{align}
We will always assume that $\scrP(\scrW,\scrL)$ is
\textit{regular}, meaning that 
$\partial\scrW$ has measure zero with respect to the Haar measure of $\scrA$.
Under these assumptions,
it is known that $\scrP$ has an asymptotic density,
i.e.\ \eqref{asymptdensityDEF} holds,
with $c_\scrP$ being given by a simple explicit expression
in terms of $\scrL$ and $\scrW$ \cite[Prop.\ 3.2]{MarklofStrombergsson2014};
see also \eqref{cPformula} below.

Let us view the unit circle $\S_1^1$ as the set of $z\in\C$ with $|z|=1$, 
and for any $z\in\S_1^1$ let us call the number
$\frac1{2\pi}\arg(z)$ the \textit{normalized angle} of  $z$;
this gives an identification between $\S_1^1$ and $\R/\Z$.
Let us order the normalized angles 
of the points in $\Delta_R$ 
in an increasing list as
\begin{align}\label{xilist}
-\tfrac12<\xi_{R,1}\leq\xi_{R,2}\leq\cdots\leq\xi_{R,N(R)}\leq\tfrac12,
\end{align}
where $N(R)=\#\Delta_R$.
Also set $\xi_{R,0}=\xi_{R,N(R)}-1$.

The following result was proved in \cite{MarklofStrombergsson2015} (see also Remark \ref{gaplimitdistrexistsTHMREM} below).
\begin{thm}\cite[{Cor.\ 3}]{MarklofStrombergsson2015}
\label{gaplimitdistrexistsTHM}
Let $\scrP$ be a regular cut-and-project set in $\R^2$.
Then there exists a decreasing function $F:\R_{\geq0}\to[0,1]$,
which is continuous on $\R_{>0}$, such that for every $s\geq0$,
\begin{align}\label{gaplimitdistrexistsTHMres}
\lim_{R\to\infty}\frac{\#\{1\leq j\leq N(R)\col N(R)(\xi_{R,j}-\xi_{R,j-1})\geq s\}}{N(R)}=F(s).
\end{align}
Furthermore,
this limit relation \eqref{gaplimitdistrexistsTHMres} remains true, with the function $F$ unchanged,
if $\scrP$ is replaced by $\scrP\,T=\{\vecv T\col \vecv\in\scrP\}$ for any fixed $T\in\GL_2(\R)$.
\end{thm}

We recall in Section \ref{FsexplSEC} below some of the key steps in the proof of the above result;
a crucial ingredient 
is Ratner's classification \cite{Ratner1991} of measures invariant under unipotent flows 
in homogeneous dynamics.
The limit distribution function $F(s)$ in Theorem \ref{gaplimitdistrexistsTHM} is explicitly given by a formula involving
the Haar measure of a certain family of subsets of a homogeneous space. 
For the particular 
planar cut-and-project sets which we consider in the present paper, this homogeneous space turns out to be of the form
$\Gamma\bs\mathrm{H}$ where $\Gamma$ is a Hilbert modular subgroup of $\mathrm{H}=\SL_2(\R)\times\SL_2(\R)$.

\begin{figure}[h!]
\begin{center}
\framebox{
\begin{minipage}{0.28\textwidth}
\unitlength0.1\textwidth
\begin{picture}(10,9.5)(0,1.2)
\put(0.5,0.8){\includegraphics[width=0.9\textwidth]{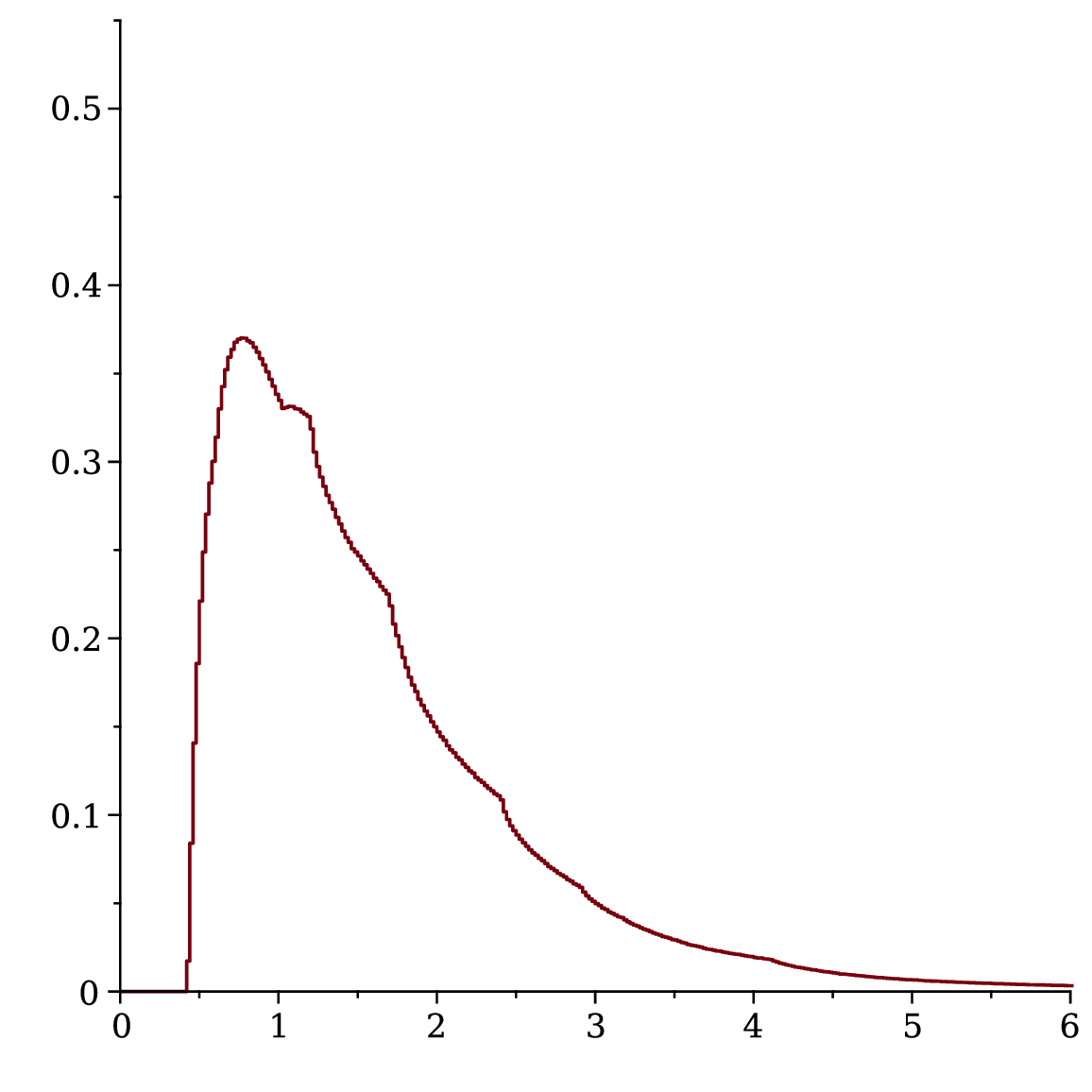}}
\put(0.1,9.7){\begin{footnotesize}\textbf{I.}\end{footnotesize}}
\end{picture}
\end{minipage}
}
\framebox{
\begin{minipage}{0.28\textwidth}
\unitlength0.1\textwidth
\begin{picture}(10,9.5)(0,1.2)
\put(0.0,9.9){\begin{footnotesize}\textbf{II.}\end{footnotesize}}
\put(0.5,0.8){\includegraphics[width=0.9\textwidth]{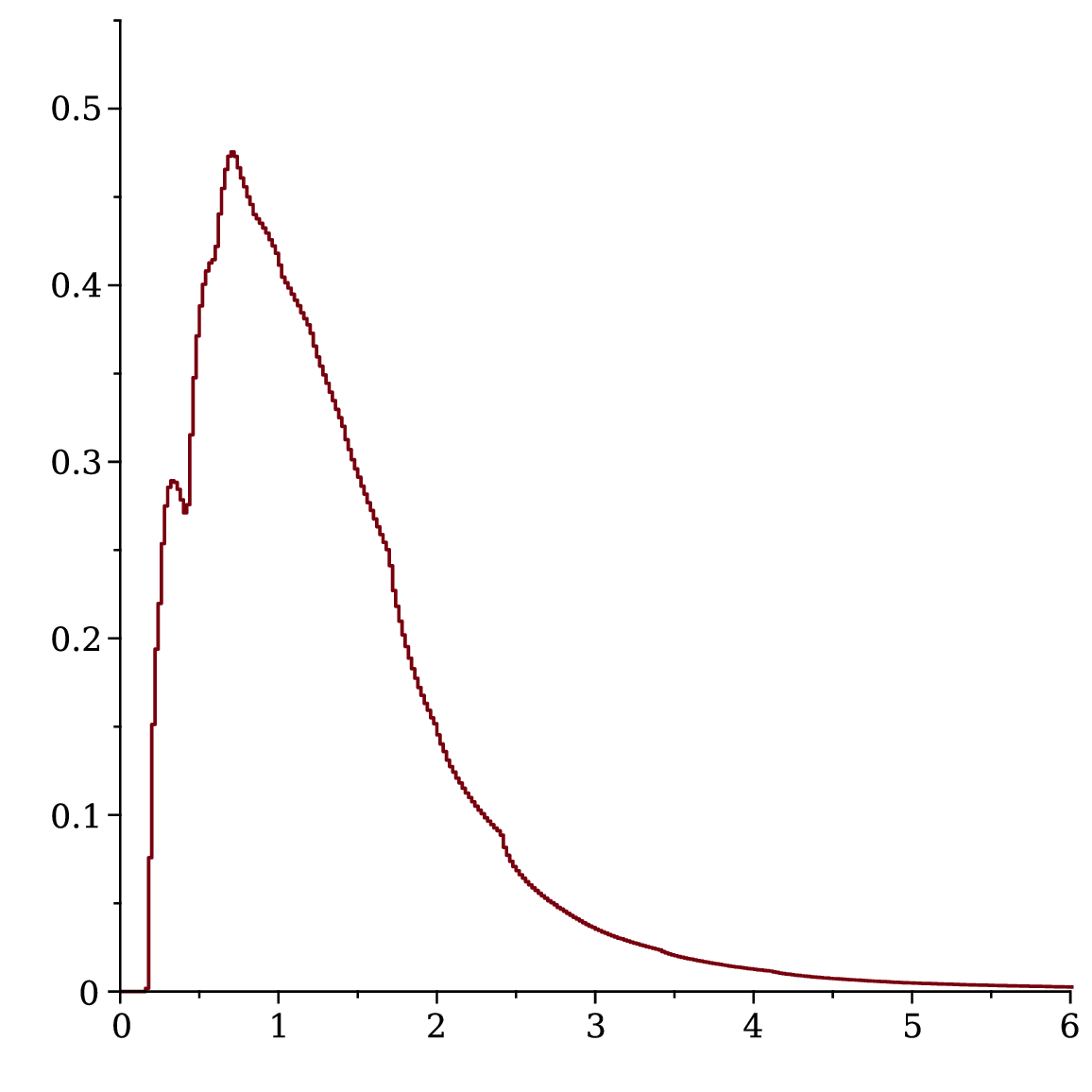}}
\end{picture}
\end{minipage}
}
\framebox{
\begin{minipage}{0.28\textwidth}
\unitlength0.1\textwidth
\begin{picture}(10,9.5)(0,1.2)
\put(0.0,9.9){\begin{footnotesize}\textbf{III.}\end{footnotesize}}
\put(0.5,0.8){\includegraphics[width=0.9\textwidth]{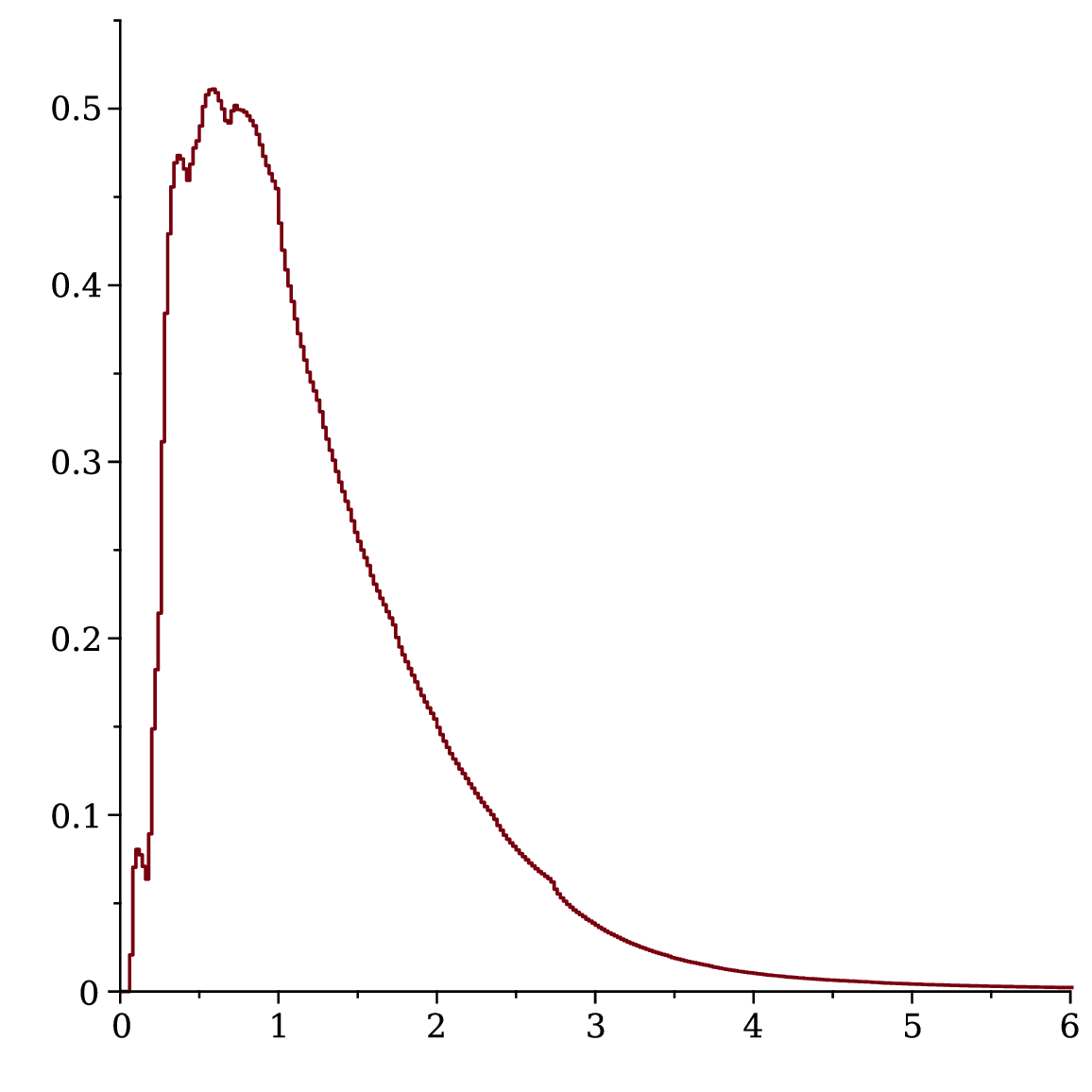}}
\end{picture}
\end{minipage}
}
\framebox{
\begin{minipage}{0.28\textwidth}
\unitlength0.1\textwidth
\begin{picture}(10,9.5)(0,1.2)
\put(0.0,9.9){\begin{footnotesize}\textbf{IV.}\end{footnotesize}}
\put(0.5,0.8){\includegraphics[width=0.9\textwidth]{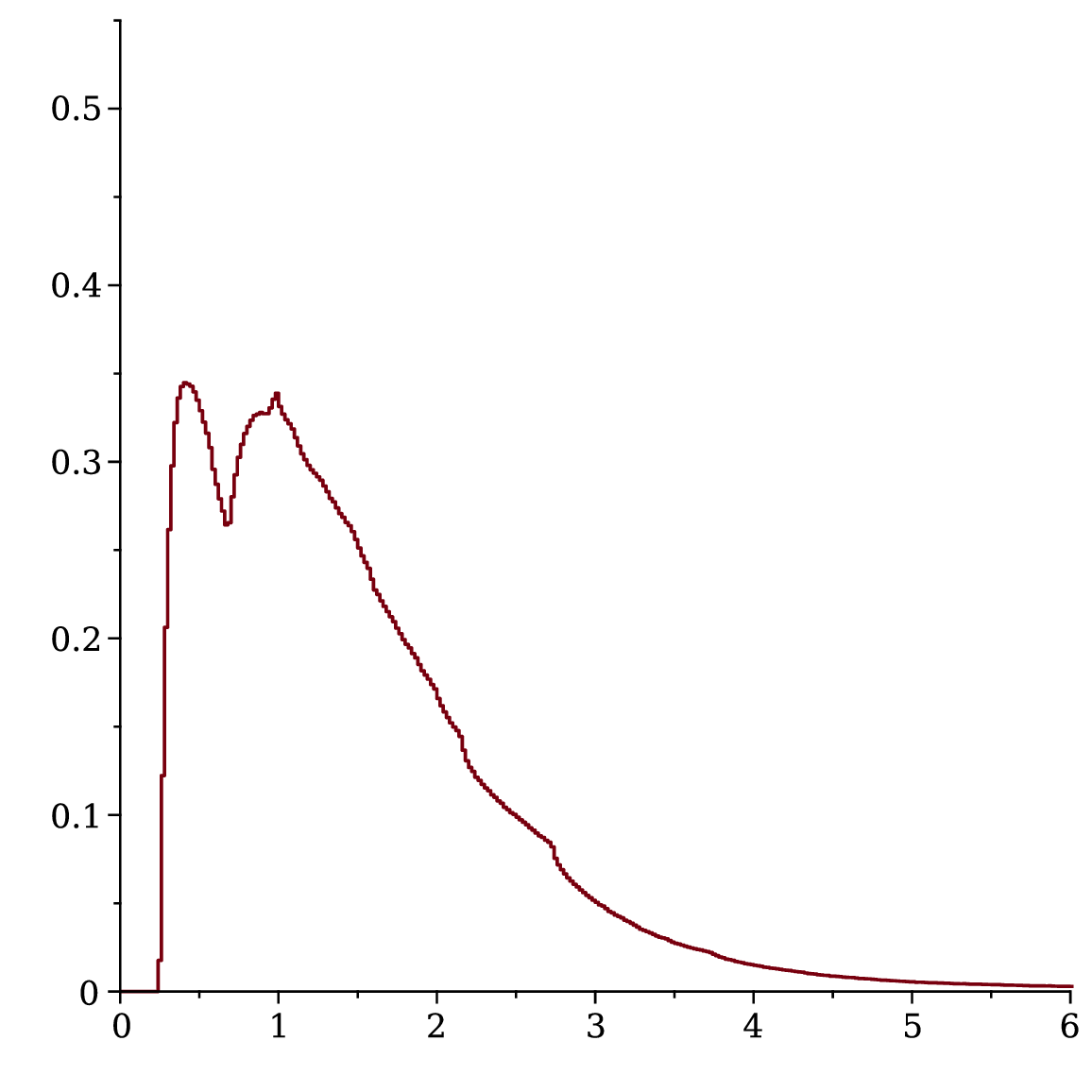}}
\end{picture}
\end{minipage}
}
\framebox{
\begin{minipage}{0.28\textwidth}
\unitlength0.1\textwidth
\begin{picture}(10,9.5)(0,1.2)
\put(0.0,9.9){\begin{footnotesize}\textbf{V.}\end{footnotesize}}
\put(0.5,0.8){\includegraphics[width=0.9\textwidth]{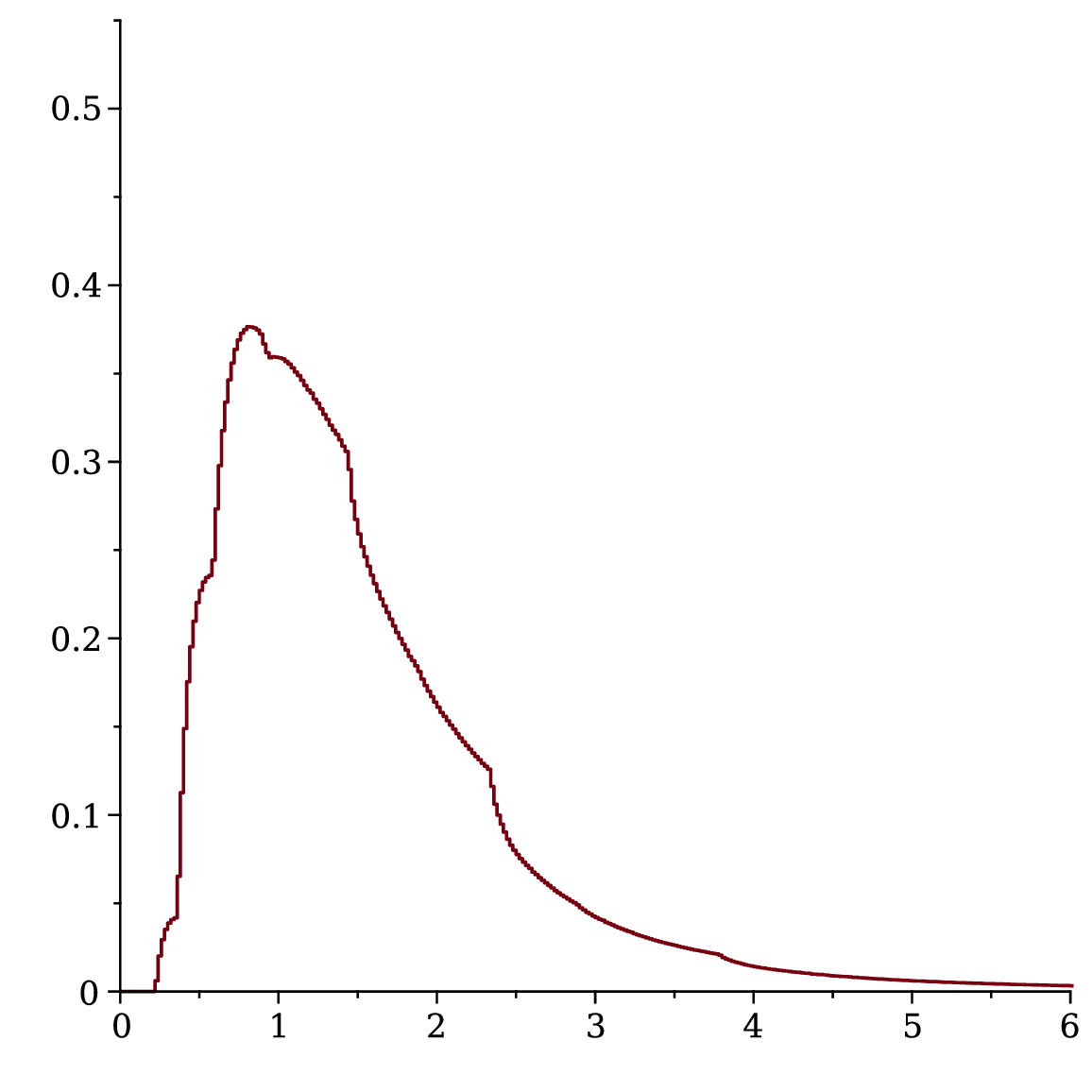}}
\end{picture}
\end{minipage}
}
\framebox{
\begin{minipage}{0.28\textwidth}
\unitlength0.1\textwidth
\begin{picture}(10,9.5)(0,1.2)
\put(0.0,9.9){\begin{footnotesize}\textbf{VI.}\end{footnotesize}}
\put(0.5,0.8){\includegraphics[width=0.9\textwidth]{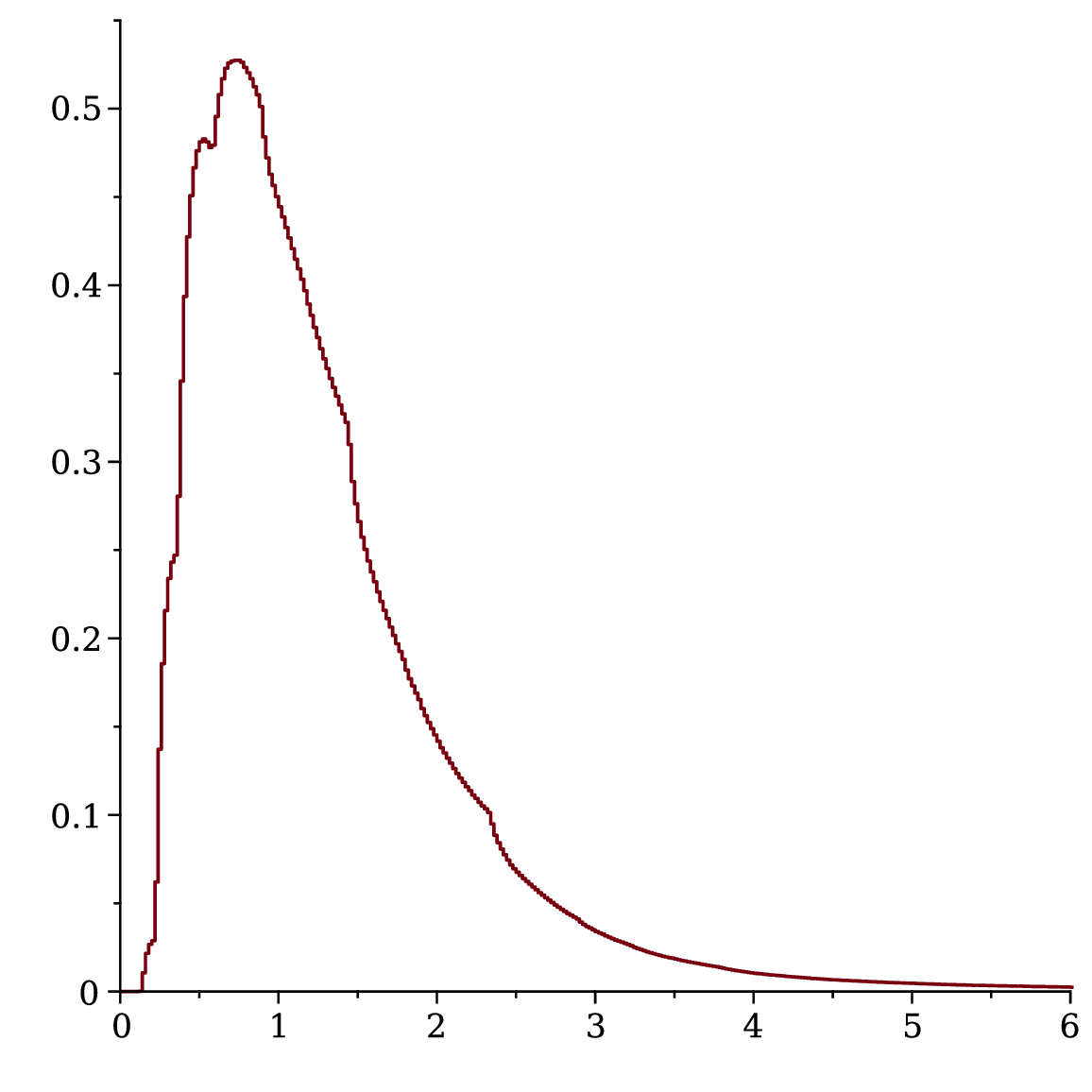}}
\end{picture}
\end{minipage}
}
\end{center}
\caption{Experimental graphs of the limiting density of normalized gaps, $-F'(s)$,
for the vertex sets of Ammann-Beenker tilings
$\scrP_{\AB,\vecw}$ with (I) $\vecw=\bn$ and (II) $\vecw=(0.9,0.3)$;
G\"ahler's shield tilings $\scrP_{\Gh,\vecw}$ with
(III) $\vecw=(1.7,0.6)$ and (IV) $\vecw=(0.1,0)$;
and T\"ubingen triangle tilings $\scrP_{\TT,\vecw}$ with 
(V) $\vecw=(0.5,0.4)$ and (VI) $\vecw=(1.3,0.4)$
(notation as in Section \ref{examplessec}).
These graphs were obtained from numerical computations described in
Section \ref{numcompSUBSEC}.
}
\label{densityplot1}
\end{figure}

Figure \ref{densityplot1}
shows conjectural graphs of the function $-F'(s)$,
i.e.\ the limiting density of normalized gaps,
for some examples of 
vertex sets of
Ammann-Beenker, G\"ahler's shield, and T\"ubingen triangle tilings.
\begin{remark}\label{gaplimitdistrexistsTHMREM}
Theorem \ref{gaplimitdistrexistsTHM} is a special case of
\cite[{Cor.\ 3}]{MarklofStrombergsson2015},
since in 
\cite[{Cor.\ 3}]{MarklofStrombergsson2015}
we allow $\scrL$ in $\scrP=\scrP(\scrW,\scrL)$
to be a \textit{translate} of a lattice in $\R^n$,
whereas in the present paper we always require $\scrL$ to be a genuine lattice,
i.e.\ $\bn\in\scrL$.
In fact, in Theorem \ref{MAINTHM1} below
we will also assume that $\bn\in\scrW$,  
and so $\bn\in\scrP$.

We also point out that 
the proof of \cite[{Cor.\ 3}]{MarklofStrombergsson2015}
immediately extends (by utilizing the freedom of choice of the measure ``$\lambda$'' in 
\cite[{Thm.\ 2}]{MarklofStrombergsson2015})
to show that the limiting gap distribution for $\scrP$ remains the
same if we restrict attention to the directions lying in any fixed subinterval of 
$\S_1^1$.
That is, the following
more general version of \eqref{gaplimitdistrexistsTHMres} holds:
For any fixed $-\frac12\leq \alpha_1<\alpha_2\leq\frac12$,
\begin{align}\label{gaplimitdistrexistsTHMresGEN}
\lim_{R\to\infty}\frac{\#\{1\leq j\leq N(R)\col \xi_{R,j}\in (\alpha_1,\alpha_2],\: N(R)(\xi_{R,j}-\xi_{R,j-1})\geq s\}}
{(\alpha_2-\alpha_1)\cdot N(R)}=F(s).
\end{align}

The last statement of Theorem \ref{gaplimitdistrexistsTHM}, 
regarding the invariance of $F(s)$ when replacing $\scrP$ by $\scrP\, T$,
is easily derived from 
\eqref{gaplimitdistrexistsTHMresGEN} by a limit argument letting $\alpha_2-\alpha_1\to0$.
(In the special cases when $T$ is simply a scaling by a constant or a rotation or a reflection,
the statement even follows directly from \eqref{gaplimitdistrexistsTHMres}.)
\end{remark}

\begin{remark}\label{visibleREM} 
Our notation differs from 
the notation used in \cite{mBfGcHtJ2014}, \cite{MarklofStrombergsson2015} and \cite{Hammarhjelm2022}:
In the present paper we do not remove the points in $\scrP$ which
are 'invisible' from the origin; 
thus we allow $\Delta_R$ to be a multi-set,
and we allow  
equalities in the list in \eqref{xilist}.
Also, our function $F(s)$ 
equals the function in the right hand side of 
\cite[(1.15)]{MarklofStrombergsson2015},
and is thus in general not the same as 
``$F(s)$'' in 
\cite[{Cor.\ 3}]{MarklofStrombergsson2015}.
Clearly by \eqref{gaplimitdistrexistsTHMres}, our function $F(s)$
satisfies $F(0)=1$.
In general $F$ has a jump discontinuity at $0$;
it follows from 
\cite[{Cor.\ 3}]{MarklofStrombergsson2015}
that $\lim_{s\to0+}F(s)=\kappa_{\scrP}$
where $\kappa_{\scrP}$ is the relative density of visible points in $\scrP$,
a quantity which was defined and proved to exist in\label{kappaPDEF}
\cite{MarklofStrombergsson2015}.
Note that it is immediate to translate between  
the two versions of ``$F(s)$'',
except that it requires knowledge of the constant
$\kappa_{\scrP}$, which is non-trivial to compute.
See \cite{Hammarhjelm2022} 
for explicit formulas for $\kappa_{\scrP}$ in special cases.
See also the discussion at the end of Section \ref{numcompSUBSEC} below
regarding the effect of the different normalizations when 
comparing the graph in the top left panel of Figure \ref{densityplot1}
with the graphs in \cite[Fig.~9]{mBfGcHtJ2014}
and
\cite[Fig.\ 2]{Hammarhjelm2022}.
\end{remark}

\subsection{Main result; tail asymptotics}
Our main goal in the present paper is to describe the tail asymptotics
of the distribution function  $F(s)$ in Theorem \ref{gaplimitdistrexistsTHM},
for a particular class of cut-and-project sets in $\R^2$,
which includes several classical examples.
We expect that the methods which we develop can be extended to more general cases as well
(in particular see Remark \ref{PenroseREM} below).
We now give the description of the class which we will consider.
Let $K$ be a real quadratic field;\label{p:qufield}
let $\scrO_K$ be its ring of integers;
let $\sigma$ be the unique non-trivial automorphism of $K$,
and let $\lambda\in \cO_K^{\times}$ be the fundamental unit\label{p:unitgp}\label{p:lambda}
(thus $\lambda>1$ and $\cO_K^{\times}=\left\{\pm \lambda^n\col n\in\Z\right\}$).
Let $\scrL_K$ be the Minkowski embedding of $\scrO_K^2$ in $\R^4$, 
viz.,
\begin{align}\label{scrLKdef}
\scrL_K:=\left\{(\alpha,\beta, \sigma(\alpha),\sigma(\beta))\in\R^4\col (\alpha,\beta)\in \cO_K^2\right\}.
\end{align}
We set $d=m=2$, that is,
we view $\R^4$ as the product of a 2-dimensional physical space and a 2-dimensional internal
space,
and $\pi$ and $\pi_{\intl}$ denote the orthogonal projection of $\R^4$
onto the first 2 coordinates and the last 2 coordinates, respectively.
Note that in this case $\pi_{\intl}(\scrL_K)$ is dense in $\R^2$, i.e.\ we have $\scrA=\R^2$,
and so we take the window $\scrW$ to be a bounded subset of $\R^2$ with nonempty interior.
(In the present case
the restriction of $\pi$ to $\scrL_K$ is injective,
and so the property \eqref{LWproperty} is automatically fulfilled.)
We will always assume that $\scrP:=\scrP(\scrW,\scrL_K)$ is regular,
i.e.\ that $\partial\scrW$ has Lebesgue measure zero.
\vspace{2pt}

For this class of cut-and-project sets,  
the formula for the asymptotic density of $\scrP$,
\cite[(1.7)]{MarklofStrombergsson2014}, 
becomes (see also Section \ref{HilbertprelSEC} below):
\begin{align}\label{cPformula}
c_{\cP}=\Delta_K^{-1}\Area(\scrW),
\end{align}
where $\Delta_K$ is the discriminant of $K$.
\vspace{2pt}

As we prove in Remark \ref{rmk:nullsetinva} below, 
the limiting gap distribution function is unchanged if $\cW$ is modified by any measure zero set.
In particular, without loss of generality we may assume $\cW$ is open (indeed, otherwise replace $\cW$ by its interior).
Finally, in the present paper we will make the key assumption that 
$\bn\in\scrW$.
\vspace{2pt}

In our approach we actually work with the integral of the distribution function $F$,
i.e.\ with
\begin{align}\label{Gdef}
G(s):=\int_s^\infty F(t)\,\text{d}t.
\end{align}
As we will note below (see Lemma \ref{prop:gtof}),
using the fact that $F$ is decreasing,
an easy interpolation argument allows us 
to deduce an asymptotic formula for $F$
once we know an asymptotic formula for $G$.
\vspace{2pt}

Our main result is the following:

\begin{thm}\label{MAINTHM1}
Let $\scrP=\scrP(\scrW,\scrL_K)$
where $K$ is a real quadratic field
and $\scrW$ is a bounded open subset of $\R^2$ such that
$\bn\in\scrW$ and $\partial\scrW$ has Lebesgue measure zero.
Let $F(s)$ be the associated limiting gap distribution function
as in Theorem \ref{gaplimitdistrexistsTHM}, and let $G(s)=\int_s^\infty F(t)\,\textup{d}t$.
Then there exists a positive constant $a_{\cP}$ such that 
\begin{align}\label{equ:gsestgenwin}
G(s)\sim a_{\scrP}s^{-1}\quad
\text{and} \quad F(s)\sim a_{\scrP}s^{-2}\qquad \text{as }\: s\to\infty.
\end{align}
If we further assume $\cW$ to be convex, then we have
\begin{align}\label{equ:polyest}
G(s)=a_{\scrP}s^{-1}+O\bigl(s^{-2}\log s\bigr)
\quad\text{and}\quad
F(s)=a_{\scrP}s^{-2}+O\Bigl(s^{-\frac52}\sqrt{\log s}\Bigr)
\qquad\text{as }\: s\to\infty.
\end{align}
\end{thm}

\begin{remark}\label{rmk:classonefor}
As mentioned above, in order to prove Theorem \ref{MAINTHM1},
it suffices to prove the two asymptotic estimates for $G(s)$ in \eqref{equ:gsestgenwin} and \eqref{equ:polyest}. 
Both these estimates will in fact follow from a more precise,
general asymptotic formula for $G(s)$ 
with an explicit error term; see \eqref{MAINTHM1res1} below. 
We also mention that our analysis leads to an explicit formula 
for the leading coefficient $a_{\cP}$; see \eqref{MAINTHM1res2}. 
\end{remark}

\subsection{Formula for $a_{\cP}$} 
\label{sec:forap}
As we will now describe, after imposing one more assumption on the window set,
we can further evaluate the aforementioned formula for $a_{\cP}$.
To state our result, we need to introduce some more notation. 
For $\theta\in\R$ we write $\text{k}_{\theta}:=\left(\begin{smallmatrix}
\cos\theta & -\sin\theta\\
\sin\theta & \cos\theta\end{smallmatrix}\right)\in\SL_2(\R)$. Then for any subset $\cW\subset \R^2$ we
let $\ell_{\cW}(\theta)$ be the projection of $\cW \text{k}_{-\theta}$ on the $x$-axis. 
On top of the assumptions in Theorem \ref{MAINTHM1},
we will now assume that for each $\theta$, $\ell_{\cW}(\theta)$ is an \textit{interval.}
Note that this assumption is always satisfied if $\scrW$ is \textit{connected},
but it also holds for many non-connected window sets $\scrW$.

Note that since $\cW$ is open and $\bn\in\cW$, 
$\ell_{\cW}(\theta)$ is an open subset of $\R$ containing $0$.
Hence if $\ell_{\cW}(\theta)$ is assumed to be an interval, then we can parametrize it 
by two positive numbers $r(\theta), \nu(\theta)$ via 
\begin{align}\label{equ:ellw}
\ell_{\cW}(\theta)=r(\theta)(-\nu(\theta), 1).
\end{align} 

\begin{thm}\label{thm:moreexforclassone}
Retain the notation and assumptions in Theorem \ref{MAINTHM1} and further assume that $K$ is of class number one, and 
that $\ell_{\cW}(\theta)$ is an interval for each $\theta$. Then there exist a finite partition $\R_{>0}=\bigsqcup_{j=1}^{\, l}S_{j}$ of $\R_{>0}$
into intervals, and non-negative constants $A_{j}, B_{j}$ $(1\leq j\leq l)$ depending only on $K$, such that
\begin{align}\label{equ:genforapclassno1}
a_{\scrP}=\frac{\Area(\scrW)}{4 \Delta_K^2\zeta_K(2)}\sum_{j=1}^{l}\int_{\tilde{S}_{j}}r(\theta)^{-2}\left(A_{j}+B_{j}\nu(\theta)^{-2}\right)\,\mathrm{d}\theta,
\end{align}
where  
$\zeta_K$ is the Dedekind zeta function attached to $K$, 
$\tilde{S}_{j}:=\{\theta\in [0,2\pi)\col \nu(\theta)\in S_{j}\}$,
and $r(\theta)$ and $\nu(\theta)$ are defined by \eqref{equ:ellw}. 
\end{thm}
\begin{remark}
We stress that the intervals $S_{j}$ are allowed to be open, closed, or half-open,
and may be degenerate, i.e.\ of the form $[a,a]=\{a\}$ for some $a>0$.
\end{remark}
\begin{remark}\label{AjBjcomputableREM}
Our analysis applies for any real quadratic field $K$ and the class number one assumption in 
Theorem \ref{thm:moreexforclassone}
is only for simplicity of presentation; see Theorem \ref{thm:moreexfor} in Section \ref{acPmoreexplicitSEC}
for the most general version. 
We also mention that the partition $\R_{>0}=\bigsqcup_{j=1}^{l}S_{j}$ and the constants $A_j, B_j$ in the above theorem are all computable and we will illustrate it in the next section for three well-known classes of quasicrystals, namely, the vertex sets of the
Ammann-Beenker, G\"ahler's shield, and T\"ubingen triangle tilings.
\end{remark}

\begin{remark}\label{aPsimpleformulaREM}
Among the sets $\tilde{S}_j$, there is one that is always non-empty, namely, the unique $\tilde{S}_i$ such that $1\in S_i$.
In certain cases, this $\tilde{S}_i$ equals the whole interval $[0,2\pi)$
while all other $\tilde{S}_j$'s are empty. 
In this case, the formula for $a_{\mathcal{P}}$ can be simplified into \begin{align*}
a_{\cP}=C_K\Area(\cW)\Area(\cW^*),
\end{align*}
where $C_K$ is some positive constant depending only on $K$, and 
$\scrW^*$ is the \textit{polar set} of $\cW$, i.e.\
\begin{align}\label{def:polbody}
\scrW^*:=\{\vecz\in\R^2\col \vecz\cdot\vecw\leq1\:\:\forall\:\vecw\in\scrW\}.
\end{align}
Here $\cdot$ is the standard scalar product in $\R^2$. See Section \ref{UNDintegralSEC} for more details. 
If $\cW$ is also convex and centrally symmetric, the product 
$\Area(\cW)\Area(\cW^*)$ appearing in the above formula is known as the \textit{Mahler volume} of $\cW$.
\end{remark}

\begin{remark}
For any $\cW$ as in Theorem \ref{thm:moreexforclassone},
replacing $\cW$ by its \textit{convex hull} does not 
affect the intervals $\ell_{\cW}(\theta)$,
hence does not affect the functions $r(\theta)$ and $\nu(\theta)$.
Therefore, in the explicit formula \eqref{equ:genforapclassno1},
only the factor $\Area(\cW)$ is affected when replacing 
$\cW$ by its convex hull.
In this connection it should also be noted that 
the polar set of the convex hull of $\cW$ equals the polar set of $\cW$.
\end{remark}

\subsection{Examples}
\label{examplessec}
We next illustrate 
how our results apply 
in three cases of 
well-known 
planar quasicrystals. 
The three propositions below are all proved in Section \ref{sec:exmapksqrt23}.
More details, and comparison with numerics, is provided in Section \ref{numcompSEC}.

We first consider the \textit{Ammann-Beenker tiling},
which was discovered by
Robert Ammann in the 70s
and first described in
\cite{GS87} and \cite{AGS92}.
Specifically, we consider the ``A5 set, variant (b)'', in the notation of
\cite{AGS92}; 
see Figure \ref{ABandGpatches} (top left panel) above for a small patch of this tiling.
It is well-known that the set of vertices of an 
arbitrary Ammann-Beenker tiling
can be generated using the cut-and-project construction.
Specifically, 
let $K=\Q(\sqrt2)$, and let $\scrW_{\AB}\subset\R^2$ be the open regular octagon 
centered at the origin of edge length $1$, 
oriented so that four of the edges are perpendicular to a coordinate axis\label{p:waboct}.
For each $\vecw\in\R^2$ we set
\begin{align}\label{PABwDEF}
\scrW_{\vecw}^{(\AB)}:=(\scrW_{\AB}+\vecw)g_1\subset\R^2
\qquad\text{and}\qquad
\scrP_{\AB,\vecw}:=\scrP(\scrW^{(\AB)}_{\vecw},\scrL_K)g_2,
\end{align}
where $g_1:=\smatr101{\sqrt2}$
and $g_2:=\smatr10{1/\sqrt2\:}{1/\sqrt2}$.
Then for any $\vecw\in\scrW_{\AB}$ with the property that
$\pi_{\intl}(\scrL_K)\cap\partial\scrW^{(\AB)}_{\vecw}=\emptyset$,
the cut-and-project set
$\scrP_{\AB,\vecw}$
is the vertex set of an Ammann-Beenker tiling with one vertex at $\bn$,
and conversely,
the vertex set of \textit{any} Ammann-Beenker tiling having a vertex at $\bn$
is (up to scaling and rotation)
either equal to such a point set $\scrP_{\AB,\vecw}$,
or can be obtained as a limit of such point sets $\scrP_{\AB,\vecw}$
in an appropriate topology. 
(See \cite[Ch.\ 7.3]{mBuG2013} and Section \ref{ABexplicitSEC} below.)

Note that because of
\eqref{PABwDEF} and the last statement in Theorem \ref{gaplimitdistrexistsTHM},
the formula \eqref{equ:polyest} in Theorem~\ref{MAINTHM1}
applies to the cut-and-project set
$\scrP_{\AB,\vecw}$ for any $\vecw\in\scrW_{\AB}$.
Here 
the constant $a_{\scrP}$ is given by a quite simple formula,
also valid for much more general window sets:
\begin{Prop}\label{ABaPexpl}
Let $\scrP=\scrP(\scrW,\scrL_K)$
with $K=\Q(\sqrt2)$ and with $\scrW$ as in Theorem \ref{thm:moreexforclassone},
and let $\ell_{\scrW}(\theta)$ be parametrized as in \eqref{equ:ellw}.
Then 
\begin{align}\notag
a_{\cP}=\frac{3\sqrt{2}}{16\pi^4}\Area(\cW)\biggl((7+8\sqrt{2})\int_{\tilde{S}_1}r(\theta)^{-2}\,\mathrm{d}\theta+(4+6\sqrt{2})\int_{\tilde{S}_2}r(\theta)^{-2}\,\mathrm{d}\theta
\hspace{100pt}
\\ \label{ABaPexplres1}
+(2+2\sqrt{2})\int_{\tilde{S}_3}r(\theta)^{-2}\,\mathrm{d}\theta +\int_{\tilde{S}_4}r(\theta)^{-2}\,\mathrm{d}\theta\biggr),
\end{align}
where $\tilde{S}_j:=\{\theta\in [0,2\pi)\col \nu(\theta)\in S_j\}$, with
\begin{align*}
S_{1}=\bigl(0,\sqrt{2}-1\bigr], 
\quad S_{2}=\bigl(\sqrt{2}-1,\tfrac{1}{2}\sqrt{2}\bigr),
\quad S_{3}=\bigl[\tfrac12\sqrt{2}, \sqrt{2}\bigr]
\quad\text{and}\quad
 S_{4}=\bigl(\sqrt{2}, 1+\sqrt{2}\bigr).
\end{align*}
When $\tilde{S}_3=[0,2\pi)$, the above formula simplifies to
\begin{align}\label{ABaPexplres2}
a_{\cP}&=\frac{3\sqrt{2}(1+\sqrt{2})}{4\pi^4}\Area(\scrW)\Area(\scrW^*).
\end{align}
\end{Prop}

Note that Proposition \ref{ABaPexpl} applies 
when $\scrW=\scrW^{(\AB)}_{\vecw}$ for any $\vecw\in\scrW_{\AB}$,
and the formula \eqref{ABaPexplres2} holds whenever $\vecw$ lies sufficiently near $\bn$.
In particular, \eqref{ABaPexplres2} implies that
\begin{align} 
a_{\scrP_{\AB,\bn}}=24/\pi^4=0.24638\ldots.
\end{align}
In Section \ref{numcompSUBSEC}, 
we present a comparison,
for a few examples of points $\vecw\in\scrW_{\AB}$,
between the exact values of $a_{\scrP}$ computed using 
Proposition \ref{ABaPexpl},
and numerically computed approximate values of $s^2F(s)$ for large $s$.

\vspace{5pt}

Next we consider \textit{G\"ahler's shield tiling,}
which was discovered in \cite{Gah88};
these tilings are built up of a certain \label{Ghstart}
triple of tiles (a triangle, a square, and a hexagon called a 'shield'),
equipped with local matching rules.
The vertex set of a    
G\"ahler's shield tiling
can be obtained using the cut-and-project construction in the following way:
Let $K=\Q(\sqrt3)$, and let $\scrW_{\Gh}\subset\R^2$
be the open regular dodecagon 
centered at the origin of edge length $1$,
so that four of the edges are perpendicular to a coordinate axis\label{p:vecwgh}.
For each $\vecw\in\R^2$ we set
\begin{align}\label{PGhwDEF}
\scrW_{\vecw}^{(\Gh)}:=(\scrW_{\Gh}+\vecw)g_1\subset\R^2
\qquad\text{and}\qquad
\scrP_{\Gh,\vecw}:=\scrP(\scrW_{\vecw}^{(\Gh)},\scrL_K)g_2,
\end{align}
where $g_1:=\smatr10{\sqrt3}2$
and $g_2:=\smatr10{\sqrt3/2}{1/2}$.
Then for any $\vecw\in\scrW_{\Gh}$ with the property that
$\pi_{\intl}(\scrL_K)\cap\partial\scrW^{(\Gh)}_{\vecw}=\emptyset$,
the cut-and-project set
$\scrP_{\Gh,\vecw}$
is the vertex set of a G\"ahler's shield tiling with one vertex at $\bn$.
(See \cite[Ch.\ 7.3]{mBuG2013} and Section \ref{GSexplicitSEC} below.)

Again, \eqref{equ:polyest} in Theorem \ref{MAINTHM1}
applies to 
$\scrP_{\Gh,\vecw}$ for any $\vecw\in\scrW_{\Gh}$,
and regarding  
$a_{\scrP}$ we have: 
\begin{Prop}\label{GHaPexpl}
Let $\scrP=\scrP(\scrW,\scrL_K)$
with $K=\Q(\sqrt3)$ 
and with $\scrW$ as in Theorem \ref{thm:moreexforclassone},
and let $\ell_{\scrW}(\theta)$ be parametrized as in \eqref{equ:ellw}.
Then 
\begin{align}\notag
a_{\cP}=\frac{\sqrt{3}}{16\pi^4}\Area(\cW)\biggl((17+16\sqrt{3})\int_{\tilde{S}_1}r(\theta)^{-2}\,\mathrm{d}\theta
+(9+12\sqrt{3})\int_{\tilde{S}_2}r(\theta)^{-2}\,\mathrm{d}\theta
\hspace{120pt}
\\ \label{GHaPexplres1}
+(6+6\sqrt{3})\int_{\tilde{S}_3}r(\theta)^{-2}\,\mathrm{d}\theta 
+(3+4\sqrt3)\int_{\tilde{S}_4}r(\theta)^{-2}\,\mathrm{d}\theta
+(3+2\sqrt3)\int_{\tilde{S}_5}r(\theta)^{-2}\,\mathrm{d}\theta
+2\int_{\tilde{S}_6}r(\theta)^{-2}\,\mathrm{d}\theta\biggr),
\end{align}
where $\tilde{S}_j:=\{\theta\in [0,2\pi)\col \nu(\theta)\in S_j\}$, with
$S_{1}=\bigl(0,\tfrac1{1+\sqrt3}\bigr]$,
$S_{2}=\bigl(\tfrac1{1+\sqrt3},\tfrac1{\sqrt3}\bigr)$,
$S_{3}=\bigl[\tfrac1{\sqrt3},\sqrt3-1\bigr]$,
$S_{4}=\bigl(\sqrt{3}-1,\tfrac1{\sqrt3-1}\bigr)$,
$S_{5}=\bigl[\tfrac1{\sqrt3-1},\sqrt3\bigr]$
and $S_{6}=\bigl(\sqrt{3}, 1+\sqrt{3}\bigr)$.
When $\tilde{S}_4=[0,2\pi)$, the above formula simplifies to
\begin{align}\label{GHaPexplres2}
a_{\cP}&=\frac{12+3\sqrt3}{8\pi^4}\Area(\scrW)\Area(\scrW^*).
\end{align}
\end{Prop}

Note that Proposition \ref{GHaPexpl} applies 
when $\scrW=\scrW_{\vecw}^{(\Gh)}$ for any $\vecw\in\scrW_{\Gh}$,
and the formula \eqref{GHaPexplres2} holds whenever $\vecw$ lies sufficiently near $\bn$.
Again see Section \ref{numcompSUBSEC} for numerical computations 
related to the values of $a_{\scrP}$ for $\scrP=\scrP_{\Gh,\vecw}$.

\vspace{5pt}

Finally, we consider the 
\textit{T\"ubingen triangle tiling},
which was discovered and studied in \cite{BKSZ90};
\label{TTTstart}
these tilings are built up of a certain pair of isosceles triangles.
The set of vertices of a
T\"ubingen triangle tiling
can be obtained using the cut-and-project construction in the following way.
Let $K=\Q(\sqrt5)$ and $\tau=\frac12(1+\sqrt5)$,
and let $\scrW_{\TT}\subset\R^2$
be the open regular decagon centered at the origin of edge length 
$\sqrt{(2+\tau)/5}$, oriented so that two of the edges are perpendicular
to the first coordinate axis.
For each $\vecw\in\R^2$ we set
\begin{align}\label{PTTTwDEF}
\scrW_{\vecw}^{(\TT)}:=(\scrW_{\TT}+\vecw)g_1\subset\R^2
\qquad\text{and}\qquad
\scrP_{\TT,\vecw}:=\scrP(\scrW^{(\TT)}_{\vecw},\scrL_K)g_2,
\end{align}
where 
$g_1:=\smatr{1}{0}0{\sqrt{(2+\tau)/5}}\smatr10{\tau}2$
and 
$g_2:=\smatr10{(\tau-1)/2\:}{\sqrt{2+\tau}/2}$.
Then for any $\vecw\in\scrW_{\TT}$ with the property that
$\frac15\pi_{\intl}(\scrL_K)\cap\partial\scrW^{(\TT)}_{\vecw}=\emptyset$,
the cut-and-project set
$\scrP_{\TT,\vecw}$
is the vertex set of a T\"ubingen triangle tiling with one vertex at $\bn$.
(See \cite[Ch.\ 7.3]{mBuG2013}, \cite[Sec.\ 4]{BKSZ90} and Section \ref{TTTexplicitSEC} below.)

Again, \eqref{equ:polyest} in Theorem \ref{MAINTHM1}
applies to 
$\scrP_{\TT,\vecw}$ for any $\vecw\in\scrW_{\TT}$,
and regarding $a_{\scrP}$ we have: 
\begin{Prop}\label{TTaPexpl}
Let $\scrP=\scrP(\scrW,\scrL_K)$
with $K=\Q(\sqrt5)$ 
and with $\scrW$ as in Theorem \ref{thm:moreexforclassone},
and let $\ell_{\scrW}(\theta)$ be parametrized as in \eqref{equ:ellw}.
Then 
\begin{align}\label{TTaPexplres1}
a_{\cP}=\frac{3\sqrt{5}}{16\pi^4}\Area(\cW)\biggl(&(5+3\sqrt{5})\int_{\tilde{S}_{1}}r(\theta)^{-2}\,\mathrm{d}\theta+(1+\sqrt{5})\int_{\tilde{S}_{2}}r(\theta)^{-2}\,\mathrm{d}\theta
\biggr),
\end{align}
where $\tilde{S}_j:=\{\theta\in [0,2\pi)\col \nu(\theta)\in S_j\}$, with
$S_1=(0,\frac{\sqrt{5}-1}{2}]$
and $S_2=(\frac{\sqrt{5}-1}{2},\frac{1+\sqrt{5}}{2})$. 
When $\tilde{S}_2=[0,2\pi)$, the above formula simplifies to
\begin{align}\label{TTaPexplres2}
a_{\cP}&=\frac{3(5+\sqrt5)}{8\pi^4}\Area(\scrW)\Area(\scrW^*).
\end{align}
\end{Prop}

Proposition \ref{TTaPexpl} applies 
when $\scrW=\scrW_{\vecw}^{(\TT)}$ for any $\vecw\in\scrW_{\TT}$,
and the formula \eqref{TTaPexplres2} holds whenever $\vecw$ lies sufficiently near $\bn$.
Again see Section \ref{numcompSUBSEC} for related numerical computations.

\begin{remark}\label{PenroseREM}
As a concluding remark of the introduction, we mention that our main result, Theorem \ref{MAINTHM1}, does not apply to the vertex sets 
of the classical rhombic Penrose tilings, 
as these can only be realized as \textit{unions} of (four) translates of the cut-and-project type sets considered in Theorem~\ref{MAINTHM1}; 
see \cite[Sec. 2.5]{MarklofStrombergsson2014}. 
Nevertheless, 
in preliminary work, using similar methods as in the present paper 
we have proved analogous tail asymptotic formulas for the limiting gap distribution function for point sets in this generality, 
thus in particular covering 
the rhombic Penrose tilings. 
\end{remark}

\section{Preliminaries}
\label{HilbertprelSEC}

\subsection{Real quadratic fields}
In this section we give a brief review on backgrounds on real quadratic fields.   
Let $K=\Q(\sqrt{d})$ be a real quadratic field \label{p:qufield_rep}
with $d\in \N$ square-free. Let $\cO_K=\Z[\tau]\subset K$ be its \textit{ring of integers}, where  \label{p:tau}
\begin{align*}
\tau=\left\lbrace\begin{array}{ll} 
\frac{1+\sqrt{d}}{2} & d\equiv 1\Mod{4},\\[4pt]
\sqrt{d}& d\equiv 2,3\Mod{4}.
\end{array}\right.
\end{align*}
Let $J_K$ be the (abelian) group of \textit{fractional ideals} of $\cO_K$, and let $P_K$ be the subgroup of \textit{principal ideals}. 
We let
$C_K=J_K/P_K$ be the \textit{ideal class group} of $K$, and we call elements in $C_K$ \textit{ideal classes} of $K$. \label{p:clgrp} 
For any nonzero $a_1,\ldots,a_m\in K$, we denote by $\langle a_1,\ldots, a_m\rangle$ the fractional ideal generated by $a_1,\ldots, a_m$.

Let $\sigma$ be the unique non-trivial automorphism of $K$ \label{p:sigma} and consider the embedding
\begin{align}\label{equ:embd}
\iota: K\to \R^2,\qquad k\mapsto (k,\sigma(k)).
\end{align}
We note that $\iota(I)$ is a lattice in $\R^2$ for any $I\in J_K$. In particular, $\iota(\cO_K)$ is the lattice generated by $(1,1)$ and $(\tau, \sigma(\tau))$, which is of covolume $\Delta_K^{1/2}$, where 
\begin{align}\label{discriminantKdef}
\Delta_K:=\left\lbrace\begin{array}{ll} 
d & d\equiv 1\Mod{4},\\
4d& d\equiv 2,3\Mod{4},
\end{array}\right.
\end{align}
is the \textit{discriminant} of $K$. 
For any fractional ideal $I\in J_K$, 
its \textit{absolute norm} is defined by 
\begin{align}\label{equ:covolfor}
\text{Nr}(I):=\frac{\covol(\iota(I))}{\covol(\iota(\cO_K))}=\Delta_K^{-1/2}\covol(\iota(I)). 
\end{align}
Note that $\text{Nr}$ is multiplicative in $J_K$, that is, $\text{Nr}(I_1I_2)=\text{Nr}(I_1)\text{Nr}(I_2)$ for any $I_1,I_2\in J_K$. Moreover, for any $\alpha\in K\setminus\{0\}$, we have $\text{Nr}(\alpha I)=|\mathrm{N}(\alpha)|\text{Nr}(I)$, where $\mathrm{N}: K\to \R$ is the standard norm on $K$ given by $\mathrm{N}(\alpha):=\alpha\sigma(\alpha)$\label{p:stnorm}.

\subsection{Hilbert modular group}\label{HilbertModGpSEC}
Let $\bH$ be the upper half plane. The group $\SL_2(\R)$ acts on $\bH$ via the M\"{o}bius transformation: $g\cdot z=\frac{az+b}{cz+d}$ for any $g=\left(\begin{smallmatrix}
a & b\\
c & d\end{smallmatrix}\right)\in \SL_2(\R)$ and $z\in \bH$.
Let $\mathrm{H}=\SL_2(\R)\times \SL_2(\R)$\label{p:gph} (which we view as a subgroup of $\SL_4(\R)$ via the block diagonal embedding) and let $\mathbb{H}^2$ be the product of two upper half planes. The group $\mathrm{H}$ naturally acts on $\mathbb{H}^2$ via
\begin{align*}
(h_1,h_2)(z_1,z_2):=(h_1\cdot z_2, h_2\cdot z_2),\qquad \forall\, h_1, h_2\in \SL_2(\R),\, z_1,z_2\in \mathbb{H}.
\end{align*} 
Denote by $\overline{\R}:=\R\cup\{\infty\}$ which can be identified with the boundary of $\mathbb{H}$. The action of $\mathrm{H}$ on $\mathbb{H}^2$ naturally extends to $\overline{\R}^2$ via the same formula. We also denote by $\overline{K}:=K\cup\{\infty\}$. \label{p:overlinek} It embeds into $\overline{\R}^2$ via the natural extension of the embedding $\iota: K\to \R^2$ (by sending $\infty\in \overline{K}$ to $(\infty,\infty)\in \overline{\R}^2$). With slight abuse of notation, we also denote this embedding from $\overline{K}$ to $\overline{\R}^2$ by $\iota$.
We will \textit{also} write $\iota$ for the group homomorphism $\SL_2(K)\to \mathrm{H}$ given by
\begin{align}\label{equ:iotaemb}
\iota\left(\smatr abcd\right)=\left(\smatr abcd,\smatr{\sigma(a)}{\sigma(b)}{\sigma(c)}{\sigma(d)}\right).
\end{align}

The \textit{Hilbert modular group} for the field $K$ is given by \label{equ:slok} 
\begin{align*}
\SL_2(\scrO_K):=\left\{\smatr abcd\in \SL_2(K) \col a,b,c,d\in \cO_K\right\}.
\end{align*}
We write $\Gamma_K$ for the embedding of $\SL_2(\scrO_K)$ in $\mathrm{H}$:\label{p:gammak}
\begin{align*}
\Gamma_K:=\iota\bigl(\SL_2(\scrO_K)\bigr).
\end{align*} 
The discreteness of $\iota(\cO_K)$ in $\R^2$ implies that $\G_K$ is a discrete subgroup of $\mathrm{H}$. Indeed, it is a non-uniform lattice in $\mathrm{H}$, that is, the homogeneous space $\G_K\bk \mathrm{H}$ is non-compact and has finite volume with respect to a Haar measure of $\mathrm{H}$.  
As a subgroup of $\mathrm{H}$, $\G_K$ naturally acts on $\overline{\R}^2$,  
and this action preserves $\iota\bigl(\overline{K}\bigr)$. 
The orbits of $\iota\bigl(\overline{K}\bigr)$ under $\G_K$ are called the \textit{cusps} of $\G_K$.
We will often represent a cusp by an element in the corresponding $\Gamma_K$-orbit, or by an element in $\overline{K}$ (via the natural identification between $\iota(\overline{K})$ and $\overline{K}$). 
The following lemma, together with the obvious identification between $\iota(\overline{K})/\G_K$ and $\overline{K}/\SL_2(\cO_K)$, shows that the number of cusps of $\G_K$ equals the number of ideal classes of $K$.
\begin{Lem}[{\cite[Lemma 3.5]{Freitag1990}}]\label{lem:cuspclnum}
The map sending $k\in \overline{K}$ to $[\langle k,1 \rangle]\in C_K$ induces a bijection between $\overline{K}/\SL_2(\cO_K)$ and $C_K$, where $\langle k, 1\rangle$ is the ideal generated by $k, 1$ if $k\in K$ and $\langle k, 1\rangle:=\cO_K$ if $k=\infty$. In particular,\label{p:kappa}
\begin{align*}
\kappa:=\# C_K=\text{number of cusps of $\G_K$}.
\end{align*}
 \end{Lem}

Throughout the remainder of this paper, we fix $k_1=\infty, k_2,\ldots, k_{\kappa}\in \overline{K}$ to be a complete list of cusps of $\G_K$. \label{p:cusps}
Fix $a_1=1$, $a_2,\ldots, a_\kappa\in\scrO_K$ 
and $c_1=0$, $c_2,\ldots,c_\kappa\in\scrO_K$ 
so that $k_i=\frac{a_i}{c_i}$ for each $i$;
then set  
$I_i:=\langle a_i, c_i\rangle$.
Then $I_1=\cO_K,I_2,\ldots,I_\kappa$ are integral ideals which by Lemma \ref{lem:cuspclnum}
form a system of representatives of the ideal classes of $K$. \label{p:idecuspi}
(Note that $I_i$ depends on the choice of $a_i,c_i$; however the class $[I_i]$ depends only on $k_i$.)

For each $1\leq i\leq \kappa$, let 
\begin{align}\label{equ:isotropygp}
\G_i
:=\left\{
\gamma\in \G_K\col \gamma\iota(k_i)=\iota(k_i)
\right\}
\end{align}
be the \textit{isotropy group} of the cusp $k_i$. 
Below we give a more precise description of these isotropy groups.  First, the isotropy group of $\infty$ is easy to compute: For $\gamma=\iota\left(\smatr abcd\right)\in\Gamma_K$, 
by direct computation $\gamma\cdot \iota(\infty)=\iota(\infty)$ if and only if $c=0$, that is
\begin{align*}
\G_1=\G_K\cap B=\left\{\iota\left(\smatr ab0{a^{-1}}\right)
\col a\in\cO_K^{\times},\ b\in\cO_K\right\},
\end{align*}
where $B<\mathrm{H}$ is the subgroup of upper triangular matrices in $\mathrm{H}$.

To compute the other isotropy groups, we will translate the cusp $k_i$ to $\infty$.  
For any integral ideal $I\subset \cO_K$ let
\begin{align}\label{def:scoki2}
\Gamma_{K,I}
:=\left\{\iota\left(\smatr abcd\right) \col a,d\in \cO_K,\ b\in I^{-2},\ c\in I^2,\ ad-bc=1\right\}. 
\end{align}
For each $1\leq i\leq \kappa$ since $I_i=\la a_i,c_i\ra$, we can find $b_i,d_i\in I_i^{-1}$ such that $a_id_i-b_ic_i=1$. Throughout the paper we fix the \textit{scaling matrix}
\begin{align}\label{def:xidef}
{\xi}_i=\iota\left(\smatr{a_i}{b_i}{c_i}{d_i}\right)\in \mathrm{H}.
\end{align} 
Note that $\xi_i$ satisfies $\xi_i\iota(\infty)=\iota(k_i)$. 
We then have the following description of $\G_i$.
\begin{Lem}\label{lem:desisgp}
We have for each $1\leq i\leq \kappa$,
\begin{align}\label{equ:isoge}
\G_i=\xi_i\left(\Gamma_{K,I_i}\cap B\right)\xi_i^{-1}=\xi_i\left\{
\iota\left(\smatr ab0{a^{-1}}\right)
\col a\in\cO_K^{\times},\ b\in I_i^{-2} \right\}\xi_i^{-1}.
\end{align} 
\end{Lem}
\begin{proof} 
By definition, $\gamma\in \G_i$ if any only if $\gamma \iota(k_i)=\iota(k_i)$. 
Since $\xi_i\iota(\infty)=\iota(k_i)$, this is equivalent to $\xi_i^{-1}\gamma\xi_i\iota(\infty)=\iota(\infty)$.  
This shows that  
$$
{\xi}_i^{-1}\G_i{\xi}_i=({\xi}^{-1}_i\G_K{\xi}_i)\cap B. 
$$
By direction computation we have ${\xi}^{-1}_i\G_K{\xi}_i\subset \G_{K,I_i}$ and $\xi_i\G_{K,I_i}\xi_i^{-1}\subset \G_K$, implying that ${\xi}^{-1}_i\G_K{\xi}_i=\G_{K,I_i}$. Hence 
\begin{align*}
{\xi}_i^{-1}\G_i{\xi}_i=\Gamma_{K,I_i}\cap B=\left\{
\iota\left(\smatr ab0{a^{-1}}\right)
\col a\in\cO_K^{\times},\ b\in I_i^{-2} \right\}.
\end{align*}
We can then finish the proof by conjugating both sides of the above equation by ${\xi}_i$. 
\end{proof}
\subsection{Coordinates and measures}\label{sec:cormea}
Let 
\begin{align*}
\cL_K:=\left\{(\alpha,\beta, \sigma(\alpha),\sigma(\beta))\in\R^4\col (\alpha,\beta)\in \cO_K^2\right\}
\end{align*}
be the Minkowski embedding of $\cO_K^2$ in $\R^4$. 
Let $X$ be the space of lattices of the form $\cL_K h$ with $h\in \mathrm{H}$. Since $\cL_Kh=\cL_K$ if and only if $h\in \G_K$, $X$ can be parameterized by the homogeneous space $\G_K\bk \mathrm{H}$ via $\cL_K h\in X\leftrightarrow \G_K h\in \G_K\bk \mathrm{H}$. We thus identify $X$ with $\G_K\bk \mathrm{H}$. Let $\mu_K$ be the unique invariant probability measure on $X=\G_K\bk \mathrm{H}$.

Since we will be working with this measure extensively when computing $G(s)$ later in Section \ref{sec:asyofgs}, here we give a more explicit description of it in terms of coordinates from an Iwasawa decomposition of $\mathrm{H}$.
It is well known that the group $\SL_2(\R)$ has an Iwasawa decomposition 
saying that any element $g\in \SL_2(\R)$ can be written uniquely as $g=\text{n}_x\text{a}_y\text{k}_{\theta}$ for some $x\in \R, y>0$ and $\theta\in \R/2\pi \Z$, where $\text{n}_x:=\left(\begin{smallmatrix}
1 & x\\
0 & 1\end{smallmatrix}\right)$, $\text{a}_y:=\left(\begin{smallmatrix}
y & 0\\
0 & y^{-1}\end{smallmatrix}\right)$ and $\text{k}_{\theta}:=\left(\begin{smallmatrix}
\cos\theta & -\sin\theta\\
\sin\theta & \cos\theta\end{smallmatrix}\right)$. This thus induces an Iwasawa decomposition of $\mathrm{H}$: Every $h\in \mathrm{H}$ can be written uniquely in the form $h=\text{n}_{\bm{x}}\text{a}_{\bm{y}}\text{k}_{\bm{\theta}}$ with $\bm{x}\in \R^2$, $\bm{y}\in (\R_{>0})^2$ and $\bm{\theta}\in (\R/2\pi\Z)^2$, where\label{p:iwdec}
$$
\text{n}_{\bm{x}}:=\left(\text{n}_{x_1}, \text{n}_{x_2}\right), 
\quad \text{a}_{\bm{y}}:=\left(\text{a}_{y_1}, \text{a}_{y_2}\right)\quad \text{and}\quad \text{k}_{\bm{\theta}}:=(\text{k}_{\theta_1}, \text{k}_{\theta_2}).
$$ 
Under these coordinates, the Haar measure of $\mathrm{H}$ (up to scalars) is given by 
\begin{align}\label{equ:haarh}
\text{d}\mu_{\mathrm{H}}(h)=y_1^{-3}y_2^{-3}\, \text{d}\bm{x}\,\text{d}\bm{y}\,\text{d}\bm{\theta}.
\end{align}
Note that $\mu_K$ essentially comes from a Haar measure of $\mathrm{H}$. Indeed, we may identify $X$ with a fundamental domain inside $\mathrm{H}$; then $\mu_K$ is the restriction of a certain Haar measure to this fundamental domain normalized to be a probability measure.
Thus $\mu_K$ is given by 
\begin{align}\label{equ:haar11}
	\text{d}\mu_K(h)=c_Ky_1^{-3}y_2^{-3}\, \text{d}\bm{x}\,\text{d}\bm{y}\,\text{d}\bm{\theta}
\end{align}
for some normalizing factor $c_K>0$. This normalizing factor was computed by Siegel \cite{Siegel1936,Siegel1937} and is given by the following formula (see \cite[p.\ 59]{Geer1988}):
\begin{align}\label{equ:compck}
c_K=\Delta_K^{-3/2}\zeta_K(2)^{-1},
\end{align}
where $\zeta_K(s):=\sum_{I\subset \cO_K}\text{Nr}(I)^{-s}$ ($\Re(s)>1$) is the Dedekind zeta function attached to $K$. Here the summation is over all the nonzero integral ideals of $\cO_K$\label{p:dezefun}. 
\begin{remark}\label{rmk:zeta2}
To compare our formula for $a_{\cP}$ in \eqref{equ:genforapclassno1} with numerical computations, we need to express $\zeta_K(2)$ in more explicit terms. Let $\chi_K$ be the Dirichlet character defined by $\chi_K(n):=\left(\frac{\Delta_K}{n}\right)$ with $\left(\frac{\cdot}{\cdot}\right)$ the Kronecker symbol. Note that $\chi_K$ is an even primitive quadratic character of modulus $\Delta_K$; see \cite[p.\ 296-297]{MontgomeryVaughan2007}.  
Using the fact that for any prime number $p$, the ideal $(p)\subset \cO_K$ is inert, ramified or split if and only if $\chi_K(p)$ equals $-1, 0,1$ respectively, we have the following well-known formula for $\zeta_K$:
\begin{align*}
\zeta_K(s)=\zeta(s)L(s,\chi_K),
\end{align*}
where $\zeta(s)$ is the Riemann zeta function and $L(s,\chi_K)$ is the Dirichlet $L$-function associated to $\chi_K$.
Now by \cite[Prop.\ 4.1 and Thm.\ 4.2]{Washington1997} we have for any $n\in \N$,
\begin{align*}
L(1-n,\chi_K)=-\frac{\Delta_K^{n-1}}{n}\sum_{a=1}^{\Delta_K}\chi_K(a)B_n\left(\tfrac{a}{\Delta_K}\right),
\end{align*}
where $B_n(X)$ is the $n$-th Bernoulli polynomial given by the formula $B_n(X)=\sum_{i=0}^n\binom{n}{i}B_iX^{n-i}$ with $B_i$ the $i$-th Bernoulli number. 
Combined with the functional equation of $\zeta_K(s)$ (see \cite[p.\ 59]{Geer1988}), we get
\begin{align}\label{equ:zetak2}
\zeta_K(2)=4\pi^4\Delta_K^{-\frac32}\zeta_K(-1)=4\pi^4\Delta_K^{-\frac32}\zeta(-1)L(-1,\chi_K)=\frac{\pi^4}{6\sqrt{\Delta_K}}\sum_{a=1}^{\Delta_K}\chi_K(a)B_2\left(\tfrac{a}{\Delta_K}\right).
\end{align}
For instance when $K=\Q(\sqrt{2})$ we have $\Delta_K=8$ and $\chi_K(n)=\left(\frac{8}{n}\right)$ is the quadratic character of modulus $8$ with $\chi_K(1)=\chi_K(7)=1$ and $\chi_K(3)=\chi_K(5)=-1$. We also have $B_2(1/8)=B_2(7/8)=\frac{11}{192}$ and $B_2(3/8)=B_2(5/8)=-\frac{13}{192}$. Plugging all these relations we get in this case $\zeta_K(2)=\frac{\pi^4}{48\sqrt{2}}$. Similarly, one can apply \eqref{equ:zetak2} to get $\zeta_{\Q(\sqrt{3})}(2)=\frac{\pi^4}{36\sqrt{3}}$ and $\zeta_{\Q(\sqrt{5})}(2)=\frac{2\pi^4}{75\sqrt{5}}$.
\end{remark}

\subsection{Siegel domain and cusp neighborhoods}\label{sec:siegel} 
In this section we give a more precise description of the homogeneous space $X=\G_K\bk \mathrm{H}$ in terms of the coordinates given in the previous section.

First, note that under these coordinates and in view of Lemma \ref{lem:desisgp}, for each $1\leq i\leq \kappa$  
the set
\begin{align}\label{equ:cusp1}
\cF_i:=\left\{{\xi}_i\text{n}_{\bm{x}}\text{a}_{\bm{y}}\text{k}_{\bm{\theta}}\in \mathrm{H}\col \bm{x}\in \mathfrak{F}_i,\ y_1/y_2\in [1,\lambda^2),\ \bm{\theta}\in [0, \pi)\times [0,2\pi) \right\}
\end{align}
is a fundamental domain for $\G_i\bk \mathrm{H}$. Here $\mathfrak{F}_i\subset \R^2$ is a fundamental domain for $\R^2/\iota(I_i^{-2})$\label{p:fFi}. (Recall that $\iota(I)$ is a lattice in $\R^2$ for any fractional ideal $I$ of $\cO_K$.)

Now for each $1\leq i\leq \kappa$ and any $t>0$ define
\begin{align}\label{equ:cusp2}
\cF_i(t):=\left\{\xi_i\text{n}_{\bm{x}}\text{a}_{\bm{y}}\text{k}_{\bm{\theta}}\in \cF_i\col y_1y_2\geq t\right\}.
\end{align}
It follows from Shimizu's lemma \cite[Lemma 5]{Shimizu1963} that there exists some constant $t_1>0$ depending only on $K$ such that \label{p:t1}
\begin{align}\label{equ:cuspdescp}
\forall\, \gamma \in \G_K, 1\leq i,j\leq \kappa\col \gamma\cF_i(t_1)\cap \cF_j(t_1) \neq \emptyset\quad \Leftrightarrow\quad i=j\  \text{and}\ \gamma=\text{id}.
\end{align} 
Moreover, by the reduction theory of Borel and Harish-Chandra \cite{BorelHC1962} we have a Siegel fundamental domain of the form
\begin{align}\label{equ:siegefd}
\cF_{\G_K}:=\mathfrak{C}\bigcup\left(\bigsqcup_{i=1}^{\kappa}\cF_i(t_1)\right).
\end{align}
Here $\mathfrak{C}\subset \mathrm{H}$ is compact and the natural map from $\cF_{\G_K}$ to $\G_K\bk \mathrm{H}$ is surjective and finite-to-one. Note, on the other hand, it follows from \eqref{equ:cuspdescp} that the projection from $\bigsqcup_{i=1}^{\kappa}\cF_i(t_1)$ to $\G_K\bk \mathrm{H}$ is injective.

\subsection{Geometry of lattices avoiding large balls}\label{sec:geomlavo}
In this section we give necessary conditions for lattices $\cL_K h\in X$ avoiding large balls. These lattices will be the main objects we deal with when computing $G(s)$. We have the following description of these lattices. 
\begin{Prop}\label{prop:bigball} 
There exists some $R_0>0$ depending only on $K$ such that for any $R>R_0$ and for any lattice $\cL_K h\in X$ avoiding a ball of radius $R$, we have $\cL_K h=\cL_K \xi_i\text{n}_{\bm{x}}\text{a}_{\bm{y}}\text{k}_{\bm{\theta}}\in \G_K\cF_i(t_1)$ for some $1\leq i\leq \kappa$. Moreover, we have $y_1\asymp y_2\gg R$ with the bounding constants depending only on $K$.
\end{Prop}
\begin{remark}
Regarding the notation ``$\asymp$'', ``$\gg$'' and ``$\ll$'':
For two positive quantities $A$ and $B$, we write $A\ll B$ or $B\gg A$ to mean that
there is a constant $C>0$ such that $A\leq CB$,
and we will write $A\asymp B$ for $A\ll B\ll A$.
We will sometimes use subscripts to indicate the dependence of the bounding constant $C$ on 
parameters.
\end{remark}
\begin{proof}
Let $\bm{v}_1=(1,0,1,0)$, $\bm{v}_2=(\tau,0,\sigma(\tau),0)$, $\bm{v}_3=(0,1,0,1)$ and $\bm{v}_4=(0,\tau,0,\sigma(\tau))$. Note that $\{\bm{v}_j\col  j=1,2,3,4\}$ is a basis for $\cL_K$. Let $R_0>0$ be sufficiently large such that 
\begin{align}\label{equ:R0}
2\|\bm{v}_jh\|\leq R_0,\qquad \forall\,1\leq j\leq 4,\ h\in \mathfrak{C}.
\end{align}
We note that $R_0$ depends only on the compact set $\mathfrak{C}$ (hence only depends on $K$). Now take $R>R_0$ and suppose $x\in X$ avoids a ball $B_R$ of radius $R$. To prove the first claim, in view of the description of the Siegel fundamental domain $\cF_{\G_K}$ in  \eqref{equ:siegefd}, 
we want to show there does not exist any $h\in \mathfrak{C}$ such that $x=\cL_K h$. We prove by contradiction. Suppose there exists such $h \in \mathfrak{C}$. Then $\{\bm{v}_jh\col 1\leq j\leq 4\}$ is a basis for $x=\cL_K h$. 
Let $\bm{x}\in \R^4$ be the center of $B_R$ and write $\bm{x}=\sum_{j=1}^4a_j\bm{v}_jh$ ($a_j\in\R$) as an $\R$-linear combination of these basis vectors. For each $1\leq j\leq 4$, take $n_j\in\Z$ the closest integer to $a_j$ so that we have $|n_j-a_j|\leq \frac12$. Consider the lattice vector $\bm{x}':=\sum_{j=1}^4n_j\bm{v}_jh\in \cL_K h$. Since $\cL_K h\cap B_R=\emptyset$, we have 
\begin{align*}
R\leq \|\bm{x}'-\bm{x}\|\leq \sum_{j=1}^4|n_j-a_j|\|\bm{v}_jh\|\leq \frac12\sum_{i=j}^4\|\bm{v}_jh\|,
\end{align*}
implying that there exists $1\leq j\leq 4$ such that $\|\bm{v}_jh\|\geq \frac{R}{2}>\frac{R_0}{2}$, contradicting \eqref{equ:R0}. This proves the first claim, i.e.  
$\cL_K h\in \G_K\cF_i(t_1)$ for some $1\leq i\leq \kappa$. For the second claim, without loss of generality we may assume $h=\xi_i\text{n}_{\bm{x}}\text{a}_{\bm{y}}\text{k}_{\bm{\theta}}\in \cF_i(t_1)$. To show $y_1\asymp y_2\asymp R$, first note that $y_1\asymp_{\lambda}y_2$ and $y_1y_2\geq t_1\gg 1$. Thus
\begin{align*}
y_1\asymp y_2 \asymp \max\{y_1,y_1^{-1},y_2,y_2^{-1}\}=\|\text{a}_{\bm{y}}\|_{\rm op}.
\end{align*}
Here $\|\cdot\|_{\rm op}$ denotes the operator norm (with respect to the Euclidean norm).
On the other hand, for any $1\leq j\leq 4$, 
\begin{align*}
\|\bm{v}_jh\|\leq \|\bm{v}_j\|\|\xi_i\|_{\rm op}\|\text{n}_{\bm{x}}\|_{\rm op}\|\text{a}_{\bm{y}}\|_{\rm op}\ll_{\mathfrak{F}_i,\xi_i} \|\text{a}_{\bm{y}}\|_{\rm op}\asymp y_1\asymp y_2.
\end{align*}
This, together with the bound $\sum_{j=1}^4\|\bm{v}_jh\|\geq 2R$ implies that $y_1\asymp y_2\gg R$ as claimed. 
\end{proof}

For later reference, note that for each $1\leq i\leq \kappa$,
\begin{align}\label{equ:lalterdes}
\cL_K=\left\{(\alpha,\beta, \sigma(\alpha), \sigma(\beta))\xi_i^{-1}\col \alpha\in I_i,\ \beta\in I_i^{-1}\right\}
\end{align}
with $I_i=\langle a_i, c_i\rangle$ as before.
Let $\bm{\fc}_i=(-c_i,a_i)\in \cO_K^2$\label{p:frakci}. We denote by
\begin{align}\label{equ:fillplane1}
\Pi_i:=\mathbb{R}\iota(\bm{\fc}_i)\oplus \R \iota(\tau\bm{\fc}_i)=
\left(\R(0,1,0,1)\oplus \R(0, \tau, 0,\sigma(\tau))\right)\xi^{-1}_i=(\{0\}\hspace{-3pt}\times\hspace{-3pt} \R\hspace{-3pt}\times\hspace{-3pt}\{0\}\hspace{-3pt}\times\hspace{-3pt}\R)\xi^{-1}_i,
\end{align}
and 
\begin{align}\label{equ:fillplane2}
\cL_i:=\Pi_i\cap \cL_K&
=\left((\{0\}\hspace{-3pt}\times\hspace{-3pt} \R\hspace{-3pt}\times\hspace{-3pt}\{0\}\hspace{-3pt}\times\hspace{-3pt}\R)
\cap \cL_K \xi_i\right)\xi^{-1}_i
=\left\{(0,\beta,0,\sigma(\beta))\xi^{-1}_i\col \beta\in I_i^{-1} \right\}.
\end{align}
By Proposition \ref{prop:bigball} if $x\in X$ avoids some ball of radius $R>R_0$, then $x=\cL_K h$ with $h=\xi_i\text{n}_{\bm{x}}\text{a}_{\bm{y}}\text{k}_{\bm{\theta}}\in \cF_i(t_1)$ for some $1\leq i\leq \kappa$ and $y_1,y_2>0$ with $y_1\asymp y_2\gg R$. Then the two vectors 
\begin{align*}
\iota(\bm{\fc}_i)h=(0,y_1^{-1},0, y_2^{-1})\text{k}_{\bm{\theta}}\quad \text{and}\quad  \iota(\tau\bm{\fc}_i)h=(0,\tau y_1^{-1},0, \sigma(\tau)y_2^{-1})\text{k}_{\bm{\theta}}
\end{align*}
are of length $\ll R^{-1}$. Note that $\iota(\bm{\fc}_i)h$ and $\iota(\tau\bm{\fc}_i)h$ are linearly independent and $\cL_i h$ is 
the primitive rank two sublattice 
of $x$ containing them. Thus $\cL_i h$ lies $O(R^{-1})$-densely in the plane $\Pi_i h$. For this reason, we call $\Pi_i h$ the \textit{filled plane} of $x=\cL_K h$.

\section{Asymptotics of $G(s)$}
\label{sec:asyofgs}
\subsection{Main result}
The main goal of this Section \ref{sec:asyofgs} 
is to prove Theorem~\ref{MAINTHM1}, 
which gives the tail asymptotics of both the limiting gap distribution function $F(s)$ and its
integral $G(s)$,
for the particular class of cut-and-project sets which we consider.
In fact we will prove 
Theorem \ref{MAINTHM1moreexplict} below,
which will be shown (in Section \ref{MAINTHM1prooffrommoreexplSEC})
to imply Theorem~\ref{MAINTHM1}.
This Theorem \ref{MAINTHM1moreexplict} gives 
an explicit bound on the error term
in the asymptotics for $G(s)$,
and also an explicit formula for the leading coefficient $a_{\scrP}$.
As we mentioned in the introduction, 
the asymptotics for $F(s)$ will be derived as a consequence of those for $G(s)$;
see Lemma \ref{prop:gtof} below.

In order to state Theorem \ref{MAINTHM1moreexplict}, we introduce some further notation.
For any $\theta\in\R$ we write $\text{k}_{\theta}:=\left(\begin{smallmatrix}
\cos\theta & -\sin\theta\\
\sin\theta & \cos\theta\end{smallmatrix}\right)\in\SL_2(\R)$.
Then for any $\theta,x\in\R$ and any subset $\scrW\subset\R^2$,
we set
\begin{align}\label{RWdef}
R_{\scrW}(\theta,x):=\{y\in\R\col (x,y)\text{k}_\theta\in\scrW\}.
\end{align}
It follows that 
$\ell_{\cW}(\theta)$, which we defined in Section~\ref{sec:forap}, can be expressed as 
\begin{align}\label{ellWdef}
\ell_{\scrW}(\theta)=\bigl\{x\in\R\col R_{\scrW}(\theta,x)\neq \emptyset\bigr\}.
\end{align} 
Note also that $\ell_{\cW}(\theta)$ satisfies the relation $\ell_{\cW}(\theta)=-\ell_{\cW}(\pi+\theta)$ for any $\theta\in \R$.
For any integral ideal $I\subset \cO_K$, any
bounded Lebesgue measurable set $J\subset\R$ with non-empty interior,
and any $y>0$, we set
\begin{align}\label{alphapWiDEF}
\alpha_{I,J}(y):=\min\bigl\{\alpha\in I\cap\R_{>0}\col y\,\sigma(\alpha)\in J\bigr\}.
\end{align}
For later reference we record here a useful identity for the function $\alpha_{I,J}$ which can be checked directly from its definition: For any $\beta\in I^{-1}$ with $\sigma(\beta)>0$ and for any $y>0$,
\begin{align}\label{equ:alphaprop}
\alpha_{\beta I,J}(y)=|\beta|\alpha_{I,\sgn(\beta)J}(\sigma(\beta)y).
\end{align} 
Let us write $m$ for the Lebesgue measure on $\R$. \label{mLEBdef}
For any Lebesgue measurable set $A\subset\R$ and $a>0$, we define
$\ff(A,a)$ to be the infimum of $a^{-1}m(J\setminus A)$ when $J$ ranges over all intervals of length $a$:
\begin{align}\label{ffAaDEF}
\ff(A,a): 
&=a^{-1}\inf\bigl\{m\bigl((x,x+a)\setminus A\bigr)\col x\in\R\bigr\}. 
\end{align}
Note that $0\leq\ff(A,a)\leq1$ always, and 
$\ff(A,a)=0$ whenever $A$ contains some interval of length $a$.
Let us also define
\begin{align}\label{tffAaDEF}
\tff(A,a):=\begin{cases}\,\ff(A,a)&\text{if }\: A\neq \emptyset,
\\
0&\text{if }\: A=\emptyset.
\end{cases}
\end{align}

Recall from below Lemma \ref{lem:cuspclnum}
that we have made a fixed choice of 
ideals $I_1,\ldots,I_{\kappa}$ of $\scrO_K$,
and that these form a system of representatives of the ideal classes of $K$.

\begin{thm}\label{MAINTHM1moreexplict}
Let $\scrP=\scrP(\scrW,\scrL_K)$
where $K$ is a real quadratic field
and $\scrW$ is a bounded open subset of $\R^2$ such that
$\bn\in\scrW$ and $\partial\scrW$ has Lebesgue measure zero.
Let $F(s)$ be the associated limiting gap distribution function
as in Theorem \ref{gaplimitdistrexistsTHM}, and let $G(s)=\int_s^\infty F(t)\,\textup{d}t$.
Then there exists a positive constant $C=C_{\cP}$ such that
for all sufficiently large $s$,
\begin{align}\label{MAINTHM1res1}
0\leq G(s)-\frac{a_{\scrP}}{s}
\ll  s^{-2}\int_0^{Cs}\int_{0}^{2\pi}\int_{\R}\tff\biggl(R_{\scrW}(\theta,y),u^{-1}\biggr)\,\mathrm{d}y\,\mathrm{d}\theta\,\mathrm{d}u,
\end{align}
where
\begin{align}\label{MAINTHM1res2}
a_{\scrP}=\frac{\Area(\scrW)}{4 \Delta_K^{2}\zeta_K(2)}\sum_{i=1}^{\kappa}\operatorname{Nr}(I_i)^{-2}\int_0^{2\pi}\int_1^{\lambda^2}
\alpha_{I_i,\ell_{\cW}(\theta)}(y)^2\,\frac{\mathrm{d}y}{y^3}\, \mathrm{d}\theta.
\end{align}
\end{thm}

We here make a couple of remarks regarding Theorem \ref{MAINTHM1moreexplict};
for two more remarks see the end of Section \ref{MAINTHM1prooffrommoreexplSEC}.
\begin{remark}\label{idclassinvREM}
In the formula \eqref{MAINTHM1res2}, it should be noted that for any integral ideal $I$ of $\scrO_K$,
we have
\begin{align}\label{idclassinvREM2}
\int_1^{\lambda^2}\alpha_{I,\ell_{\cW}(\theta)}(y)^2\,\frac{\text{d}y}{y^3}
=\int_a^{a\lambda^2}\alpha_{I,\ell_{\cW}(\theta)}(y)^2\,\frac{\text{d}y}{y^3},
\qquad \forall\, a>0,\:\theta\in\R,
\end{align}
and furthermore, the product 
\begin{align}\label{idclassinvREMpf1}
\operatorname{Nr}(I)^{-2}
\int_0^{2\pi}\int_1^{\lambda^2}
\alpha_{I,\ell_{\cW}(\theta)}(y)^2\,\frac{\text{d}y}{y^3}\, \text{d}\theta
\end{align}
depends only on the ideal class $[I]$ of $I$.
\\[3pt]
To prove \eqref{idclassinvREM2}, it suffices to note that
the differential form $\alpha_{I,\ell_{\cW}(\theta)}(y)^2\,\frac{\text{d}y}{y^3}$ is invariant under the map sending $y$ to $\lambda^2 y$,
a fact which follows by applying \eqref{equ:alphaprop} for the special case when $\beta=\lambda^{-2}$.
Next, to prove the invariance of \eqref{idclassinvREMpf1},
assume that $I'$ is another integral ideal belonging to the same ideal class as $I$.
Then $I'=aI$ for some nonzero $a\in I^{-1}$, which we may take to satisfy $\sigma(a)>0$ (otherwise replace $a$ by $-a$).
Then the absolute norm $\operatorname{Nr}(I_i)^{-2}$ scales by $\mathrm{N}(a)^{-2}$ when replacing $I$ by $I'=aI$.
On the other hand, using \eqref{equ:alphaprop} and the fact that $\ell_{\cW}(\theta+\pi)=-\ell_{\cW}(\theta)$ for any $\theta\in \R$,
we have 
\begin{align*}
\alpha_{I',\ell_{\cW}(\theta)}(y)=\begin{cases}
	a\alpha_{I,\ell_{\cW}(\theta)}(\sigma(a)y) & \text{if $a>0$},\\[2pt]
	|a|\alpha_{I,\ell_{\cW}(\theta+\pi)}(\sigma(a)y) & \text{if $a<0$},
\end{cases}
\end{align*}
from which it follows (via the obvious substitution and then using
\eqref{idclassinvREM2} with $\sigma(a)$ in place of $a$)
that the double integral in \eqref{idclassinvREMpf1} scales by a factor
$\mathrm{N}(a)^2$ when replacing $I$ by $I'=aI$. Hence the product in \eqref{idclassinvREMpf1} is invariant as claimed.
\end{remark}

\begin{remark}\label{rmk1:sl2inva} 
The right hand side of \eqref{MAINTHM1res2}
is invariant under replacing $\scrW$ by $\scrW g$, for any $g\in\GL_2(\R)$.

\begin{proof}
We will prove the claim by showing that for any integral ideal $I$ of $\scrO_K$,
and any $g\in\GL_2(\R)$,
\begin{align}\label{rmk1:sl2invapf1}
\Area(\scrW g)\int_0^{2\pi}\int_1^{\lambda^2}\alpha_{I,\ell_{\cW g}(\theta)}(y)^2\,\frac{\text{d}y}{y^3}\, \text{d}\theta
=\Area(\scrW)\int_0^{2\pi}\int_1^{\lambda^2}\alpha_{I,\ell_{\cW}(\theta)}(y)^2\,\frac{\text{d}y}{y^3}\, \text{d}\theta.
\end{align}
Write $\vece_1=(1,0)$, $\vece_2=(0,1)$.
For fixed $g\in\GL_2(\R)$, 
define the smooth maps $\omega:\R/2\pi\Z\to\R/2\pi\Z$ and $\rho:\R/2\pi\Z\to\R_{>0}$
through
\begin{align}\label{rmk1:sl2invapf2}
\vece_2\text{k}_{\theta}g^{-1}=\rho(\theta)\vece_2 \text{k}_{\omega(\theta)}.
\end{align}
Note that $\omega$ is a diffeomorphism of  
$\R/2\pi\Z$ onto itself.
Differentiating both sides of \eqref{rmk1:sl2invapf2} 
with respect to $\theta$,
and using the fact that 
$\frac{\text{d}}{\text{d}\theta}\text{k}_\theta=\text{k}_{\theta+\frac\pi2}$,
we have
\begin{align}\label{rmk1:sl2invapf4}
\vece_1\text{k}_{\theta}g^{-1}=\bigl(f(\theta),\rho'(\theta)\bigr) \text{k}_{\omega(\theta)},
\qquad\text{where }\:
f(\theta):=\rho(\theta)\omega'(\theta).
\end{align}
It also follows from \eqref{rmk1:sl2invapf4} and \eqref{rmk1:sl2invapf2} that
\begin{align}\notag
f(\theta) =\vece_1\text{k}_{\omega(\theta)}\cdot\vece_1\text{k}_\theta g^{-1}
=\vece_2\text{k}_{\omega(\theta)}\text{k}_{\frac{\pi}2}\cdot\vece_1\text{k}_\theta g^{-1}
&=\rho(\theta)^{-1}\vece_2\text{k}_{\theta}g^{-1}\text{k}_{\frac{\pi}2}\cdot\vece_1\text{k}_\theta g^{-1}
\\\label{rmk1:sl2invapf5}
&=\rho(\theta)^{-1}\det(\text{k}_\theta g^{-1})=\rho(\theta)^{-1}\det(g)^{-1},
\end{align}
where in the fourth equality we used the fact that
$\vecw\text{k}_{\frac{\pi}2}\cdot\vecv=v_1w_2-v_2w_1$ 
for any $\vecv=(v_1,v_2)$ and $\vecw=(w_1,w_2)$ in $\R^2$.
Note that \eqref{rmk1:sl2invapf4} and \eqref{rmk1:sl2invapf5} imply that
if $\det(g)>0$ then $f(\theta)>0$ and $\omega'(\theta)>0$ for all $\theta$,
while if $\det(g)<0$ then $f(\theta)<0$ and $\omega'(\theta)<0$ for all $\theta$.

Using \eqref{rmk1:sl2invapf2} and \eqref{rmk1:sl2invapf4} we have
$(x,y)\text{k}_\theta g^{-1}=\bigl(xf(\theta),x\rho'(\theta)+y \rho(\theta)\bigr)\text{k}_{\omega(\theta)}$
for any $(x,y)\in\R^2$,
and hence via \eqref{RWdef} and \eqref{ellWdef} one gets
$\ell_{\scrW g}(\theta)=f(\theta)^{-1}\ell_{\scrW}(\omega(\theta))$.
We also have $\omega(\theta+\pi)=\omega(\theta)+\pi$, $f(\theta+\pi)=f(\theta)$,
and $\ell_{\scrW}(\omega(\theta)+\pi)=-\ell_{\scrW}(\omega(\theta))$,
and hence:
\begin{align*}
\ell_{\scrW g}(\theta)=|f(\ttheta)|^{-1}\ell_{\scrW}(\omega(\ttheta)),
\qquad\text{where $\ttheta:=\theta$ if $\det(g)>0$, $\ttheta:=\theta+\pi$ if $\det(g)<0$.}
\end{align*}
Using also the simple relation
$\alpha_{I,aJ}(y)=\alpha_{I,J}(a^{-1}y)$ ($\forall\, a>0$),
it follows that 
\begin{align}\notag
\Area(\scrW g)\int_0^{2\pi}\int_1^{\lambda^2}\alpha_{I,\ell_{\cW g}(\theta)}(y)^2\,\frac{\text{d}y}{y^3}\, \text{d}\theta
=\Area(\scrW g)\int_{\R/2\pi\Z}\int_1^{\lambda^2}\alpha_{I,\ell_{\cW}(\omega(\ttheta))}\bigl(|f(\ttheta)|y\bigr)^2\,\frac{\text{d}y}{y^3}
\, \text{d}\ttheta
\\\label{rmk1:sl2invapf3}
=|\det g|
\Area(\scrW)\int_{\R/2\pi\Z}f(\ttheta)^2
\int_1^{\lambda^2}\alpha_{I,\ell_{\cW}(\omega(\ttheta))}(u)^2\,\frac{\text{d}u}{u^3}\, \text{d}\ttheta,
\end{align}
where we substituted $u=|f(\ttheta)|y$ and then used \eqref{idclassinvREM2}. 
The two formulas for $f(\theta)$ in 
\eqref{rmk1:sl2invapf4} and \eqref{rmk1:sl2invapf5}
imply that
$f(\ttheta)^2=\det(g)^{-1}\omega'(\ttheta)$.
Using this in the last expression in  
\eqref{rmk1:sl2invapf3}, 
and then taking $\omega$ as a new variable of integration,
we obtain the identity \eqref{rmk1:sl2invapf1}. 
\end{proof}

Let us also note that 
for any $g\in\SL_2(\R)$,
the invariance statement in Remark \ref{rmk1:sl2inva}
may alternatively be deduced  
as a consequence of
\eqref{MAINTHM1res1} in Theorem \ref{MAINTHM1moreexplict}
and Remark \ref{rmk:sl2inva} below.
Once the invariance is known for all $g\in\SL_2(\R)$,
in order to extend it to all $g\in\GL_2(\R)$ it suffices to verify that 
it also holds for $g=\diag(1,-1)$ and $g=\diag(a,a)$ ($\forall\, a>0$);
for these cases the above proof applies in a significantly simplified form.  
\end{remark}

\subsection{Proof of Theorem \ref*{MAINTHM1} using Theorem \ref*{MAINTHM1moreexplict}}\label{MAINTHM1prooffrommoreexplSEC}
In this section we prove Theorem \ref{MAINTHM1} assuming Theorem \ref{MAINTHM1moreexplict}.
We will first prove the two tail asymptotics of $G(s)$ in \eqref{equ:gsestgenwin} and \eqref{equ:polyest};
then, by a standard general argument using the convexity of $G(s)$,
we will deduce the tail asymptotics of $F(s)$.

We start with \eqref{equ:gsestgenwin}. 
We need to show that
for any $\scrW$ as in Theorem \ref{MAINTHM1},\label{errortermalwaysokCORpf}
the error bound in \eqref{MAINTHM1res1} is $o(s^{-1})$ as $s\to\infty$.
Take $R>0$ so that $\scrW\subset\scrB_R^2$, where $\scrB_R^2\subset \R^2$ is the open disc with center $\bm{0}$ and radius $R$\label{p:ball}.
Then for every $y\in\R$ with $|y|\geq R$, and for every $\theta$, 
we have $R_{\scrW}(\theta,y)=\emptyset$
and hence $\tff\bigl(R_{\scrW}(\theta,y),u^{-1}\bigr)=0$ for all $u>0$.
It follows that the integral over $\R$ in 
\eqref{MAINTHM1res1} can be replaced by $\int_{-R}^R$.
Hence, if we take $C=C_{\cP}$ as in Theorem \ref{MAINTHM1moreexplict}, and set
\begin{align*}
F_{\theta,y}(s):=s^{-1}\int_0^{Cs}\tff\bigl(R_{\scrW}(\theta,y),u^{-1}\bigr)\,\text{d}u,
\end{align*}
then our task is to prove that 
$\int_0^{2\pi}\int_{-R}^R F_{\theta,y}(s)\,\text{d}y\,\text{d}\theta\to0$ as $s\to\infty$.
Using $0\leq\tff\bigl(R_{\scrW}(\theta,y),u^{-1}\bigr)\leq1$ we have
$0\leq F_{\theta,y}(s)\leq C$ for all $\theta$, $y$ and $s$;
hence by the Lebesgue bounded convergence theorem,
it suffices to prove that 
$F_{\theta,y}(s)\to0$ as $s\to\infty$
for any fixed $\theta,y$.
But if $R_{\scrW}(\theta,y)=\emptyset$ then
$F_{\theta,y}(s)=0$ for all $s$, by the definition \eqref{tffAaDEF},
while if $R_{\scrW}(\theta,y)\neq\emptyset$ then
since $\scrW$ is open, 
$R_{\scrW}(\theta,y)$ contains a non-empty interval,
and so $\tff\bigl(R_{\scrW}(\theta,y),u^{-1}\bigr)=0$ for all sufficiently large $u$;
therefore $F_{\theta,y}(s)\to0$ as $s\to\infty$.
This concludes the proof that 
the error bound in \eqref{MAINTHM1res1} is $o(s^{-1})$ as $s\to\infty$,
i.e.\ we have proved that $G(s)\sim a_{\scrP}s^{-1}$, as stated in 
\eqref{equ:gsestgenwin}.

For the proof of \eqref{equ:polyest} 
we need the following simple lemma regarding the function $\tff(A,a)$.
\begin{Lem}\label{lem:ffprop}
Let $C>0$.  
Then for any interval $A$ of length $b>0$ and
for any $s>0$, 
\begin{align*}
\int_0^{Cs}\tff(A,u^{-1})\,\mathrm{d}u\asymp_C\min(s,b^{-1}).
\end{align*}
\end{Lem}
\begin{proof}
Note that
$\ff(A,a)=0$ for $a\in(0,b]$ and
$\ff(A,a)=1-(b/a)$ for $a\in[b,\infty)$.
Also $\tff(A,a)=\ff(A,a)$.
Hence for $s\leq (Cb)^{-1}$:
\begin{align*}
\int_0^{Cs}\tff(A,u^{-1})\,\text{d}u
=\int_0^{Cs}\bigl(1-bu\bigr)\,\text{d}u\asymp Cs\asymp_C s,
\end{align*}
while for $s\geq (Cb)^{-1}$:
\begin{align*}
\int_0^{Cs}\tff(A,u^{-1})\,\text{d}u
=\int_0^{b^{-1}}\bigl(1-bu\bigr)\,\text{d}u+0\asymp b^{-1}.
\end{align*}
The desired estimate then follows from these two estimates noting that 
$\min(s, b^{-1})\asymp_C s$ for $s\leq (Cb)^{-1}$ and $\min(s, b^{-1})\asymp_C b^{-1}$ for $s\geq (Cb)^{-1}$.
\end{proof}

Now let us further assume that $\cW$ is convex.
Then for all $\theta$ and $y$,
$R_{\scrW}(\theta,y)$ is either empty or an interval.
Hence, using Lemma \ref{lem:ffprop} and
the fact that $m(R_{\scrW}(\theta+\pi,y))=m(R_{\scrW}(\theta,-y))$,
it follows that the error bound in \eqref{MAINTHM1res1} is
\begin{align}\label{MAINTHM1pf1}
\asymp s^{-2}\int_0^{2\pi}\int_0^\infty I\bigl(R_{\scrW}(\theta,y)\neq\emptyset\bigr)
\cdot\min\Bigl(s,m(R_{\scrW}(\theta,y))^{-1}\Bigr)\, \text{d}y\,\text{d}\theta.
\end{align}
Choose $0<r_0\leq R_0$ such that $\scrB_{r_0}^2\subset\cW\subset\scrB_{R_0}^2$.
Let us temporarily fix $\theta\in[0,2\pi]$,
and set $y_0:=\sup\{y>0\col R_{\scrW}(\theta,y)\neq\emptyset\}$;
then the function $f(y):=m(R_{\scrW}(\theta,y))$ vanishes for $y>y_0$, 
but is positive and concave for $y\in[0,y_0)$.
Hence for our fixed $\theta$,
the inner integral in \eqref{MAINTHM1pf1} equals
$\int_0^{y_0}\min\bigl(s,f(y)^{-1}\bigr)\,\text{d}y$.
We also have $r_0\leq y_0\leq R_0$ and $f(0)\geq r_0$,
and so, by the concavity,
$f(y)\geq (r_0/R_0)\cdot (y_0-y)$ for all $y\in[0,y_0]$.
Hence
\begin{align*}
\int_0^{y_0}\min\bigl(s,f(y)^{-1}\bigr)\,\text{d}y
\leq\int_0^{R_0}\min\Bigl(s,\frac{R_0}{r_0\,t}\Bigr)\,\text{d}t,
\end{align*}
and by a simple computation, the last integral is seen to equal
${\displaystyle \frac{R_0}{r_0}}\bigl(\log(r_0s)+1\bigr)$
for all $s\geq 1/r_0$.
Since this holds for every $\theta\in[0,2\pi]$,
we conclude that the expression in \eqref{MAINTHM1pf1} 
is $\ll s^{-2}\log s$ for all large $s$,
i.e.\ we have proved the asymptotics for $G(s)$ stated in 
\eqref{equ:polyest}.

\vspace{3pt}

The next lemma
shows how to use the convexity of $G(s)$ to deduce the tail asymptotics of
$F(s)=-G'(s)$ from that of $G(s)$.
\begin{Lem}\label{prop:gtof}
Let $G: \R_{>0}\to \R_{>0}$ be a differentiable function with $F(s):=-G'(s)$ continuous and decreasing. 
Suppose further that there exist constants $a>0$ and $C_1,C_2>1$ 
and a function $E:(C_1,\infty)\to\R_{>0}$, such that
\begin{align}\label{Gsasymptgenform}
\bigl|G(s)-as^{-1}\bigr|\leq E(s)\leq\frac1{s},\qquad\forall\, s>C_1,
\end{align}
and
\begin{align}\label{equ:assupones}
C_2^{-1}\leq\frac{E(s)}{E(s')}\leq C_2\qquad \text{whenever $C_1<s\leq s'\leq 2s$.}
\end{align}
Then we have
\begin{align}\label{equ:fsasymp}
F(s)=as^{-2}+O(s^{-\frac{3}{2}}E(s)^{\frac12})\quad \text{as $s\to\infty$},
\end{align}
where the implied constant depends only on $a$ and $C_2$.
\end{Lem}
\begin{proof}
Since by assumption $F(s)$ is continuous and decreasing, we have for all 
$0<s_1<s_2$, $G(s_1)-G(s_2)=\int_{s_1}^{s_2}F(s)\,\mathrm{d}s\leq(s_2-s_1)F(s_1)$.
Hence, assuming $C_1<s_1<s_2\leq 2s_1$,
writing $h:=s_2-s_1$,
and using \eqref{Gsasymptgenform}
and \eqref{equ:assupones},
we have
\begin{align*}
F(s_1)- \frac{a (s_1^{-1}-s_2^{-1})}{h}\geq
-\frac{E(s_1)+E(s_2)}{h}
\geq-\frac{C_2+1}hE(s_1).
\end{align*}
From this, using also 
$s_1^{-1}-s_2^{-1}=h/(s_1s_2)$, 
we get
\begin{align*}
F(s_1)\geq \frac a{s_1s_2}-\frac{C_2+1}h E(s_1)
\geq\frac a{s_1^2}-\frac{ah}{s_1^3}-\frac{C_2+1}h E(s_1).
\end{align*}
Here we optimize by choosing 
$h=\sqrt{E(s_1)s_1^3}$; 
note that this choice of $h$ is admissible, 
i.e.\ yields $s_2\leq2s_1$,
because of our assumption $E(s)\leq s^{-1}$.
The conclusion is that $F(s_1)-as_1^{-2}\geq- (a+C_2+1) s_1^{-\frac32}E(s_1)^{\frac12}$
holds for any $s_1>C_1$.

Similarly, 
using
$G(s_1)-G(s_2)=\int_{s_1}^{s_2}F(s)\,\mathrm{d}s\geq(s_2-s_1)F(s_2)$,
we have for all $C_1<s_1<s_2\leq 2s_1$:
\begin{align*}
F(s_2)\leq\frac a{s_1s_2}+\frac{E(s_1)+E(s_2)}h
\leq\frac a{s_2^2}+\frac{2ah}{s_2^3}+\frac{C_2+1}hE(s_2).
\end{align*}
Here choose $h=\sqrt{4^{-1}E(s_2)s_2^3}$;
then $h\leq\sqrt{4^{-1}s_2^2}=2^{-1}s_2$,
so that our choice is admissible for any $s_2>2C_1$.
For these $s_2$, we obtain
$F(s_2)-as_2^{-2}\leq  (a+2C_2+2)s_2^{-\frac32}E(s_2)^{\frac12}$.
Combining this with the bound proved above, we have proved \eqref{equ:fsasymp} (with an implied constant $C=a+2C_2+2$).
\end{proof}

Let us now use Lemma \ref{prop:gtof}
to conclude the proof of \eqref{equ:gsestgenwin}.
Recall that in the situation of Theorem \ref{MAINTHM1},
$F(s)$ is continuous and decreasing for $s>0$ 
(see Theorem~\ref{gaplimitdistrexistsTHM}),
and $G(s)=\int_s^\infty F(t)\,\text{d}t$.
Furthermore, we have proved the first half of
\eqref{equ:gsestgenwin}, i.e.\ that $G(s)\sim a_{\scrP}s^{-1}$
as $s\to\infty$.
Hence, for any given constant $0<c<1$,
Lemma \ref{prop:gtof} applies 
with $E(s)=cs^{-1}$, $C_2=2$, 
and an appropriate constant $C_1$ (depending on $c$),
to yield $F(s)=a_{\cP}s^{-2}+O(c^{1/2}s^{-2})$ as $s\to\infty$.
Here the implied constant is independent of $c$,
and by taking $c$ arbitrarily small we conclude that
$F(s)=a_{\cP}s^{-2}+o(s^{-2})$.
This completes the proof of 
\eqref{equ:gsestgenwin}.

Similarly, the last part of \eqref{equ:polyest}
follows from the first part of \eqref{equ:polyest},
by Lemma \ref{prop:gtof} applied with
appropriate constants $C_1,C_2$ and 
$E(s)=C s^{-2}\log s$ with an appropriate $C>0$.
This concludes the proof of Theorem \ref{MAINTHM1}.
\hfill$\square$

\begin{remark}\label{niceWslowdecayREM}
Let us note that there exist quite ``nice" (but non-convex)
windows $\scrW$ for which the relative error bound in \eqref{MAINTHM1res1}
in Theorem \ref{MAINTHM1moreexplict}
tends to zero more slowly than any prescribed rate, as $s\to\infty$. 
One way in which this can happen is if $\scrW$ has a \textit{cusp} of an appropriate asymptotic shape.

To give a concrete statement, let us define
\begin{align}\label{niceWslowdecayREMpf2}
H_{\scrW}(s):=\frac 1s\int_0^{C_{\scrP}s}\int_{\R/2\pi\Z}\int_{\R}
\tff\bigl(R_{\scrW}(\theta,y),u^{-1}\bigr)
\,\text{d}y\,\text{d}\theta\,\text{d}u.
\end{align}
Then \eqref{MAINTHM1res1} says that
$G(s)=a_{\scrP} s^{-1}\bigl(1+O(H_{\scrW}(s))\bigr)$ as $s\to\infty$,
viz., $H_\scrW(s)$ is the relative error bound in 
this asymptotics,
and our proof of the relation $G(s)\sim a_{\scrP}s^{-1}$ in \eqref{equ:gsestgenwin}
shows that $H_\scrW(s)\to0$ as $s\to \infty$.
Now let $H_0:[1,\infty)\to(0,\infty)$ be an arbitrary decreasing function
with $\lim_{s\to\infty}H_0(s)=0$.
Choose a sequence $1>t_2>t_3>\cdots$ satisfying
both $\lim_{n\to\infty}t_{n+1}/H_0(n)=\infty$ and $\lim_{n\to\infty}t_n=0$,\footnote{For example,
one may choose $t_{n+1}:=3^{-1}(1+n^{-1})\sqrt{H_0(n)/H_0(1)}$.}
and then let $f:[0,1]\to[0,1]$ be a continuous, increasing function 
satisfying $f(0)=0$, $f(1)=1$,
and $f(t_n)=n^{-1}$ for all $n\geq2$.
(Note that we may even take $f(s)$ to be smooth for $0<s\leq1$.)
Set
\begin{align*}
\scrW:=\scrW_0+\vecw
\quad\text{with}\quad
\scrW_0:=\{(w_1,w_2)\in\R^2\col 0<w_1<1,\:0<w_2<f(w_1)\},
\end{align*}
and with $\vecw$ being a fixed vector in $\R^2$ chosen so that $\bn\in\scrW=\scrW^\circ$.
Finally let $\scrP$ be the cut-and-project set $\scrP=\scrP(\scrW,\scrL_K)$,
for any fixed real quadratic field $K$.
We claim that in this case, $\lim_{s\to\infty}H_{\scrW}(s)/H_0(s)=\infty$,
viz., $H_\scrW(s)$ tends to zero more slowly than the given function $H_0(s)$!

To prove the claim, we first note that 
for every 
$-\frac{\pi}4<\theta<0$, 
writing $w_{\theta,1}$ for the first coordinate of $\vecw \text{k}_{-\theta}$,
one easily verifies that for every $0<y<1/2$,
$R_{\scrW}(\theta,y+w_{\theta,1})$  
is a non-empty open interval of length
$\leq (\cos\theta)^{-1}f(y/\cos\theta)\leq\sqrt 2 f(y/\cos\theta)$.
Hence by Lemma \ref{lem:ffprop}, for each such $\theta$, 
\begin{align*}
\int_{\R}\int_0^{C_{\scrP}s}\tff(R_{\scrW}(\theta,y),u^{-1})\,\text{d}u\,\text{d}y
\gg_{\scrP}
\int_0^{1/2}\min(s,f(y/\cos\theta)^{-1})\,\text{d}y
\gg \int_0^{1/2}\min(s,f(y)^{-1})\,\text{d}y.
\end{align*}
This implies that
\begin{align}\label{niceWslowdecayREMpf1}
H_\scrW(s)
\gg_{\scrP} \: \frac1s\int_0^{1/2}\min(s,f(y)^{-1})\,\text{d}y
=\int_0^{1/2}\min\Bigl(1,\frac1{s f(y)}\Bigr)\,\text{d}y=:H_1(s).
\end{align}
Clearly $H_1(s)$ is a decreasing function of $s>0$.
Also, for each $n\geq2$ we have
$f(y)\leq f(t_n)=n^{-1}$ for all $y\in[0,t_n]$,
and so $H_1(n)\geq\int_0^{t_n}1\,\text{d}y=t_n$.
Hence by our choice of $\{t_n\}$,
we have
$\lim_{n\to\infty} H_1(n+1)/H_0(n)=\infty$.
Note also that for every $n\geq1$ and every $s\in[n,n+1]$,
we have $H_1(s)/H_0(s)\geq H_1(n+1)/H_0(n)$;
hence $\lim_{s\to\infty}H_1(s)/H_0(s)=\infty$.
In view of \eqref{niceWslowdecayREMpf1},
this implies that $\lim_{s\to\infty}H_\scrW(s)/H_0(s)=\infty$.
\hfill$\square$

\vspace{3pt}

As an interesting open problem,
we mention that it seems plausible to us that, by using windows $\scrW$ as above
(perhaps with some further restrictions),
also the \textit{actual} relative difference, $G(s)s-a_{\scrP}$,
can be proved to tend to zero more slowly than any prescribed rate.
\end{remark}

\vspace{3pt}

\begin{remark}\label{WfractalbdryREM}
On the other hand, there exist $\scrW$ with \textit{fractal} boundary for which
\eqref{MAINTHM1res1} in Theorem \ref{MAINTHM1moreexplict} implies as strong an error bound as for convex window sets, 
i.e.\ $G(s)=a_{\scrP}s^{-1}+O\bigl(s^{-2}\log s\bigr)$
as in \eqref{equ:polyest} in Theorem \ref{MAINTHM1}.
For example, this holds if $\partial\scrW$ is the standard Koch snowflake;
and it appears to also hold for the window sets appearing
in, e.g., 
\cite{dDbBzZ2017},
\cite{Niizeki2007}
and for several of the specific examples in
\cite{Niizeki2008}.
Indeed, the proof of Lemma \ref{lem:ffprop} shows that
if the set $A\subset\R$
\textit{contains} an interval of length $b>0$ then
$\int_0^{Cs}\tff(A,u^{-1})\,\text{d}u\ll_C \min(s,b^{-1})$ for all $s>0$.
Hence by the argument leading to \eqref{MAINTHM1pf1}, 
the error bound in \eqref{MAINTHM1res1} is
\begin{align*}
\ll s^{-2}\int_0^{2\pi}\int_0^\infty I\bigl(m(R_{\scrW}(\theta,y))>0\bigr)
\cdot\min\Bigl(s,j(R_{\scrW}(\theta,y))^{-1}\Bigr)\, \text{d}y\,\text{d}\theta,
\end{align*}
where for $A\subset\R$ we write $j(A)$ for the supremum of the lengths of all intervals contained in $A$.
Now temporarily fix $\theta\in[0,2\pi]$ and write $y_0:=\sup\{y>0\col R_{\scrW}(\theta,y)\neq\emptyset\}$.
One verifies that,
if $\scrW$ is an open set (with $\bn\in\scrW$)
such that $\partial\scrW$ is the standard Koch snowflake,
then there exists a constant $c>0$ independent of $\theta$
such that $j(R_{\scrW}(\theta,y))>c\cdot(y_0-y)$ for all $y\in[0,y_0]$.
This leads to $G(s)=a_{\scrP}s^{-1}+O\bigl(s^{-2}\log s\bigr)$, as claimed.
\end{remark}

\subsection{A preliminary integral formula for $G(s)$}
\label{FsexplSEC}
We will now recall 
from \cite{MarklofStrombergsson2015} 
some of the key steps in the proof of Theorem~\ref{gaplimitdistrexistsTHM},
and in particular an
expression for the function $G(s):=\int_s^{\infty}F(t)\, \text{d}t$
as the Haar measure of a certain set in a homogeneous space. 

Let $\scrP$ be as in Theorem~\ref{gaplimitdistrexistsTHM}; thus $\scrP$ is given by
\eqref{scrPWLdef} where $d=2$ and where $\scrL$ is a lattice in $\R^n$ with $n=2+m$.
Set 
$\mathrm{G}=\SL_n(\R)$ and $\Gamma=\SL_n(\Z)$,
and choose $g\in \mathrm{G}$ and $\delta>0$
so that $\scrL=\delta^{1/n}\Z^ng$.
Let $\varphi_g$ be the embedding of $\SL_2(\R)$ in $\mathrm{G}$ given by
\begin{align*}
\varphi_g(A)=g\matr A00{I_m} g^{-1}.
\end{align*}
It follows from Ratner's work 
\cite{Ratner1991,Ratner1991a}
that there exists a unique closed connected subgroup $\mathrm{H}_g$ of $\mathrm{G}$
such that $\Gamma\cap \mathrm{H}_g$ is a lattice in $\mathrm{H}_g$,
$\varphi_g(\SL_2(\R))\subset \mathrm{H}_g$,
and the closure of $\Gamma\bs\Gamma\varphi_g(\SL_2(\R))$
in $\GaG$ equals $\Gamma\bs\Gamma \mathrm{H}_g$.
Let $\mu_g$ be the Haar measure on $\mathrm{H}_g$ normalized so that
$\mu_g(\Gamma\bs\Gamma \mathrm{H}_g)=1$.

Recall that for any $R>0$
we have defined $\Delta_R$ in
\eqref{deltaRDEF}, and 
we have ordered the normalized angles of the points in $\Delta_R$ as
$$
-\tfrac12<\xi_{R,1}\leq\xi_{R,2}\leq\cdots\leq\xi_{R,N(R)}\leq\tfrac12,
$$ 
where 
$N(R)=\#\Delta_R$;
also $\xi_{R,0}:=\xi_{R,N(R)}-1$. 
The proof of 
Theorem~\ref{gaplimitdistrexistsTHM} in \cite{MarklofStrombergsson2015}
starts from the simple geometric observation 
that the events $\xi_{R,j}-\xi_{R,j-1}\geq s/N(R)$ can be detected by rotating the following open sector:
$$
S_R(s):=\left\{(r\cos\theta, r\sin\theta): 0<r<R,\, |\theta|<\frac{\pi s}{N(R)}\right\}
\qquad (s>0).
$$
Indeed, we have
\begin{align}\label{equ:gaprela}
\sum_{j=1}^{N(R)} 
\biggl(\xi_{j,R}-\xi_{j-1,R}-\frac s{N(R)}\biggr)^+
&=\frac{1}{2\pi}\int_0^{2\pi}I(\cP\cap S_R(s) \mathrm{k}_{\theta}=\emptyset)\,\mathrm{d}\theta,
\end{align}
where $x^+:=\max(0,x)$. 
Noticing here that $\cP\cap S_R(s) \mathrm{k}_{\theta}=\emptyset$ holds if and only if
$\cP\mathrm{k}_{-\theta}\mathrm{a}_{1/R}\cap S_R(s)\mathrm{a}_{1/R} =\emptyset$,
and then using \eqref{scrPWLdef},
it follows that the last expression is
\begin{align}
\label{equ:gaprelb}
&=\frac{1}{2\pi}\int_0^{2\pi}I\Bigl(\delta^{1/n}\Z^n\varphi_g(\mathrm{k}_{-\theta}\mathrm{a}_{1/R})g\cap 
(S_R(s)\mathrm{a}_{1/R}\times\scrW) =\emptyset\Bigr)\,\mathrm{d}\theta.
\end{align}
Recalling that $N(R)\sim c_{\scrP}\pi R^2$ as $R\to\infty$ with $c_{\scrP}$ the asymptotic density of $\scrP$, 
one verifies immediately that as $R\to\infty$,
the rescaled sector $S_R(s)\mathrm{a}_{1/R}$
approaches the open triangle $\fC(s)$ in $\R^2$ with vertices at $(0,0)$
and $(1,\pm s/c_{\scrP})$;
\begin{align*}
\fC(s):=\biggl\{(x_1,x_2)\col 0<x_1<1,\: |x_2|<\frac s{c_{\scrP}} x_1\biggr\},
\end{align*}
in the sense that the area of the symmetric difference of
$S_R(s)\mathrm{a}_{1/R}$ and $\fC(s)$ tends to zero.
Also as $R\to\infty$,
the expanding $\operatorname{SO}(2)$-orbit appearing in
\eqref{equ:gaprelb} equidistributes in the homogeneous space
$\G\bk \G\mathrm{H}_g$ \cite[Theorem 4.1]{MarklofStrombergsson2014};
the proof of this fact makes crucial use of
Shah \cite[Thm.\ 1.4]{Shah1996}
and ultimately hinges on Ratner's classification of 
measures invariant under unipotent flows \cite{Ratner1991}.
As is shown in \cite{MarklofStrombergsson2015},
these facts imply
that the expression in \eqref{equ:gaprelb} tends to
\begin{align}\notag
G(s) :&=
\mu_g\bigl(\bigl\{\Gamma h\in\Gamma\bs\Gamma \mathrm{H}_g\col
\delta^{1/n}\Z^nhg\cap(\fC(s)\times\scrW)=\emptyset\bigr\}\bigr)
\\\label{GenGsformula}
&=\mu_g\bigl(\bigl\{\Gamma h\in\Gamma\bs\Gamma \mathrm{H}_g\col \scrP(\scrW,\delta^{1/n}\Z^nhg)\cap\fC(s)=\emptyset\bigr\}\bigr)
\end{align}
as $R\to\infty$.
It is verified in
\cite[Sec.\ 11]{MarklofStrombergsson2015}
that this fact in turn leads to the conclusions of Theorem~\ref{gaplimitdistrexistsTHM},
with the limit distribution $F(s)$ satisfying
\begin{align}\label{FeqmGp}
F(s)=-G'(s)\qquad (s>0).
\end{align}
(In particular the proof in
\cite[Sec.\ 11]{MarklofStrombergsson2015} shows that
$G(s)$ is $C^1$ on $\R_{>0}$, and so $F$ is continuous on $\R_{>0}$ as stated in 
Theorem \ref{gaplimitdistrexistsTHM}.)


\begin{remark}\label{visibleREM2}
The relation \eqref{FeqmGp} corresponds to \cite[(3.6) and (11.1)]{MarklofStrombergsson2015}; however
recall from Remark~\ref{visibleREM}
that our function $F(s)$ 
equals the function in the right hand side of 
\cite[(1.15)]{MarklofStrombergsson2015}.
This is the reason why we have $\fC(s)$ in
\eqref{GenGsformula}, and not $\fC(\kappa_{\scrP}^{-1}s)$
as in \cite[(3.6)]{MarklofStrombergsson2015}.
\end{remark}

\begin{remark}\label{rmk:nullsetinva}
It follows from the formula \eqref{GenGsformula} that $G(s)$ (and hence also $F(s)$) 
remains unchanged if $\cW$ is modified by a null set 
(with respect to the Haar measure $\mu_{\cA}$ of $\cA=\overline{\pi_{\intl}(\scrL)}\subset \R^m$). 
Indeed, let $\cW_1\subset \cA$ be any bounded subset of $\cA$ with boundary of measure zero,
satisfying $\mu_{\cA}(\cW\triangle\cW_1)=0$,
and let $G_1(s)$ be given by 
\eqref{GenGsformula} but with $\cW_1$ in place of $\cW$.
Let $s>0$ be given, and let 
$f_s:\R^d\times\scrA\to\{0,1\}$ be the indicator function 
of the set $\delta^{-\frac1n}\bigl(\fC(s)\times(\cW\triangle\cW_1)\bigr)$.
Then
\begin{align}\notag
|G_1(s)-G(s)|\hspace{390pt}
\\\notag
\hspace{10pt}\leq\mu_g\bigl(\bigl\{\Gamma h\col\text{exactly one of the sets
$\scrP(\scrW,\delta^{1/n}\Z^nhg)$
and $\scrP(\scrW_1,\delta^{1/n}\Z^nhg)$
intersects $\fC(s)$}\bigr\}\bigr)
\\\label{rmk:nullsetinvapf1}
\leq \int_{\G\bk \G \mathrm{H}_g}\sum_{\substack{\vecm\in\Z^nhg}}f_s(\vecm)\,\text{d}\mu_g(h)=0,
\hspace{230pt}
\end{align}
where the final equality holds by 
the Siegel formula 
\cite[Thm.\ 5.1]{MarklofStrombergsson2014},\footnote{Note the correction of this formula given in the erratum to
\cite{MarklofStrombergsson2014}; note also that we may restrict the summation over $\vecm$ in
\eqref{rmk:nullsetinvapf1}
by requiring $\pi(\vecm)\neq\bn$, since $\bn\notin\fC(s)$.}
since $\mu_{\cA}(\cW\triangle\cW_1)=0$.
\end{remark}

Because of the $\mathrm{H}_g$-invariance of $\mu_g$,
the formula \eqref{GenGsformula} remains valid if
in the right hand side we replace
``$\delta^{1/n}\Z^nhg$'' by
``$\delta^{1/n}\Z^nh\varphi_g(A)g$'',
for any given $A\in\SL_2(\R)$;
and applying this with $A=\diag[(s/c_\scrP)^{-1/2},(s/c_\scrP)^{1/2}]$
we obtain:
\begin{align}\label{GenGsformula2}
G(s)=\mu_g\bigl(\bigl\{\Gamma h\in\Gamma\bs\Gamma \mathrm{H}_g\col \scrP(\scrW,\delta^{1/n}\Z^nhg)\cap T(s)=\emptyset\bigr\}\bigr),
\end{align}
where $T(s)\subset \R^2$ is the open triangle with vertices at $(0,0)$ and $(s/c_\scrP)^{1/2}(1,\pm 1)$\label{p:tstriangle}. 

Next we specialize to the setting of Theorem \ref{MAINTHM1moreexplict}.
Thus we take $m=2$ and let $\scrL_K$ be as in 
\eqref{scrLKdef};
this means that $n=4$, $\mathrm{G}=\SL_4(\R)$ and $\Gamma=\SL_4(\Z)$\label{p:sl4z},
and that
$g\in \mathrm{G}$ and $\delta>0$
are such that $\scrL_K=\delta^{1/4}\Z^4g$.
Then by \cite[Sec.\ 2.2.1]{MarklofStrombergsson2014}
we have $\mathrm{H}_g=g\mathrm{H}g^{-1}$\label{p:hg},
where $\mathrm{H}$ is as in Section~\ref{HilbertModGpSEC}.
Furthermore,
\begin{align}\label{GammaKeqGammaconjintH}
\Gamma_K=g^{-1}\Gamma g\cap \mathrm{H}.
\end{align}
Indeed, we have $\Gamma_K\subset \mathrm{H}$,
and every $\gamma\in\Gamma_K$ satisfies
$\scrL_K\gamma=\scrL_K$, hence
$\Z^4g\gamma=\Z^4g$ and so $\gamma\in g^{-1}\Gamma g$.
Conversely,
assume $\gamma\in g^{-1}\Gamma g\cap \mathrm{H}$,
and write $\gamma=\left(\smatr abcd,\smatr{a'}{b'}{c'}{d'}\right)$.
Then $\gamma\in g^{-1}\Gamma g$ implies that
$\Z^4g\gamma g^{-1}=\Z^4$,
and so $\scrL_K \gamma=\scrL_K$.
It follows in particular that $(1,0,1,0)\gamma$ and $(0,1,0,1)\gamma$ lie in $\scrL_K$,
and this in turn implies that $a,b,c,d\in\scrO_K$
and $a'=\sigma(a)$,
$b'=\sigma(b)$,
$c'=\sigma(c)$,
$d'=\sigma(d)$.
Hence $\gamma\in\Gamma_K$, and \eqref{GammaKeqGammaconjintH} is proved.

It follows that the map $\Gamma_K h\mapsto \Gamma ghg^{-1}$ ($h\in \mathrm{H}$)
is a diffeomorphism from $\Gamma\bs\Gamma \mathrm{H}_g$ onto $\Gamma_K\bs \mathrm{H}$,
carrying $\mu_g$ to $\mu_K$.
Hence \eqref{GenGsformula2} can be rewritten as
\begin{align}\label{SpecGsformula}
G(s)=\mu_K\bigl(\bigl\{\Gamma_K h\in\Gamma_K\bs \mathrm{H}\col \scrP(\scrW,\scrL_K h)\cap T(s)=\emptyset\bigr\}\bigr).
\end{align}
Here let us also note that 
$\scrP(\scrW,\scrL_K h)\cap T(s)=\emptyset$
is equivalent with $\scrL_K h\cap\scrT'(s)=\emptyset$, where\label{p:scrt's}
\begin{align*}
\scrT'(s):=T(s)\times\scrW\subset\R^4.
\end{align*}
Hence \eqref{SpecGsformula} can be rewritten as follows:
\begin{align}\label{GsDef}
G(s)&:=\int_{X}I\left(\cL_K h\cap \cT'(s)=\emptyset\right)\, \text{d}\mu_K(h).
\end{align}
This formula will be the starting point of our proof of Theorem \ref{MAINTHM1moreexplict}.
Since we will be only concerned about the asymptotics of $G(s)$ for large s,
in the remainder of this section we will always assume $s > 1$.
\begin{remark}\label{rmk:sl2inva}
Using \eqref{GsDef} and the fact that $(I_2,g)\in \mathrm{H}$ for any $g\in\SL_2(\R)$, 
it is easy to verify that $G(s)$ is unchanged if $\cW$ is replaced by $\cW g$ for any $g\in \SL_2(\R)$.
\end{remark}

\subsection{Separating the main term and error term}\label{sec:smet}

The remainder of Section \ref{sec:asyofgs} is devoted to proving Theorem \ref{MAINTHM1moreexplict}. 
We will start by applying results from Section \ref{sec:geomlavo} to further analyze the condition $\cL_K h\cap \cT'(s)=\emptyset$
in \eqref{GsDef}.
For this we first renormalize the sets $\cT'(s)$ to produce sets containing large balls (for large $s$).
\begin{Lem}\label{renormLEM}
	For $s>1$ let $r\in\N$ be such that $\lambda^{2r}\leq s^{1/4}<\lambda^{2r+2}$, and let \label{p:cts}
\begin{align*}
\cT(s):=(\lambda^{-2r}T(s))\times (\lambda^{2r}\cW).
\end{align*}
Then for any $h\in \mathrm{H}$, $\cL_K h\cap \cT'(s)=\emptyset$ holds if and only if $\cL_K h\cap \cT(s)=\emptyset$.
\end{Lem}
\begin{proof}
	Let $d_{2r}=\diag(\lambda^{-2r}, \lambda^{-2r},\lambda^{2r},\lambda^{2r})$. 
The lemma then follows immediately by noting that $\cT'(s)d_{2r}=\cT(s)$ and $\cL_K hd_{2r}=\cL_K d_{2r}h=\cL_K h$.
	\end{proof}
	Since $\cT(s)$ contains a ball of radius $\gg s^{1/4}$, 
Proposition \ref{prop:bigball} gives that for $s$ sufficiently large, 
the condition $\cL_K h\cap \cT(s)=\emptyset$ forces $\Gamma_K h\in \G_K\cF_i(t_1)$
for some $1\leq i\leq\kappa$,
that is, the point $\Gamma_K h$ in $X$ belongs to the $i$-th cuspidal neighborhood.  
Using also the fact that these cuspidal neighborhoods are pairwise disjoint, by \eqref{equ:cuspdescp}, 
it follows that for $s$ sufficiently large we have
\begin{align*}
G(s)=\sum_{i=1}^{\kappa}G_i(s)
\end{align*}
where for each $1\leq i\leq\kappa$,
\begin{align}\label{GisDEF}
	G_i(s):=\int_{X} I\bigl(\cL_K h\cap \cT(s)=\emptyset\, \text{ and }\, \cL_K h\in \G_K\cF_i(t_1)\bigr)\, \text{d}\mu_K(h).
\end{align}
For the remainder of this section, we fix 
an index $i\in\{1,\ldots,\kappa\}$, and seek an asymptotic formula for the function $G_i(s)$ as $s\to\infty$. 
	
First we note that since $\scrF_i(t_1)$ projects injectively into $X$
(by \eqref{equ:cuspdescp} applied with $j=i$),
we may express $G_i(s)$ as an integral over $\scrF_i(t_1)$.
In view of 
\eqref{equ:haar11}, 
\eqref{equ:cusp1} and \eqref{equ:cusp2},
we get: 
\begin{align}\label{Gisformula1}
G_i(s)=c_K\int_{Z}\int_{Y_{t_1}}\int_{\mathfrak{F}_i}I\left(\cL_K \xi_i\text{n}_{\bm{x}}\text{a}_{\bm{y}}\text{k}_{\bm{\theta}}\cap \cT(s)=\emptyset\right)y_1^{-3}y_2^{-3}\, \text{d}\bm{x}\,\text{d}\bm{y}\,\text{d}\bm{\theta},
	\end{align}
where $c_K=\Delta_K^{-3/2}\zeta_K(2)^{-1}$ as in \eqref{equ:compck},\label{p:setz} 
\begin{align*}
Z:=(0,\pi)\times(0,2\pi),
\end{align*}
and 
\begin{align}\label{Yt1DEF}
Y_{t_1}:=\left\{(y_1,y_2)\in(\R_{>0})^2\col 1\leq y_1/y_2<\lambda^2,\ y_1y_2>t_1\right\}.
\end{align}
Now for $s$ sufficiently large, let $h=\xi_i\text{n}_{\bm{x}}\text{a}_{\bm{y}}\text{k}_{\bm{\theta}}\in \cF_i(t_1)$ be such that $\cL_K h\cap \cT(s)=\emptyset$. Then by Proposition \ref{prop:bigball} we have $y_1,y_2\gg s^{1/4}$.
As discussed in Section \ref{sec:geomlavo}, since $\cL_K h$ avoids a large ball and $h\in \cF_i(t_1)$, it contains a rank two sublattice $\cL_ih$ (cf. \eqref{equ:fillplane2}) which sits densely in the corresponding filled plane $\Pi_i h$ (cf. \eqref{equ:fillplane1}). 
Let us write
\begin{align}\label{tscrLiDEF}
\widetilde{\cL}_i:=\bigcup_{\bm{v}\in \cL_K}(\Pi_i+\bm{v})=\left\{(\alpha, t_1,\sigma(\alpha), t_2)\col \alpha\in I_i,\ t_1, t_2\in \R\right\}\xi^{-1}_i.
\end{align}
Then 
$\cL_K h\subset\widetilde{\cL}_i h$, and $\widetilde{\cL}_i h$ 
is the union of all $\cL_K h$-translates of the filled plane $\Pi_i h$. Here $I_i=\langle a_i,c_i\rangle$ is as before and the second equality 
in \eqref{tscrLiDEF} follows from 
\eqref{equ:lalterdes}
and \eqref{equ:fillplane1}.
Our strategy of computing $G_i(s)$ is to replace the condition $\cL_K h\cap \cT(s)=\emptyset$ by the slightly stronger (since $\widetilde{\cL}_ih$ is very close to $\cL_K h$ when $s$ is sufficiently large) and more manageable condition $\widetilde{\cL}_i h\cap \cT(s)=\emptyset$. 
Thus we define, for $s>1$:
\begin{align}\label{GMisDEF}
G_{M,i}(s):=c_K\int_Z\int_{Y_{t_1}}\int_{\mathfrak{F}_i}I\left(\widetilde{\cL}_ih\cap \cT(s)=\emptyset \right)y_1^{-3}y_2^{-3}\, \text{d}\bm{x}\,\text{d}\bm{y}\,\text{d}\bm{\theta},
\end{align}
and 
\begin{align}\label{GEisDEF}
G_{E,i}(s):=c_K\int_Z\int_{Y_{t_1}}\int_{\mathfrak{F}_i}I\left(\cL_K h\cap \cT(s)=\emptyset,\ \widetilde{\cL}_i h\cap \cT(s)\neq\emptyset\right)y_1^{-3}y_2^{-3}\, \text{d}\bm{x}\,\text{d}\bm{y}\,\text{d}\bm{\theta},
\end{align}
where in both integrals we use the notation
$h=\xi_i\text{n}_{\bm{x}}\text{a}_{\bm{y}}\text{k}_{\bm{\theta}}$. 
It is clear from the definitions that $G_i(s)=G_{M,i}(s)+G_{E,i}(s)$.

We start with a few auxiliary lemmas.
For $\theta\in\R$, let
$\ell(\theta)$ be the projection of $T(1)\text{k}_{-\theta}$ on the $x$-axis\label{p:elltheta}.
Recall also that $\ell_{\scrW}(\theta)$ is the projection of $\cW\text{k}_{-\theta}$ on the $x$-axis.
\begin{Lem}\label{crewheneptLEM}
Fix $s>1$ and let $r\in\N$ be such that $\lambda^{2r}\leq s^{1/4}< \lambda^{2r+2}$ as before. 
For $h=\xi_i\text{n}_{\bm{x}}\text{a}_{\bm{y}}\text{k}_{\bm{\theta}}\in \cF_i(t_1)$, 
the condition
$\widetilde{\cL}_i h\cap \cT(s)=\emptyset$ holds if and only if 
\begin{align}\label{crewheneptLEMres}
\iota(I_i) \cap (y_1^{-1}s^{1/2}\lambda^{-2r}\ell(\theta_1))\times (y_2^{-1}\lambda^{2r}\ell_{\scrW}(\theta_2))=\emptyset,
\end{align}
where $\iota(I_i)\subset \R^2$ is the Minkowski embedding of $I_i$ as before.
\end{Lem}
\begin{proof}
Note that 
\begin{align*}
\widetilde{\cL}_i h=\left\{(\alpha y_1, t_1, \sigma(\alpha)y_2, t_2)\text{k}_{\bm{\theta}} \col \alpha\in I_i,\ t_1, t_2\in \R\right\},
\end{align*}
and recall that $\cT(s)=(\lambda^{-2r}T(s))\times (\lambda^{2r}\cW)$.
Hence $\widetilde{\cL}_i h\cap \cT(s)=\emptyset$ holds if and only if 
\begin{align*}
\forall\, \alpha\in I_i,\, \forall\, (t_1,t_2)\in \R^2 :\quad
(\alpha y_1, t_1, \sigma(\alpha)y_2, t_2)\notin s^{1/2}\lambda^{-2r}T(1)\text{k}_{-\theta_1}\times \lambda^{2r}\cW\text{k}_{-\theta_2},
\end{align*}
which is easily seen to be equivalent to \eqref{crewheneptLEMres}.
\end{proof}

\begin{remark}\label{forbangleREM} 
Similarly, using
$\cL_K=\left\{(\alpha,\beta, \sigma(\alpha), \sigma(\beta))\xi^{-1}_i\col  \alpha\in I_i,\ \beta\in I_i^{-1}\right\}$
(see \eqref{equ:lalterdes}), one verifies that $\cL_K h\cap \cT(s)\neq\emptyset$ 
holds if and only if there exist $\alpha\in I_i$ and $\beta\in I_i^{-1}$ such that
\begin{align}\label{disinclud}
(y_1\alpha, y_1^{-1}(\alpha x_1+\beta), y_2\sigma(\alpha), y_2^{-1}(\sigma(\alpha)x_2+\sigma(\beta)))\in (s^{1/2}\lambda^{-2r}T(1)\text{k}_{-\theta_1})\times (\lambda^{2r}\cW \text{k}_{-\theta_2}).
\end{align} 
\end{remark}

As a consequence of Lemma \ref{crewheneptLEM}
we have a simple necessary condition for $\widetilde{\cL}_i h\cap \cT(s)=\emptyset$.
\begin{Lem}\label{necondLEM}
Keep the assumptions as in Lemma \ref{crewheneptLEM}. 
If $(\theta_1, \theta_2)\in (\frac{\pi}{4}, \frac{3\pi}{4})\times [0, 2\pi)$, then $\widetilde{\cL}_i h\cap \cT(s)\neq \emptyset$.
\end{Lem}
\begin{proof}
Recall that we are assuming $\bn\in\scrW^\circ=\cW$, i.e.\ that $\cW$ contains a small disc centered at the origin;
therefore $0\in\ell_{\scrW}(\theta_2)$ for every $\theta_2$.
Now if $\theta_1\in (\frac{\pi}{4}, \frac{3\pi}{4})$, then $0\in \ell(\theta_1)$,
and hence $\bm{0}\in \iota(I_i)\cap (y_1^{-1}s^{1/2}\lambda^{-2r}\ell(\theta_1))\times (y_2^{-1}\lambda^{2r}\ell_{\scrW}(\theta_2))$,
and so $\widetilde{\cL}_i h\cap \cT(s)\neq\emptyset$ by Lemma \ref{crewheneptLEM}.
\end{proof}

Because of Lemma \ref{necondLEM},
we will often be able to reduce our discussion to the case
$\theta_1\in J$, where\label{p:intervalj}
\begin{align*}
J:=(0,\pi)\setminus \bigl[\tfrac{\pi}{4}, \tfrac{3\pi}{4}\bigr]
=\bigl(0,\tfrac{\pi}{4}\bigr) \cup \bigl(\tfrac{3\pi}{4},\pi\bigr).
\end{align*}
In other words, we will be able to reduce the domain for $\vectheta$ to be\label{p:zj}
\begin{align*}
Z_J:=J\times(0,2\pi)\subset Z.
\end{align*}

\subsection{Computing the main term}
The goal of this section is to obtain simpler formulas for  
the function $G_{M,i}(s)$ for large $s$. We first record the following useful and simple 
sufficient condition for a rectangle in the plane to intersect  
a lattice generated by an ideal of $\cO_K$.
\begin{Lem}\label{diainvaLEM}
For any fractional ideal $I$ of $\cO_K$, there exists a constant $C=C(I)>0$ 
such that for any two intervals $R_1, R_2\subset \R$ with $|R_1|\cdot|R_2|\geq C$,
we have $\iota(I)\cap (R_1\times R_2)\neq\emptyset$.
Here $|R_i|$ denotes the length of the interval $R_i$.
\end{Lem}
\begin{proof}
This is a direct consequence of the fact that $\iota(I)$ is invariant under $\diag(u, \sigma(u))$ for any $u\in \cO_K^{\times}$,
so that we can 
renormalize the rectangle $R_1\times R_2$ to be such that $|R_1|/|R_2|\in [1, \lambda)$.
\end{proof}

We next give a first explicit formula for $G_{M,i}(s)$.
Recall that
$\ell(\theta)$ denotes the projection of $T(1)\text{k}_{-\theta}$ on the $x$-axis.

\begin{Prop}\label{preforgmPROP}
For all sufficiently large $s>1$ we have 
\begin{align}\label{preforgmPROPres}
G_{M,i}(s)=\frac{c_K\Area(\mathfrak{F}_i)}{2s}\int_{Z_J}\int_YI\left(\iota(I_i)\cap (y_1^{-1}\ell(\theta_1)\times y_2^{-1}\ell_{\scrW}(\theta_2))=\emptyset\right)y_1^{-3}y_2^{-3}\, \mathrm{d}\bm{y}\,\mathrm{d}\bm{\theta},
\end{align}
where $Y:=\left\{(y_1,y_2)\in (\R_{>0})^2\col y_2\in [1,\lambda^2)\right\}$\label{p:sety}.
\end{Prop}
\begin{proof}
Using \eqref{GMisDEF}, Lemma \ref{crewheneptLEM} and Lemma \ref{necondLEM},
and the fact that
$\iota(I_i)$ is invariant under the action of $\diag(\lambda^{-2r}, \lambda^{2r})$,
we have
\begin{align}\label{preforgmPROPpf1}
G_{M,i}(s):=c_K\Area(\fF_i)\int_{Z_J}\int_{Y_{t_1}}
I\Bigl(\iota(I_i) \cap (y_1^{-1}s^{1/2}\ell(\theta_1))\times (y_2^{-1}\ell_{\scrW}(\theta_2))=\emptyset\Bigr)
\,y_1^{-3}y_2^{-3}\, \text{d}\bm{y}\,\text{d}\bm{\theta}.
\end{align}
Because of $\bn\in\scrW^\circ=\cW$,
there is an $r_{\cW}>0$ such that $(-r_{\cW},r_{\cW})\subset\ell_{\scrW}(\theta_2)$ for all $\theta_2$. 
Furthermore, $\ell(\theta_1)$ is an open interval of length
$|\ell(\theta_1)|\asymp 1$ for all $\theta_1$. 
Hence by Lemma \ref{diainvaLEM},
there is a constant $C>0$ such that
$\iota(I_i) \cap (y_1^{-1}s^{1/2}\ell(\theta_1))\times (y_2^{-1}\ell_{\scrW}(\theta_2))\neq\emptyset$
holds for all $\vectheta\in Z_J$  
and all $s,y_1,y_2\in\R_{>0}$
subject to $y_1^{-1}s^{1/2}\cdot y_2^{-1}>C$.
It follows that if $s$ is sufficiently large,
then the formula in \eqref{preforgmPROPpf1} remains valid if the
domain of integration $Y_{t_1}$ (defined in \eqref{Yt1DEF})
is replaced by the larger set
\begin{align*}
Y':=\left\{(y_1,y_2)\in (\R_{>0})^2\col y_1/y_2\in [1,\lambda^2)\right\},
\end{align*}
i.e.\ we have (after also changing the order of integration)
\begin{align}\label{preforgmPROPpf2}
G_{M,i}(s):=c_K\Area(\fF_i)\int_{Y'}\int_{Z_J}
I\Bigl(\iota(I_i) \cap (y_1^{-1}s^{1/2}\ell(\theta_1))\times (y_2^{-1}\ell_{\scrW}(\theta_2))=\emptyset\Bigr)
\,\text{d}\vectheta\,\frac{\text{d} y_1}{y_1^3}\,\frac{\text{d} y_2}{y_2^3}.
\end{align}

Going through the same arguments but replacing the condition $y_1/y_2\in [1,\lambda^2)$ by $y_1/y_2\in [\lambda^2,\lambda^4)$ in the definition of $\cF_i$ in \eqref{equ:cusp1} (this amounts to choosing a different fundamental domain for $\G_i\bk \mathrm{H}$),
one verifies that \eqref{preforgmPROPpf2} remains valid for all sufficiently large $s>1$
with $Y'$ replaced by the set
\begin{align*}
\left\{(y_1,y_2)\in (\R_{>0})^2\col y_1/y_2\in [\lambda^2, \lambda^4)\right\}.
\end{align*}
Equivalently, we have for $s>1$ sufficiently large
\begin{align}\label{preforgmPROPpf3}
G_{M,i}(s)&=\frac12c_K\Area(\fF_i)\int_{Y''}\int_{Z_J}
I\Bigl(\iota(I_i) \cap (y_1^{-1}s^{1/2}\ell(\theta_1))\times (y_2^{-1}\ell_{\scrW}(\theta_2))=\emptyset\Bigr)
\,\text{d}\vectheta\,\frac{\text{d} y_1}{y_1^3}\,\frac{\text{d} y_2}{y_2^3},
\end{align}
where 
\begin{align*}
Y'':=\left\{(y_1,y_2)\in (\R_{>0})^2\col y_1/y_2\in [1, \lambda^4)\right\}.
\end{align*}

Next, using the fact that the lattice $\iota(I_i)$ is invariant under the action of
$S=\diag(\lambda^2,\lambda^{-2})$, it follows that the integrand in \eqref{preforgmPROPpf3} 
is invariant under $\vecy\mapsto\vecy S$, for each fixed $\vectheta$. Also the measure $y_1^{-3}y_2^{-3}\,\text{d}\vecy$
is invariant under this map and both $Y''$ and $Y$ are fundamental domains for $(\R_{>0})^2/\langle\vecy\mapsto \vecy S\rangle$. Hence the formula 
\eqref{preforgmPROPpf3} remains valid if the domain of integration $Y''$ is replaced by $Y$, giving us \eqref{preforgmPROPres}. 
\end{proof}

Next we further simplify the above empty intersection condition to get an even more explicit formula in terms of the function $\alpha_{I,J}$ defined in \eqref{alphapWiDEF}.

\begin{Prop}\label{moreexplicitPROP}
For all sufficiently large $s>1$ we have 
\begin{align}\label{newGMi3}
G_{M,i}(s)=\frac{c_K\Area(\mathfrak{F}_i)\,c_{\scrP}}{4s}\int_0^{2\pi}\int_1^{\lambda^2}
\alpha_{I_i,\ell_{\cW}(\theta)}(y)^2\,\frac{\mathrm{d}y}{y^3}\, \mathrm{d}\theta.
\end{align}
In particular, we have for all sufficiently large $s>1$,
\begin{align}\label{moreexplicitPROPres2}
\sum_{i=1}^{\kappa}G_{M,i}(s)=\frac{a_{\cP}}{s},
\end{align}
where $a_{\cP}$ is as in \eqref{MAINTHM1res2}.
\end{Prop}
\begin{proof}
The formula in Proposition \ref{preforgmPROP} can be expressed as
\begin{align}\label{newGMi2}
G_{M,i}(s)=\frac{c_K\Area(\mathfrak{F}_i)}{2s}\int_0^{2\pi}\int_1^{\lambda^2}
\scrI(y_2,\theta_2)\,\frac{\text{d}y_2}{y_2^3}\,\text{d}\theta_2,
\end{align}
where
\begin{align*}
\scrI(y_2,\theta_2):=\int_J\int_0^\infty I\left(\iota(I_i)\cap (y_1^{-1}\ell(\theta_1)\times y_2^{-1}\ell_{\scrW}(\theta_2))=\emptyset\right)
\,\frac{\text{d}y_1}{y_1^3}\,\text{d}\theta_1.
\end{align*}
Here recall that $\ell(\theta_1)$ equals the projection of $T(1)\text{k}_{-\theta_1}$ on the $x$-axis,
which we can compute to be:
\begin{align}\label{moreexplicitPROPpf1}
\ell(\theta_1)=\begin{cases}
\bigl(0,c_{\scrP}^{-1/2}(\cos\theta_1+\sin\theta_1)\bigr)
&\text{if }\:0\leq\theta_1\leq\pi/4,
\\
\bigl(c_{\scrP}^{-1/2}(\cos\theta_1-\sin\theta_1),0\bigr)
&\text{if }\:3\pi/4\leq\theta_1\leq\pi.
\end{cases}
\end{align}
Recall the definition of $\alpha_{I,J}(y)$
in
\eqref{alphapWiDEF}. We also give a companion definition here. For any integral ideal $I\subset \cO_K$ and any bounded Lebesgue measurable set $J\subset \R$ with non-empty interior, 
we set for any $y>0$,
\begin{align}\label{moreexplicitPROPpf2}
\tilde\alpha_{I,J}(y):=\max\bigl\{\alpha\in I\cap\R_{<0}\col y\,\sigma(\alpha)\in J\bigr\}.
\end{align}
With these definitions the formula for $\scrI(y_2,\theta_2)$ can be rewritten as
\begin{align*}
&\int_0^{\pi/4}\int_0^\infty I\left(\alpha_{I_i,\ell_{\cW}(\theta_2)}(y_2)\notin y_1^{-1}\ell(\theta_1)\right)\,\frac{\text{d}y_1}{y_1^3}\,\text{d}\theta_1
+\int_{3\pi/4}^{\pi}\int_0^\infty I\left(\tilde\alpha_{I_i,\ell_{\cW}(\theta_2)}(y_2)\notin y_1^{-1}\ell(\theta_1)\right)\,\frac{\text{d}y_1}{y_1^3}\,\text{d}\theta_1
\\
&=\frac12c_{\scrP}\alpha_{I_i,\ell_{\cW}(\theta_2)}(y_2)^2\int_0^{\pi/4}\frac{\text{d}\theta_1}{(\cos\theta_1+\sin\theta_1)^2}
+\frac12c_{\scrP}\tilde\alpha_{I_i,\ell_{\cW}(\theta_2)}(y_2)^2\int_{3\pi/4}^{\pi}\frac{\text{d}\theta_1}{(\cos\theta_1-\sin\theta_1)^2}
\\
&=\frac{c_{\scrP}}{4}\bigl(\alpha_{I_i,\ell_{\cW}(\theta_2)}(y_2)^2+\tilde\alpha_{I_i,\ell_{\cW}(\theta_2)}(y_2)^2\bigr),
\end{align*}
implying that
\begin{align*}
\int_0^{2\pi}\int_1^{\lambda^2}
\scrI(y_2,\theta_2)\,\frac{\text{d}y_2}{y_2^3}\,\text{d}\theta_2&=\frac{c_{\cP}}4\int_0^{2\pi}\int_1^{\lambda^2}
\bigl(\alpha_{I_i,\ell_{\cW}(\theta)}(y)^2+\tilde\alpha_{I_i,\ell_{\cW}(\theta)}(y)^2\bigr)\,\frac{\text{d}y}{y^3}\, \text{d}\theta.
\end{align*}
Now note that $\ell_{\cW}(\theta+\pi)=-\ell_{\cW}(\theta)$ from which we can deduce that $\alpha_{I_i,\ell_{\cW}(\theta+\pi)}(y)=-\tilde\alpha_{I_i,\ell_{\cW}(\theta)}(y)$. Hence
\begin{align*}
\int_0^{2\pi}\int_1^{\lambda^2}\tilde\alpha_{I_i,\ell_{\cW}(\theta)}(y)^2\,\frac{\text{d}y}{y^3}\,\text{d}\theta&=\int_0^{2\pi}\int_1^{\lambda^2}
\alpha_{I_i,\ell_{\cW}(\theta+\pi)}(y)^2\,\frac{\text{d}y}{y^3}\, \text{d}\theta\\
&=\int_0^{2\pi}\int_1^{\lambda^2}
\alpha_{I_i,\ell_{\cW}(\theta)}(y)^2\,\frac{\text{d}y}{y^3}\, \text{d}\theta,
\end{align*}
implying that 
\begin{align*}
\int_0^{2\pi}\int_1^{\lambda^2}
\scrI(y_2,\theta_2)\,\frac{\text{d}y_2}{y_2^3}\,\text{d}\theta_2&=\frac{c_{\cP}}2\int_0^{2\pi}\int_1^{\lambda^2}
\alpha_{I_i,\ell_{\cW}(\theta)}(y)^2\,\frac{\text{d}y}{y^3}\, \text{d}\theta.
\end{align*}
Plugging this back to \eqref{newGMi2} we get the desired formula for $G_{M,i}(s)$, \eqref{newGMi3}.

The formula \eqref{moreexplicitPROPres2} 
follows directly from \eqref{newGMi3} by using
$$
\Area(\mathfrak{F}_i)=\covol(\iota(I^{-2}_i))=
\text{Nr}(I_i)^{-2}\Delta_K^{1/2}
$$ 
(which holds by \eqref{equ:covolfor}) and since $\mathfrak{F}_i$ is a fundamental domain for $\R^2/\iota(I_i^{-2})$),
and the formulas for $c_K$ and $c_{\cP}$ respectively; see \eqref{equ:compck} and \eqref{cPformula}.
\end{proof}

\subsection{Estimating the error term}\label{sec:estimerror}
Recall that $G(s)=\sum_{i=1}^{\kappa}G_i(s)=\sum_{i=1}^{\kappa} \bigl(G_{M,i}(s)+G_{E,i}(s)\bigr)$
for all sufficiently large $s$,
and by definition we have $G_{E,i}(s)\geq0$; see \eqref{GEisDEF}.
Recalling also Proposition~\ref{moreexplicitPROP}, it follows that,
in order to complete the proof of Theorem \ref{MAINTHM1moreexplict},
it only remains
to prove that $G_{E,i}(s)$ is majorized by the bound in
\eqref{MAINTHM1res1} for all sufficiently large $s$.

We first prove two auxiliary lemmas.
\begin{Lem}\label{lem:diainvaVAR}
For any fractional ideal $I$ of $\cO_K$, there exists a constant $C=C(I)>0$
such that for any two intervals $R_1, R_2\subset \R$ 
with $|R_1|\cdot|R_2|\geq C$,
and for any Lebesgue measurable subset $R_2'\subset \R$, we have
\begin{align*}
\Area(X_I\setminus\pi(R_1\times R_2'))\leq |R_1|\cdot m(R_2\setminus R_2'),
\end{align*}
where $X_I:=\R^2/\iota(I)$, and $\pi:\R^2\to X_I$ is the quotient map.
\end{Lem}
(Recall that $m$ denotes Lebesgue measure on $\R$; also for an interval $R$ we write $|R|=m(R)$ for its length.)

\begin{proof}
Let us fix a fundamental parallelogram $\mathfrak{F}$  
for the lattice $\iota(I)$ in $\R^2$,
and then fix intervals $J_1,J_2\subset\R$ such that $\mathfrak{F}\subset J_1\times J_2$.
We claim that the statement of the lemma holds with
$C:=\lambda |J_1||J_2|$.
Indeed, assume that $R_1, R_2$ are intervals with $|R_1|\cdot|R_2|\geq\lambda |J_1||J_2|$.
Then there exists some $m\in\Z$ such that
$\lambda^m|R_1|\geq|J_1|$ and $\lambda^{-m}|R_2|\geq|J_2|$,
viz., $|\sigma(\lambda^m)R_2|\geq|J_2|$.
This means that a translate of $(R_1\times R_2)\diag(\lambda,\sigma(\lambda))^m$
contains $F$;
and hence $R_1\times R_2$ contains a translate $\mathfrak{F}'$ of
$\mathfrak{F} \diag(\lambda,\sigma(\lambda))^{-m}$.
But using the fact that $\iota(I)$ is invariant under 
$\diag(\lambda,\sigma(\lambda))$,
it follows that $\mathfrak{F}'$ is again an fundamental parallelogram for the lattice $\iota(I)$ in $\R^2$.
Hence
\begin{align*}
\Area(X_I\setminus\pi(R_1\times R_2'))
=\Area(\mathfrak{F}'\setminus\pi^{-1}(\pi(R_1\times R_2'))).
\end{align*}
Using $\mathfrak{F}'\subset R_1\times R_2$ and $\pi^{-1}(\pi(R_1\times R_2'))\supset R_1\times R_2'$,
the above is
\begin{align*}
\leq\Area((R_1\times R_2)\setminus(R_1\times R_2'))
=\Area(R_1\times(R_2\setminus R_2'))
= |R_1|\cdot m(R_2\setminus R_2'). 
\end{align*} 
\end{proof}

\begin{Lem}\label{lem:esincfa}
For any fixed set $A\subset\R$ we have
\begin{align}\label{ffessincr}
\forall\: 0<a_1\leq a_2:
\qquad
\ff(A,a_1)\leq 2\ff(A,a_2).
\end{align}
\end{Lem}
\begin{proof}
Let $0<a_1\leq a_2$ be given,
and let $J$ be an arbitrary interval of length $a_2$.
Take $k\in\Z_{>0}$ so that $ka_1\leq a_2<(k+1)a_1$;
then there exist $k$ pairwise disjoint (open) intervals 
$J_1,\ldots,J_k\subset J$ which all have length $a_1$.
Now $m(J\setminus A)\geq\sum_{\ell=1}^km(J_\ell\setminus A)\geq ka_1\ff(A,a_1)$,
and since this holds for all intervals $J$ of length $a_2$,
it follows that $\ff(A,a_2)\geq a_2^{-1}\cdot ka_1\ff(A,a_1)\geq\frac k{k+1}\ff(A,a_1)\geq\frac12\ff(A,a_1)$.
\end{proof}

We now turn to bounding the error term
$G_{E,i}(s)$, which we defined in \eqref{GEisDEF},
for large $s$.

First of all, note that if $\theta_1\in(\frac{\pi}{4}, \frac{3\pi}{4})$
then $T(1)\text{k}_{-\theta_1}$ contains the line segment between
$(0,0)$ and $(0,c_{\scrP}^{-1/2})$.
Recall also that since $\bn\in\scrW^\circ=\cW$, there is an 
$r_{\cW}>0$ such that $\scrB_{r_{\cW}}^2\subset\scrW$.
It follows that if $\theta_1\in(\frac{\pi}{4}, \frac{3\pi}{4})$
then the condition \eqref{disinclud} is satisfied with $\alpha=0$ whenever
$(\beta,\sigma(\beta))$ belongs to the rectangle
$\bigl(0,y_1s^{1/2}\lambda^{-2r}c_{\scrP}^{-1/2}\bigr)\times\bigl(-y_2\lambda^{2r}r_{\cW},y_2\lambda^{2r}r_{\cW}\bigr)$.
This rectangle has area $\gg s^{1/2}$, since $\vecy\in Y_{t_1}$.
Hence by Lemma \ref{diainvaLEM},
for $s$ sufficiently large there always exists some $\beta\in I_i^{-1}$
such that $(\beta,\sigma(\beta))$ belongs to the rectangle,
and so, by Remark \ref{forbangleREM}, $\scrL_K h\cap\scrT(s)\neq\emptyset$.
This proves that the contribution from all $\theta_1\in(\frac{\pi}{4}, \frac{3\pi}{4})$
to the integral in \eqref{GEisDEF} is zero.
Thus, for $s$ sufficiently large, we may replace the range of integration
for $\vectheta$ in \eqref{GEisDEF} 
by $Z_J$.

Next, as in the discussion in Section \ref{sec:smet}
(near \eqref{Yt1DEF}),
because of the condition $\cL_K h\cap \cT(s)=\emptyset$ in the integrand,
any $\vecy\in Y_{t_1}$ which contribute to the
integral in \eqref{GEisDEF}
must satisfy $y_1,y_2\gg s^{1/4}$.
Furthermore, by Lemma \ref{crewheneptLEM},
given any $\langle\vectheta,\vecy,\vecx\rangle$ which makes the integrand in \eqref{GEisDEF} nonzero,
there exists some $\alpha\in I_i$ such that 
$(\alpha,\sigma(\alpha))$ belongs to the set 
\begin{align}\label{GEisbounddisc2}
\fA_{\vectheta,\vecy}(s):=(y_1^{-1}s^{1/2}\lambda^{-2r}\ell(\theta_1))\times (y_2^{-1}\lambda^{2r}\ell_{\scrW}(\theta_2)).
\end{align}
Since we have restricted to $\theta_1\in J$,
we have $0\notin\ell(\theta_1)$, and hence $\alpha\neq0$.
Note that $\ell(\theta_1)$ and $\ell_{\scrW}(\theta_2)$
are contained in bounded intervals independent of $\theta_1,\theta_2$;
hence for $y_1,y_2\gg s^{1/4}$, the set $\fA_{\vectheta,\vecy}(s)$
is contained in a fixed bounded region $B\subset\R^2$.
Set
\begin{align*}
\Lambda:=\{\alpha\in I_i\setminus\{0\}\col (\alpha,\sigma(\alpha))\in B\}.
\end{align*}
This is a finite subset of $I_i\setminus\{0\}$,
and our discussion shows that for any 
$\langle\vectheta,\vecy,\vecx\rangle\in Z_J\times Y_{t_1}\times \fF_i$ 
which makes the integrand in \eqref{GEisDEF} nonzero,
there must exist some $\alpha\in\Lambda$ satisfying 
$(\alpha,\sigma(\alpha))\in\fA_{\vectheta,\vecy}(s)$,
or equivalently, $\iota(\Lambda)\cap\fA_{\vectheta,\vecy}(s)\neq\emptyset$.
Using now the fact that both the minima of $|\alpha|$ and $|\sigma(\alpha)|$
for $\alpha\in\Lambda$ are bounded away from zero,
and also the fact that both $|\ell(\theta_1)|$ and $|\ell_{\scrW}(\theta_2)|$
are bounded away from zero independently of $\theta_1,\theta_2$,
it follows that we must have $y_1,y_2\ll s^{1/4}$.
We have thus proved that there exist some constants $c_2>c_1>0$
(which depend on $K$ and $\scrW$, but not on $s$),
such that for all sufficiently large $s$,
every $\langle\vectheta,\vecy,\vecx\rangle\in Z_J\times Y_{t_1}\times \fF_i$ 
which makes the integrand in \eqref{GEisDEF} nonzero,
must satisfy $c_1s^{1/4}\leq y_1,y_2\leq c_2 s^{1/4}$.

In view of the above discussion, we have for $s$ sufficiently large:
\begin{align}\notag
G_{E,i}(s)\asymp s^{-3/2}\int_{Z_J}\int_{c_1s^{1/4}}^{c_2s^{1/4}}\int_{c_1s^{1/4}}^{c_2s^{1/4}}
I\Bigl(\iota(\Lambda)\cap\fA_{\vectheta,\vecy}(s)\neq\emptyset\Bigr) 
\int_{\mathfrak{F}_i}I\left(\cL_K h\cap \cT(s)=\emptyset\right)\, \text{d}\bm{x}\,\text{d}\bm{y}\,\text{d}\bm{\theta}
\\\label{GEisbounddisc3}
\asymp s^{-3/2}\sum_{\alpha\in\Lambda}\int_{Z_J}\int_{c_1s^{1/4}}^{c_2s^{1/4}}\int_{c_1s^{1/4}}^{c_2s^{1/4}}
I\Bigl((\alpha,\sigma(\alpha))\in\fA_{\vectheta,\vecy}(s)\Bigr)
\int_{\mathfrak{F}_i}I\left(\cL_K h\cap \cT(s)=\emptyset\right)\, \text{d}\bm{x}\,\text{d}\bm{y}\,\text{d}\bm{\theta}.
\end{align}

Next, the condition $\cL_K h\cap \cT(s)=\emptyset$ in 
\eqref{GEisbounddisc3} implies by Remark \ref{forbangleREM}
that \eqref{disinclud} fails for all $(\alpha,\beta)\in I_i\times I_i^{-1}$.
This is equivalent to saying that for every $\alpha\in I_i$
we have $\iota(I_i^{-1})\cap\fB_{\vectheta,\vecy,\vecx,\alpha}(s)=\emptyset$,
where
\begin{align}\label{condfinEQU}
\fB_{\vectheta,\vecy,\vecx,\alpha}(s):=\left(y_1\mathfrak{J}_1-\alpha x_1\right)\times \left(y_2\mathfrak{J}_2-\sigma(\alpha)x_2\right) 
\end{align}
with $\mathfrak{J}_1=\mathfrak{J}_1(s,y_1,\theta_1,\alpha)$ and $\mathfrak{J}_2=\mathfrak{J}_2(s,y_2,\theta_2,\alpha)$ defined by
\begin{align}\label{J1J2def}
\mathfrak{J}_1:=\left\{t\in \R\col (y_1\alpha,t)\in s^{1/2}\lambda^{-2r}T(1)\text{k}_{-\theta_1}\right\}
\quad\text{and}\quad
\mathfrak{J}_2:=\left\{t\in\R\col (y_2\sigma(\alpha),t)\in \lambda^{2r}\cW \text{k}_{-\theta_2}\right\}.
\end{align}
Hence from \eqref{GEisbounddisc3} we get,
using also the fact that if 
$\iota(I_i^{-1})\cap\fB_{\vectheta,\vecy,\vecx,\alpha}(s)=\emptyset$
holds for every $\alpha\in I_i$,
then in particular it holds for 
the $\alpha\in\Lambda$ which is our summation variable:
\begin{align}\label{GEisbounddisc5}
G_{E,i}(s)\ll s^{-3/2}\sum_{\alpha\in\Lambda}\int_{Z_J}\int_{c_1s^{1/4}}^{c_2s^{1/4}}\int_{c_1s^{1/4}}^{c_2s^{1/4}}
I\Bigl((\alpha,\sigma(\alpha))\in\fA_{\vectheta,\vecy}(s)\Bigr)
\hspace{100pt}
\\\notag
\times
\int_{\mathfrak{F}_i}
I\Bigl(\iota(I_i^{-1})\cap \fB_{\vectheta,\vecy,\vecx,\alpha}(s)=\emptyset\Bigr)
\, \text{d}\bm{x}\,\text{d}\bm{y}\,\text{d}\bm{\theta}.
\end{align}
Note here that 
$\fB_{\vectheta,\vecy,\vecx,\alpha}(s)=\fB_{\vectheta,\vecy,\bn,\alpha}(s)-(\alpha x_1,\sigma(\alpha)x_2)$,
and recall that $\fF_i$ is a fundamental domain for $\R^2/\iota(I_i^{-2})$.
Note also that the map $\vecx\mapsto \vecv:=(\alpha x_1,\sigma(\alpha)x_2)$
induces a diffeomorphism from $\R^2/\iota(I_i^{-2})$
onto $\R^2/\iota(\alpha I_i^{-2})$,
carrying $\text{d}\vecx$
to $|\mathrm{N}(\alpha)|^{-1}\text{d}\vecv$.
Hence the integral over $\vecx$ in 
\eqref{GEisbounddisc5} equals
\begin{align}\label{GEisbounddisc5a}
\frac1{|\mathrm{N}(\alpha)|}\int_{\R^2/\iota(\alpha I_i^{-2})}
I\Bigl(\iota(I_i^{-1})\cap (\fB_{\vectheta,\vecy,\bn,\alpha}(s)-\vecv)=\emptyset\Bigr)
\, \text{d}\vecv.
\end{align}
Let $\pi$ be the quotient map from $\R^2$ onto the torus $X:=\R^2/\iota(I_i^{-1})$.
Note that the condition
$\iota(I_i^{-1})\cap (\fB_{\vectheta,\vecy,\bn,\alpha}(s)-\vecv)=\emptyset$
is equivalent with 
$\pi(\vecv)\notin\pi(\fB_{\vectheta,\vecy,\bn,\alpha}(s))$.
Also, since $\alpha\in \Lambda\subset I_i$,
$\iota(\alpha I_i^{-2})$ is a sublattice of
$\iota(I_i^{-1})$; 
hence $\pi$ induces to a covering map 
from
$\R^2/\iota(\alpha I_i^{-2})$ onto $X$,
which preserves Lebesgue measure.
Using \eqref{equ:covolfor} it follows that this covering map has degree $|\mathrm{N}(\alpha)|/\text{Nr}(I_i)$.
Hence the integral in \eqref{GEisbounddisc5a} equals
\begin{align*}
\Nr(I_i)^{-1}\Area(X\setminus\pi(\fB_{\vectheta,\vecy,\bn,\alpha}(s))).
\end{align*}
But we have $\fB_{\vectheta,\vecy,\bn,\alpha}(s)=y_1\mathfrak{J}_1\times y_2\mathfrak{J}_2$,
where $\mathfrak{J}_1$ is an open interval which is non-empty 
if and only if $\alpha\in y_1^{-1}s^{1/2}\lambda^{-2r}\ell(\theta_1)$.
This last condition 
is guaranteed to hold because of the condition $(\alpha,\sigma(\alpha))\in\fA_{\vectheta,\vecy}(s)$
appearing in \eqref{GEisbounddisc5}.
It now follows from Lemma \ref{lem:diainvaVAR}
that there is a constant $C=C(I_i^{-1})$ such that,
for any interval $R_2$ of length $|R_2|=C/|y_1\mathfrak{J}_1|$, 
the integral in \eqref{GEisbounddisc5a} is bounded from above by
$\Nr(I_i)^{-1}|y_1\mathfrak{J}_1|\cdot m\bigl(R_2\setminus y_2\mathfrak{J}_2\bigr)$.
Furthermore, comparing \eqref{J1J2def} and \eqref{RWdef} we see that
$\mathfrak{J}_2=\lambda^{2r} R_{\scrW}(\theta_2,\lambda^{-2r}\sigma(\alpha)y_2)$.
Hence, writing $R_2=y_2\lambda^{2r}\cdot(x,x+a)$ with $x\in\R$
and $a=C/\bigl(\lambda^{2r}y_1y_2|\mathfrak{J}_1|\bigr)$,
it follows that the integral in \eqref{GEisbounddisc5a} is bounded from above by
\begin{align*}
\Nr(I_i)^{-1}\lambda^{2r}y_1y_2|\mathfrak{J}_1|\cdot m\bigl((x,x+a)\setminus R_{\scrW}(\theta_2,\lambda^{-2r}\sigma(\alpha)y_2)\bigr).
\end{align*}
This is true for every $x\in\R$.
Hence, using the notation introduced in \eqref{ffAaDEF},
we conclude that the integral in \eqref{GEisbounddisc5a} is bounded from above by
\begin{align*}
\Nr(I_i)^{-1} C\cdot
\ff\biggl(R_{\scrW}(\theta_2,\lambda^{-2r}\sigma(\alpha)y_2),\frac{C}{\lambda^{2r}y_1y_2|\mathfrak{J}_1|}\biggr).
\end{align*}
Note also that 
$(\alpha,\sigma(\alpha))\in\fA_{\vectheta,\vecy}(s)$ holds if and only if
both $\alpha\in y_1^{-1}s^{1/2}\lambda^{-2r}\ell(\theta_1)$
and $\sigma(\alpha)\in y_2^{-1}\lambda^{2r}\ell_{\scrW}(\theta_2)$,
and we recall that the first of these conditions is equivalent with $\mathfrak{J}_1\neq\emptyset$,
while, by \eqref{ellWdef}, the second condition is equivalent with
$R_{\scrW}(\theta_2,\lambda^{-2r}y_2\sigma(\alpha))\neq\emptyset$.
Using these facts together with our bound on the integral in \eqref{GEisbounddisc5a},
we obtain:
\begin{align}\notag
G_{E,i}(s)\ll s^{-3/2}\sum_{\alpha\in\Lambda}
\int_{Z_J}
\int_{c_1s^{1/4}}^{c_2s^{1/4}}
I(\mathfrak{J}_1\neq\emptyset)
\int_{c_1s^{1/4}}^{c_2s^{1/4}}
I\bigl(R_{\scrW}(\theta_2,\lambda^{-2r}y_2\sigma(\alpha))\neq \emptyset\bigr)
\hspace{80pt}
\\\label{GEisbounddisc7}
\times\ff\biggl(R_{\scrW}(\theta_2,\lambda^{-2r}\sigma(\alpha)y_2),\frac{C}{\lambda^{2r}y_1y_2|\mathfrak{J}_1|}\biggr)
\, \text{d}y_2\,\text{d}y_1\,\text{d}\bm{\theta}.
\end{align}
Here it should be recalled that $\mathfrak{J}_1=\mathfrak{J}_1(s,y_1,\theta_1,\alpha)$; see \eqref{J1J2def}.
Now note that for any $y_1,y_2$ appearing in the above integral,
we have $\lambda^{2r}y_1y_2\geq c_1^2\lambda^{-2} s^{3/4}$;
hence if we set $C':=C\lambda^2 c_1^{-2}$ then
${\displaystyle \frac{C}{\lambda^{2r}y_1y_2|J_1|}\leq\frac{C'}{s^{3/4}|J_1|}}$,
and so, using Lemma \ref{lem:esincfa}
and then substituting $y_2=\lambda^{2r}\sigma(\alpha)^{-1}y$,
we obtain (since $\Lambda$ is a fixed finite subset of $I_i\setminus\{0\}$):
\begin{align}\label{GEisbounddisc8}
G_{E,i}(s)
\ll s^{-5/4}\sum_{\alpha\in\Lambda}
\int_{Z_J}
\int_{c_1s^{1/4}}^{c_2s^{1/4}}
I(\mathfrak{J}_1\neq\emptyset)
\int_{\R}
\tff\biggl(R_{\scrW}(\theta_2,y),\frac{C'}{s^{3/4}|\mathfrak{J}_1|}\biggr)
\,\text{d}y\,\text{d}y_1\,\text{d}\bm{\theta},
\end{align}
where we use the notation $\tff$ defined in \eqref{ffAaDEF}.

Now we have
\begin{align*}
\mathfrak{J}_1\neq\emptyset\quad \Leftrightarrow \quad
\left\lbrace\begin{array}{ll} 
0<y_1\alpha<s^{\frac12}\lambda^{-2r}c_{\scrP}^{-\frac12}(\cos\theta_1+\sin\theta_1) & \text{if $\theta_1\in (0, \tfrac{\pi}{4})$},\\
s^{\frac12}\lambda^{-2r}c_{\scrP}^{-\frac12}(\cos\theta_1-\sin\theta_1)<y_1\alpha< 0& \text{if $\theta_1\in (\tfrac{3\pi}{4}, \pi)$},
\end{array}\right.
\end{align*}
and it follows from the definition in \eqref{J1J2def} that
for any $\alpha\in\Lambda$, $\theta_1\in J$ and $y_1\geq c_1s^{1/4}$ 
for which $\mathfrak{J}_1\neq\emptyset$ holds,
we have $|\mathfrak{J}_1|\geq C'' u_1$ with
$u_1:=s^{\frac12}\lambda^{-2r}c_{\scrP}^{-\frac12}(|\cos\theta_1|+\sin\theta_1)-y_1|\alpha|\in\R_{>0}$,
where $C''>0$ is a certain absolute constant.
Furthermore, $u_1<c_3 s^{1/4}$ 
with $c_3:=2c_{\scrP}^{-\frac12}\lambda^2$.
Hence, using Lemma \ref{lem:esincfa},
\begin{align}\notag
G_{E,i}(s)
&\ll s^{-5/4}\sum_{\alpha\in\Lambda}
\int_{Z_J}
\int_0^{c_3s^{1/4}}
\int_{\R}
\tff\biggl(R_{\scrW}(\theta_2,y),\frac{C'}{C''s^{3/4}u_1}\biggr)
\,\text{d}y\,\text{d}u_1\,\text{d}\bm{\theta}.
\end{align}
Here the integrand is independent of $\theta_1$ and $\alpha$.
Hence, using the fact that $\Lambda$ is a fixed finite set,
and substituting $u_1=(C'/C'') s^{-3/4}u$,
we conclude:
\begin{align}\label{GEibound15}
G_{E,i}(s)\ll s^{-2}\int_0^{(c_3C''/C')s}\int_{0}^{2\pi}\int_{\R}
\tff\bigl(R_{\scrW}(\theta_2,y),u^{-1}\bigr)
\,\text{d}y\,\text{d}\theta_2\,\text{d}u.
\end{align}
Hence we have proved the error bound in \eqref{MAINTHM1res1}
(with $C:=c_3C''/C'$).
As we have already noted, 
because of $G(s)=\sum_{i=1}^{\kappa}G_i(s)=\sum_{i=1}^{\kappa} \bigl(G_{M,i}(s)+G_{E,i}(s)\bigr)$
and Proposition \ref{moreexplicitPROP},
the bound \eqref{GEibound15}
completes the proof of Theorem \ref{MAINTHM1moreexplict}.
\hfill$\square$

\section{More explicit formula for the leading coefficient}\label{acPmoreexplicitSEC}
The main goal of this section is to prove the following more explicit formula for the leading coefficient $a_{\cP}$ in Theorem \ref{MAINTHM1},
when we further assume that 
$\scrW$ is such that $\ell_{\scrW}(\theta)$
is an interval for every $\theta$.
Recall here that for any $\theta\in \R/2\pi \Z$, $\ell_{\cW}(\theta)$ is the projection of $\cW\text{k}_{-\theta}$ on the $x$-axis.

\begin{thm}\label{thm:moreexfor}
Keep the notation and assumptions as in Theorem \ref{MAINTHM1} and further assume that 
$\ell_{\cW}(\theta)$ is an interval for each $\theta$,
parametrized as in \eqref{equ:ellw},
i.e.\ $\ell_{\cW}(\theta)=r(\theta)(-\nu(\theta), 1)$.
Then for each $I_i$ ($1\leq i\leq \kappa$), there exist a finite partition $\R_{>0}=\bigsqcup_{j=1}^{l_{I_i}}S_{I_i,j}$ of $\R_{>0}$ into intervals,
and non-negative constants $A_{I_i,j}, B_{I_i,j}$ $(1\leq j\leq l_{I_i})$, such that
\begin{align}\label{equ:genforap}
a_{\scrP}=\frac{\Area(\scrW)}{4 \Delta_K^2\zeta_K(2)}\sum_{i=1}^{\kappa}\operatorname{Nr}(I_i)^{-2}\sum_{j=1}^{l_{I_i}}\int_{\tilde{S}_{I_i,j}}r(\theta)^{-2}\left(A_{I_i,j}+B_{I_i,j}\nu(\theta)^{-2}\right)\,\mathrm{d}\theta,
\end{align}
where $\tilde{S}_{I_i,j}:=\{\theta\in [0,2\pi)\col \nu(\theta)\in S_{I_i,j}\}$.
\end{thm}

\begin{remark}\label{thm:moreexforREM}
For each $I_i$, the partition $\R_{>0}=\bigsqcup_{j=1}^{l_{I_i}}S_{I_i,j}$ and the constants $A_{I_i,j}, B_{I_i,j}$ are all computable from our analysis; see Section \ref{sec:exmapksqrt23} for three explicit examples.
We stress that the intervals $S_{I_i,j}$ are allowed to be open, closed, or half-open,
and may be degenerate, i.e.\ of the form $[a,a]=\{a\}$ for some $a>0$.
\end{remark}

\subsection{Relations to the polar set of $\cW$}
\label{UNDintegralSEC}
Let $\cW\subset \R^2$ be as above.  
The integral appearing in \eqref{equ:genforap} is closely related to the polar set $\cW^*$ of $\cW$,
which is defined by \eqref{def:polbody}.
Note that $\cW^*$ is closed and convex;
furthermore, since $\cW$ is bounded and contains $\bn$ in its interior,  
$\cW^*$ is also bounded and contains $\bn$ in its interior. 
Moreover, if $\scrW$ is symmetric about the origin, then so is $\scrW^*$;
and if $\scrW$ is a convex polygon, then also $\scrW^*$ is a convex polygon.
 
The following lemma shows that $\cW^*$ can be parameterized in polar coordinates 
in terms of the function $r(\theta)$.
\begin{Lem}\label{lem:polarbody}
Let $\scrW$ be as in Theorem \ref{thm:moreexfor}.
Then $\cW^*$ is parameterized in polar coordinates by
\begin{align*}
\cW^*=\left\{(t,\theta)\in \R_{\geq 0}\times \R/2\pi \Z\col 0\leq t\leq r(-\theta)^{-1}\right\}.
\end{align*}
\end{Lem}
\begin{proof}
Note that for any $\theta\in \R/2\pi \Z$, 
\begin{align*}
\cW\text{k}_{-\theta}=\left\{(w_1\cos\theta-w_2\sin\theta, w_1\sin\theta+w_2\cos\theta)\col (w_1,w_2)\in \cW\right\}.
\end{align*}
Thus
\begin{align*}
\ell_{\cW}(\theta)=\left\{w_1\cos\theta-w_2\sin\theta\col (w_1,w_2)\in \cW\right\}.
\end{align*}
It then follows from the definition of polar set and the relation $\ell_{\cW}(\theta)\cap \R_{\geq 0}=[0, r(\theta))$, that
the set
\begin{align*}
\left\{t\in\R_{\geq0}\col t(\cos\theta,-\sin\theta)\in\scrW^*\right\}
\end{align*}
equals the interval $[0,r(\theta)^{-1}]$. 
The lemma then follows by a substitution $\theta\mapsto -\theta$. 
\end{proof}

\begin{remark}
As a direct consequence of Lemma \ref{lem:polarbody} and the relation 
\begin{align}\label{equ:rela1}
r(\theta+\pi)=r(\theta)\nu(\theta),\quad \forall\,\theta\in \R/2\pi \Z,
\end{align}
(which in turn follows from the symmetry $\ell_{\cW}(\theta+\pi)=-\ell_{\cW}(\theta)$),
we have
\begin{align}\label{INTell2inv2formula}
2\Area(\scrW^*)=\int_0^{2\pi}r(\theta)^{-2}\,\text{d}\theta=\int_0^{2\pi}\nu(\theta)^{-2}r(\theta)^{-2}\,\text{d}\theta.
\end{align}
For each $1\leq i\leq \kappa$, let $1\leq j_i\leq l_{I_i}$ be the unique index such that $1\in S_{I_i,j_i}$.
In the special case when $\tilde{S}_{I_i,j_i}=[0,2\pi)$ for every $1\leq i\leq \kappa$ 
(this happens if $\cW$ is symmetric with respect to the origin, 
or a sufficiently small translate of such a set),  then by \eqref{INTell2inv2formula} the formula \eqref{equ:genforap} simplifies to
\begin{align}\label{equ:censymfor}
a_{\cP}&= 
\frac{\Area(\scrW) \Area(\scrW^*)}{\Delta_K^2\,\zeta_K(2)}
\sum_{i=1}^{\kappa}C_{I_i}\text{Nr}(I_i)^{-2},
\end{align}
with $C_{I_i}:=\frac12(A_{I_i,j_i}+B_{I_i,j_i})$.
More generally, the sum of integrals appearing in \eqref{equ:genforap} can be understood as a \textit{weighted} area of $\cW^*$. 
\end{remark}

The remainder of this section is devoted to proving Theorem \ref{thm:moreexfor}.

\subsection{Some preliminary computations}
We start our proof with the formula for $a_{\cP}$ given in Theorem \ref{MAINTHM1moreexplict} which states that $a_{\cP}=\sum_{i=1}^{\kappa}a_{\cP,i}$ with 
\begin{align}\label{equ:expsmper2}
a_{\cP,i}&=\frac{\Area(\scrW)}{4 \Delta_K^{2}\zeta_K(2)}\text{Nr}(I_i)^{-2}\int_{0}^{2\pi}\int_1^{\lambda^{2}}
\alpha_{I_i,\ell_{\cW}(\theta)}(y)^2\,\frac{\text{d}y}{y^3}\, \text{d}\theta.
\end{align}
It thus suffices to compute the above double integral which we denote by $\cI_{\cP,i}$. 
It is immediate from the definition \eqref{alphapWiDEF} that 
$\alpha_{I,aJ}(ay)=\alpha_{I,J}(y)$ for all $a,y>0$.
Thus, defining 
$J_{\nu}:=(-\nu,1)$\label{p:jnu}
for $\nu>0$ and applying \eqref{idclassinvREM2}, we have:
\begin{align}\label{cIPDEF}
\cI_{\cP,i}&=\int_{0}^{2\pi}\int_1^{\lambda^{2}}
\alpha_{I_i,J_{\nu(\theta)}}\Bigl(\frac{y}{r(\theta)}\Bigr)^2\,\frac{\text{d}y}{y^3}\, \text{d}\theta
=\int_{0}^{2\pi}r(\theta)^{-2}\int_1^{\lambda^{2}}
\alpha_{I_i,J_{\nu(\theta)}}(y)^2\,\frac{\text{d}y}{y^3}\, \text{d}\theta.
\end{align}
We wish to compute the innermost integral in the last expression.

\subsection{Extremal points}\label{sec:extpts}
Fix an integral ideal $I\subset \cO_K$. As mentioned above, we are interested in computing the integral 
\begin{align}\label{equ:desinte}
\int_1^{\lambda^2}\alpha_{I,J_{\nu}}(y)^2\,\frac{\text{d}y}{y^3}
\end{align}
for any $\nu>0$. For this we introduce a notion which captures the possible values of the function $\alpha_{I,J_{\nu}}$. For any $\nu>0$ and $\alpha\in I\cap\R_{>0}$, let $b_{\alpha,\nu}:=\sigma(\alpha)$ if $\sigma(\alpha)>0$ and $b_{\alpha,\nu}:=\nu^{-1}\sigma(\alpha)$ if $\sigma(\alpha)<0$\label{p:balpha1}. Note that for any $y>0$, $y\sigma(\alpha)\in J_{\nu}$ if and only if $y<|b_{\alpha,\nu}|^{-1}$. 
We say $\alpha\in I\cap \R_{>0}$ is \textit{$\nu$-extremal with respect to $I$} 
if $$
\iota(I)\cap R_{\alpha,\nu}=\emptyset,
$$ 
where\label{p:ralphanu} 
\begin{align*}
R_{\alpha,\nu}:=\left\{\begin{array}{ll}   
	(0,\alpha)\times [-\nu\sigma(\alpha), \sigma(\alpha)] & \text{if $\sigma(\alpha)>0$},\\
	(0, \alpha)\times [\sigma(\alpha), \nu^{-1}|\sigma(\alpha)|] & \text{if $\sigma(\alpha)<0$}.
	\end{array}\right.
\end{align*}
We denote by $E_{I,\nu}\subset I\cap \R_{>0}$ 
the set of all $\nu$-extremal points with respect to $I$\label{p:einu}. Below we will simply call elements in $E_{I,\nu}$ $\nu$-extremal when there is no ambiguity.
Since $\iota(I)$ is a lattice in $\R^2$, $E_{I,\nu}$ is always nonempty. Indeed, $E_{I,\nu}$ 
is always infinite, as it is invariant under a certain multiplication map as shown in the following lemma.
\begin{Lem}\label{lem:invar}
The set $E_{I,\nu}$ is invariant under the map $S: \R_{>0}\to \R_{>0}$ defined by $S(y):= \lambda y$ if $\sigma(\lambda)>0$ and $S(y):= \lambda^2 y$ if $\sigma(\lambda)<0$.
\end{Lem}
\begin{proof}
Suppose that $\alpha\in I\cap \R_{>0}$ is $\nu$-extremal, i.e. $\iota(I)\cap R_{\alpha,\nu}=\emptyset$ with $R_{\alpha,\nu}$ the rectangle given as above. Let $d_{\lambda}=\diag(\lambda, \sigma(\lambda))$ and note that $\iota(I)$ is invariant under the right multiplication action of $d_{\lambda}$. The lemma then follows by noting that $R_{\alpha,\nu}d_{\lambda}=R_{\lambda\alpha,\nu}$ if $\sigma(\lambda)>0$ and
$R_{\alpha,\nu}d^2_{\lambda}=R_{\lambda^2\alpha,\nu}$ if $\sigma(\lambda)<0$.
\end{proof}

\begin{remark}
The set $E_{I,1}$ is invariant under the map $y\mapsto \lambda y$ even if $\sigma(\lambda)<0$. 
This is because we always have $R_{\alpha,1}d_{\lambda}=R_{\lambda\alpha,1}$,
since $R_{\alpha,1}=(0,\alpha)\times [-|\sigma(\alpha)|, |\sigma(\alpha)|]$.
\end{remark}
Note that $\alpha\in I\cap\R_{>0}$ is $\nu$-extremal if and only if
 \begin{align}\label{equ:alterndefnuex}
 \text{$|b_{\beta,\nu}|>|b_{\alpha,\nu}|$ for any $\beta\in I\cap\R_{>0}$ with $\beta<\alpha$.}
 \end{align}
In particular, if $\alpha<\beta$ are two $\nu$-extremal points, then $|b_{\alpha,\nu}|>|b_{\beta,\nu}|$. Moreover, one can check from \eqref{equ:alterndefnuex} that $\alpha_{I,J_{\nu}}$ takes values only on $\nu$-extremal points. 
Using these facts,  
we can now prove the following lemma which enables us to compute the integral \eqref{equ:desinte}.
\begin{Lem}\label{lem:relwithin}
For any $\nu>0$, let $0<\beta<\alpha$ be two consecutive elements in $E_{I,\nu}$. Then 
$\alpha_{I,J_{\nu}}(y)=\alpha$
for any $|b_{\beta,\nu}|^{-1}\leq y<|b_{\alpha,\nu}|^{-1}$.
\end{Lem}
\begin{remark}\label{rmk:discreteness}
By discreteness of $\iota(I)$ one sees that $\#(E_{I,\nu}/S)<\infty$ for any $\nu>0$; thus we can take two consecutive elements from $E_{I,\nu}$ as in the above lemma.
\end{remark}
\begin{proof}[Proof of Lemma \ref*{lem:relwithin}]
Take any $|b_{\beta,\nu}|^{-1}\leq y<|b_{\alpha,\nu}|^{-1}$ and let $\beta'=\alpha_{I,J_{\nu}}(y)$. In particular, $\beta'$ is $\nu$-extremal and $y\sigma(\beta')\in J_{\nu}$, and the last condition is equivalent to $y<|b_{\beta',\nu}|^{-1}$.
Moreover, $y<|b_{\alpha,\nu}|^{-1}$ implies that $y\sigma(\alpha)\in J_{\nu}$. Thus $\beta'\leq \alpha$. 
If $\beta'<\alpha$, then since $\beta<\alpha$ are two consecutive $\nu$-extremal points, 
we must have $\beta'\leq \beta$ and therefore $|b_{\beta',\nu}|\geq |b_{\beta,\nu}|$, 
which contradicts the fact that
$|b_{\beta,\nu}|^{-1}\leq y<|b_{\beta',\nu}|^{-1}$.
Hence $\beta'=\alpha$ as desired.
\end{proof}

We can now compute the integral $\int_1^{\lambda^2}\alpha_{I,J_{\nu}}(y)^2\ \frac{\text{d}y}{y^3}$ in terms of the set $E_{I,\nu}$. 

\begin{Lem}\label{lem:alphaintegral}
For any $\nu>0$ and $\alpha_0\in E_{I,\nu}$,
if $\alpha_0<\alpha_1<\cdots< \alpha_l=\lambda^2\alpha_0$ is the complete list of points in 
$E_{I,\nu}\cap[\alpha_0,\lambda^2\alpha_0]$, then
\begin{align*}
\int_1^{\lambda^2}\alpha_{I,J_{\nu}}(y)^2\,\frac{\mathrm{d}y}{y^3}
&=A_{I,\nu}+B_{I,\nu}\nu^{-2},
\end{align*}
where 
$$
A_{I,\nu}:=\frac12\sum_{\substack{j=0\\(\sigma(\alpha_j)>0)}}^{l-1}\sigma(\alpha_j)^2\bigl(\alpha_{j+1}^2-\alpha_j^2\bigr)
\quad \text{and} \quad B_{I,\nu}:=\frac12\sum_{\substack{j=0\\(\sigma(\alpha_j)<0)}}^{l-1}\sigma(\alpha_j)^2\bigl(\alpha_{j+1}^2-\alpha_j^2\bigr).
$$
\end{Lem}

\begin{proof}
We have $|b_{\alpha_0,\nu}|^{-1}<|b_{\alpha_1,\nu}|^{-1}<\cdots<|b_{\alpha_l,\nu}|^{-1}=|b_{\alpha_0,\nu}|^{-1}\lambda^2$.
Hence by \eqref{idclassinvREM2} and Lemma \ref{lem:relwithin},
\begin{align*}
\int_1^{\lambda^2}\alpha_{I,J_{\nu}}(y)^2\,\frac{\text{d}y}{y^3}&=\sum_{j=1}^{l}\int_{|b_{\alpha_{j-1},\nu}|^{-1}}^{|b_{\alpha_j,\nu}|^{-1}}\alpha_j^2\, \frac{\text{d}y}{y^3}=\frac12\sum_{j=1}^{l}\alpha_j^2\left({b^{2}_{\alpha_{j-1},\nu}}-{b^{2}_{\alpha_j,\nu}}\right).
\end{align*}
The lemma then follows by rewriting the above as (using the fact that
$\alpha_l^2b_{\alpha_l,\nu}^2=\alpha_0^2 b_{\alpha_0,\nu}^2$):
\begin{align*}
\frac12\sum_{j=0}^{l-1}b_{\alpha_j,\nu}^2\bigl(\alpha_{j+1}^2-\alpha_j^2\bigr)
=\frac12\sum_{\substack{j=0\\(\sigma(\alpha_j)>0)}}^{l-1}\sigma(\alpha_j)^2\bigl(\alpha_{j+1}^2-\alpha_j^2\bigr)
+\frac{\nu^{-2}}2\sum_{\substack{j=0\\(\sigma(\alpha_j)<0)}}^{l-1}\sigma(\alpha_j)^2\bigl(\alpha_{j+1}^2-\alpha_j^2\bigr).
\end{align*}
\end{proof}
\begin{remark}
If $\alpha<\beta$ are two consecutive $\nu$-extremal points, then so are $\lambda^{2n}\alpha<\lambda^{2n}\beta$ for any integer $n$. Using the equality $\sigma({\alpha})^2(\beta^2-\alpha^2)=\sigma({\lambda^{2n}\alpha})^2\left((\lambda^{2n}\beta)^2-(\lambda^{2n}\alpha)^2\right)$ one easily sees that both $A_{I,\nu}$ and $B_{I,\nu}$ are independent of the choice of $\alpha_0$.
\end{remark}

Next, we study how $E_{I,\nu}$ varies in $\nu$. Our goal is to show there are only finitely many possibilities for $E_{I,\nu}$ (and hence also for the coefficients $A_{I,\nu}$ and $B_{I,\nu}$). For this we introduce another related notion. We say $\alpha\in I\cap \R_{>0}$ is \textit{positively} (resp. \textit{negatively}) \textit{extremal} with respect to $I$ if $\sigma(\alpha)>0$ (resp. $\sigma(\alpha)<0$) and $\iota(I)\cap (0,\alpha)\times [0,\sigma(\alpha)]=\emptyset$ (resp. $\iota(I)\cap (0,\alpha)\times [\sigma(\alpha), 0]=\emptyset$). We denote by $E^+_{I}$ and $E_I^-$ the set of positively and negatively extremal points in $I$\label{p:eipm}. 
Recall that $S$ is the map $y\mapsto \lambda y$ if $\sigma(\lambda)>0$ and the map $y\mapsto \lambda^2 y$ if $\sigma(\lambda)<0$. It is clear that $E_{I}^{\pm}$ are invariant under $S$, and from the discreteness of $\iota(I)$ we see that $\#(E_I^{\pm}/S)<\infty$. 
Note that a necessary condition for $\alpha\in I\cap \R_{>0}$ to be $\nu$-extremal for some $\nu>0$ is that $\alpha\in E_I^+\sqcup E_I^-$. Indeed this is also a sufficient condition in the sense that for any $\alpha\in E_I^+\sqcup E_I^-$ there exists some $\nu>0$ such that $\alpha\in E_{I,\nu}$; see Lemma \ref{lem:vainnu} below.

For any $\nu>0$, set $E_{I,\nu}^{\pm}:=E_{I,\nu}\cap E_{I}^{\pm}$
so that $E_{I,\nu}=E^+_{I,\nu}\sqcup E_{I,\nu}^-$. Note that $E_{I,\nu}^+$ (resp. $E_{I,\nu}^-$) is decreasing (resp.\ increasing) in $\nu$,
and there is a constant $C=C_I>0$
such that $E_{I,\nu}^+=\emptyset$ for all $\nu>C$
and $E_{I,\nu}^-=\emptyset$ for all $\nu<1/C$,
because of Lemma \ref{diainvaLEM} and the fact that  
$$
\Area(R_{\alpha,\nu})=(1+\nu^{\sgn(\sigma(\alpha))})|\alpha\sigma(\alpha)|>\nu^{\sgn(\sigma(\alpha))}
\qquad (\forall\, \alpha\in I\setminus\{0\}).
$$ 
To understand how $E_{I,\nu}$ varies in $\nu$ we study how $E_{I,\nu}^{\pm}$ varies in $\nu$ which is easier as it is monotone. 
We first study this property for a \textit{single} element in $I$.
The following lemma is immediate from our definitions.
\begin{Lem}\label{lem:semicont}
For any $\alpha\in E_I^+$ let 
$$
\nu_{\alpha}:=\sup\{\nu>0\col \iota(I)\cap (0,\alpha)\times [-\nu\sigma(\alpha),\sigma(\alpha)]=\emptyset\}.
$$
Then $\alpha$ is $\nu$-extremal if and only if $0<\nu<\nu_{\alpha}$. Similarly, for any $\alpha\in E_I^-$ define
$$
\nu_{\alpha}:=\inf\{\nu>0\col \iota(I)\cap (0,\alpha)\times [\sigma(\alpha),\nu^{-1}|\sigma(\alpha)|]=\emptyset\}.
$$  
Then $\alpha$ is $\nu$-extremal if and only if $\nu>\nu_{\alpha}$. 
\end{Lem}
 
\begin{remark}\label{rmk:resonva}
For each $\alpha\in E_I^+\sqcup E_I^-$, when $\nu$ attains the critical value $\nu_{\alpha}$, $\alpha$ \textit{fails} to be $\nu_{\alpha}$-extremal. This fact, together with the finiteness of $(E_I^+\sqcup E_I^-)/S$ and the fact that $E_{I,\nu}^+$ (resp. $E_{I,\nu}^-$) is decreasing (resp. increasing) in $\nu$ implies that $E_{I,\nu}^+$ (resp. $E_I^-$) is right (resp. left) continuous in the sense that for any $\nu>0$, there exists some $\e>0$ such that $E_{I,\nu'}^+=E_{I,\nu}^+$ for any $\nu'\in (\nu,\nu+\e)$ (resp. $\nu'\in (\nu-\e,\nu)$). 
Note also that 
\begin{align*}
\nu_{\alpha}=\left\{\begin{array}{ll}   
	\inf\left\{\tfrac{-\sigma(\beta)}{\sigma(\alpha)}\col \beta\in I,\, 0<\beta<\alpha,\, \sigma(\beta)<0 \right\} & \text{if $\alpha\in E_I^+$},\\
	\sup\left\{\tfrac{-\sigma(\alpha)}{\sigma(\beta)}\col \beta\in I, 0<\beta<\alpha,\, \sigma(\beta)>0 \right\} & \text{if $\alpha\in E_I^-$}.
	\end{array}\right.
\end{align*}
From this formula we see that the $\beta$ attaining the above infimum (resp. supremum) is an element in $E_I^-$ (resp. $E_I^+$). 
This gives very strict restrictions on the possible values of these $\nu_{\alpha}$'s. In particular, $\nu_{\alpha}$ can not be $1$ since otherwise we would have $\beta=-\alpha$, contradicting $0<\beta<\alpha$. 
Moreover, the set of values $\{\nu_\alpha\col \alpha\in E_{I}^+\sqcup E_I^-\}$
is finite, since $\#(E_{I}^+\sqcup E_I^-)/S<\infty$ and 
$\nu_{\lambda^2\alpha}=\nu_{\alpha}$ for every $\alpha\in E_{I}^+\sqcup E_I^-$.
\end{remark}

\begin{Lem}\label{lem:vainnu}
There exist finitely many constants $0<\nu_1^+<\cdots<\nu_l^+<\infty$ such that $E_I^+\supsetneq E_{I,\nu_1^+}^+\supsetneq E_{I,\nu_2^+}^+\supsetneq \cdots\supsetneq E_{I,\nu_l^+}^+=\emptyset$ and
\begin{align*}
E_{I,\nu}^+&=\left\{\begin{array}{ll}   
	E_I^{+} & \text{if $0<\nu<\nu_1^+$},\\
	E^+_{I,\nu_i^+} & \text{if $\nu_i^+\leq \nu<\nu_{i+1}^+$\quad ($1\leq i< l$)},\\
	\emptyset & \text{if $\nu\geq \nu_l^+$}.
	\end{array}\right.
\end{align*}
Similarly, there exist $\infty>\nu_1^->\cdots>\nu_{l'}^->0$ such that $E_I^-\supsetneq E_{I,\nu_1^-}^-\supsetneq E_{I,\nu_2^-}^-\supsetneq \cdots\supsetneq E_{I,\nu_{l'}^-}^-=\emptyset$ and
\begin{align*}
E_{I,\nu}^-&=\left\{\begin{array}{ll}   
	E_I^{-} & \text{if $\nu>\nu_1^-$},\\
	E^-_{I,\nu_i^-} & \text{if $\nu_{i+1}^-< \nu\leq \nu_{i}^-$\quad ($1\leq i< l'$)},\\
	\emptyset & \text{if $0<\nu\leq \nu_{l'}^-$}.
	\end{array}\right.
\end{align*}
\end{Lem}
\begin{proof}
We noted in Remark \ref{rmk:resonva} that the set
$\{\nu_\alpha\col \alpha\in E_I^+\}$ is finite;
let us order its elements as 
$0<\nu_1^+<\cdots<\nu_l^+$.
Then if $0<\nu<\nu_1^+$,
it follows from Lemma \ref{lem:semicont} that
every $\alpha\in E_I^+$ is $\nu$-extremal, i.e.\ $E_{I,\nu}^+=E_I^+$;
similarly, if $\nu\geq\nu_l^+$ then $E_{I,\nu}^+=\emptyset$.
Finally if $\nu_i^+\leq \nu<\nu_{i+1}^+$ for some $1\leq i< l$,
then for every $\alpha\in E_I^+$ the condition
$\nu<\nu_\alpha$ is equivalent with $\nu_\alpha\in\{\nu_{i+1}^+,\ldots,\nu_l^+\}$;
hence by Lemma \ref{lem:semicont},
$E_{I,\nu}^+=\{\alpha\in E_I^+\col \nu_\alpha\in\{\nu_{i+1}^+,\ldots,\nu_l^+\}\}=E_{I,\nu_i^+}^+$.
We have thus proved the part of the lemma which concerns $E_{I,\nu}^+$;
the proof of the part concerning $E_{I,\nu}^-$ is entirely similar.
\end{proof}

\begin{remark}\label{rmk:construt}
From the proof we see  
that the break-points in the statement of 
Lemma \ref{lem:vainnu} are exactly the values
$\nu_\alpha$ with $\alpha\in E_I^+\sqcup E_I^-$.
Explicitly, 
$\{\nu_i^+\col 1\leq i\leq l\}=\{\nu_{\alpha}\col \alpha\in E_I^+\}$ and $\{\nu_i^-\col 1\leq i\leq l'\}=\{\nu_{\alpha}\col \alpha\in E_I^-\}$.
\end{remark}
As a corollary we have the following.
\begin{Cor}\label{cor:fulldes}
There exist a finite partition of $\R_{>0}=\sqcup_{j=1}^{l_I}S_{I,j}$ into intervals, 
and $l_I$ pairwise distinct subsets 
$B_1=E_I^+, B_2,\ldots, B_{l_I-1}, B_{l_I}=E_I^-$ of $E_I^+\sqcup E_I^-$, such that $E_{I,\nu}=B_j$ for any $\nu\in S_{I,j}$.
\end{Cor}
(The intervals $S_{I,j}$ are allowed to be open, closed, half-open, and degenerate; see Remark \ref{thm:moreexforREM}.)
\begin{proof}
The existence of such a partition follows from Lemma \ref{lem:vainnu} and the relation $E_{I,\nu}=E_{I,\nu}^+\sqcup E_{I,\nu}^-$. The pairwise disjointness of $B_1,\ldots,B_{l_I}\subset E_I^+\sqcup E_I^-$ follows again from the relation $E_{I,\nu}=E_{I,\nu}^+\sqcup E_{I,\nu}^-$ and also the fact that $E_{I,\nu}^+$ and $E_{I,\nu}^-$ are decreasing and increasing in $\nu$ respectively.
\end{proof}

\begin{remark}\label{rmk:disonone}
We see from the above proof that if $\nu\notin \{\nu_{\alpha}\col \alpha\in E_I^+\sqcup E_I^-\}$, then $\nu$ is an interior point of one of the intervals $S_{I,j}$.  
In particular, since $1\notin \{\nu_{\alpha}\col \alpha\in E_I^+\sqcup E_I^-\}$, 
one of the intervals $S_{I,j}$ contains $1$ as an interior point.
\end{remark}
We can now give the proof of Theorem \ref{thm:moreexfor}.
\begin{proof}
For each $I_i$ ($1\leq i\leq \kappa$), let $\R_{>0}=\sqcup_{j=1}^{l_{I_i}}S_{I_i,j}$ be the partition as in Corollary \ref{cor:fulldes} 
applied with $I=I_i$, 
and define $A_{I_i,\nu}\geq0$ and $B_{I_i,\nu}\geq0$ for $\nu>0$ as in Lemma \ref{lem:alphaintegral}.
For each $1\leq j\leq l_{I_i}$, it follows from Corollary~\ref{cor:fulldes} that 
both $A_{I_i,\nu}$ and $B_{I_i,\nu}$ are constant as $\nu$ varies over $S_{I_i,j}$;
thus we may define $A_{I_i,j}:=A_{I_i,\nu}$ and $B_{I_i,j}:=B_{I_i,\nu}$ for any $\nu\in S_{I_i,j}$.
Then by \eqref{cIPDEF} and Lemma \ref{lem:alphaintegral} we have
\begin{align*}
\cI_{\cP,i}&=\sum_{j=1}^{l_{I_i}}\int_{\tilde{S}_{I_i,j}}r(\theta)^{-2}\left(A_{I_i,j}+B_{I_i,j}\nu(\theta)^{-2}\right)\,\frac{\text{d}y}{y^3}\text{d}\theta,
\end{align*}
where $\tilde{S}_{I_i,j}:=\{\theta\in [0,2\pi)\col \nu(\theta)\in S_{I_i,j}\}$.
Plugging this formula into \eqref{equ:expsmper2} and applying the relation $a_{\cP}=\sum_{i=1}^{\kappa}a_{\cP,i}$ we get the desired formula for $a_{\cP}$.
\end{proof}

\subsection{Examples; proofs of Propositions \ref*{ABaPexpl}--\ref*{TTaPexpl}}\label{sec:exmapksqrt23}
In this section we work out the formula \eqref{equ:genforap} explicitly
in the three cases $K=\Q(\sqrt{2})$, $\Q(\sqrt{3})$ and $\Q(\sqrt{5})$, thus proving Propositions \ref{ABaPexpl}--\ref{TTaPexpl}.
Note that in all cases,
$K$ is of class number one;
thus $\kappa=1$ and $I_1=\cO_K$,
and we set $I:=I_1$.
For details of the following computations, see also the program
\cite[explicit\_aP.mpl]{HSYsupplement}.

First we treat the case when $K=\Q(\sqrt{2})$. 
Then $\lambda=1+\sqrt{2}$, and one verifies that
\begin{align*}
E_I^+=\{1,2+\sqrt2\}\cdot\lambda^{2\Z}
\quad\text{and}\quad
E_I^-=\{\sqrt2,1+\sqrt2\}\cdot\lambda^{2\Z},
\end{align*}
with the associated values $\nu_\alpha$
given by 
\begin{align*}
\nu_1=1+\sqrt2,
\qquad
\nu_{2+\sqrt2}=\tfrac12\sqrt2
\qquad \text{and}\qquad
\nu_{\sqrt2}=\sqrt2,
\quad 
\nu_{1+\sqrt2}=\sqrt2-1.
\end{align*}
Hence 
\begin{align*}
E_{I,\nu}^+=\begin{cases}
\{1,2+\sqrt2\}\cdot\lambda^{2\Z}&\text{if }\: \nu<\tfrac12\sqrt2,
\\
\lambda^{2\Z}&\text{if }\: \tfrac12\sqrt2\leq\nu<1+\sqrt2,
\\
\emptyset&\text{if }\: 1+\sqrt2\leq\nu,
\end{cases}
\end{align*}
and
\begin{align*}
E_{I,\nu}^-=\begin{cases}
\emptyset&\text{if }\: \nu\leq\sqrt2-1,
\\
(1+\sqrt2)\lambda^{2\Z}&\text{if }\: \sqrt2-1<\nu\leq\sqrt2,
\\
\{\sqrt2,1+\sqrt2\}\cdot\lambda^{2\Z}&\text{if }\: \sqrt2<\nu.
\end{cases}
\end{align*}
Using also the relation $E_{I,\nu}=E_{I,\nu}^+\sqcup E_{I,\nu}^-$, 
we get the corresponding partition $\R_{>0}=\sqcup_{j=1}^5S_{I,j}$ where $S_{I,j}=S_j$ 
with $S_{1}=(0,\sqrt{2}-1], S_{2}=(\sqrt{2}-1, \frac{1}{2}\sqrt{2}), S_{3}=[\frac12\sqrt{2}, \sqrt{2}], S_{4}=(\sqrt{2}, 1+\sqrt{2})$ and $S_{5}=[1+\sqrt{2},\infty)$.
The corresponding constants $A_j:=A_{I,j}$ and $B_j:=B_{I,j}$
are 
\begin{align*}
(A_1, B_1)=(B_5, A_5)=(\tfrac72+4\sqrt{2},0),\quad (A_2, B_2)=(B_4, A_4)=(2+3\sqrt{2}, \tfrac12)\quad\text{and}\quad A_3=B_3=1+\sqrt{2}.
\end{align*}
Plugging these values into \eqref{equ:genforap} and using also $\Delta_K=8$ 
and $\zeta_K(2)=\frac{\pi^4}{48\sqrt2}$ (see Remark \ref{rmk:zeta2}), we get 
\begin{align}\notag
a_{\cP}&=\frac{3\sqrt{2}}{16\pi^4}\Area(\cW)\biggl(\int_{\tilde{S}_{2}}r(\theta)^{-2}\left((2+3\sqrt{2})+\tfrac12\nu(\theta)^{-2}\right)\,\text{d}\theta+\int_{ \tilde{S}_{4}}r(\theta)^{-2}\left(\tfrac12+(2+3\sqrt{2})\nu(\theta)^{-2}\right)\,\text{d}\theta
\\ \label{ABaPexplres1al}
&\hspace{20pt}+(\tfrac{7}{2}+4\sqrt{2})\left(\int_{\tilde{S}_{1}}r(\theta)^{-2}\,\text{d}\theta+\int_{\tilde{S}_{5}}r(\theta)^{-2}\nu(\theta)^{-2}\,\text{d}\theta\right)
+(1+\sqrt{2})\int_{\tilde{S}_{3}}r(\theta)^{-2}\left(1+\nu(\theta)^{-2}\right)\,\text{d}\theta\biggr).
\end{align}
Moreover, because of \eqref{equ:rela1} and $\nu(\theta+\pi)=\nu(\theta)^{-1}$,
we have the general symmetry relation that for any
bounded measurable set $S\subset \R_{>0}$,
\begin{align}\label{Gensplitting}
\int_{\{\theta\col\nu(\theta)\in S\}}r(\theta)^{-2}\nu(\theta)^{-2}\,\text{d}\theta
=\int_{\{\theta\col \nu(\theta)\in S^{-1}\}}r(\theta)^{-2}\,\text{d}\theta.
\end{align}
Using \eqref{Gensplitting}
and the fact that $S_{j}^{-1}=S_{6-j}$ for $1\leq j\leq 5$,
we see that \eqref{ABaPexplres1al} is equivalent to the formula \eqref{ABaPexplres1} in Proposition \ref{ABaPexpl}. 
Finally, if 
$\tilde{S}_3=[0,2\pi)$
then the sum inside the parenthesis in \eqref{ABaPexplres1} equals
$2(1+\sqrt2)\int_0^{2\pi}r(\theta)^{-2}\,\mathrm{d}\theta=4(1+\sqrt2)\Area(\scrW^*)$
by \eqref{INTell2inv2formula},
and hence \eqref{ABaPexplres2} holds.
This completes the proof of Proposition \ref{ABaPexpl}.

\vspace{5pt}

Next, we consider the case when $K=\Q(\sqrt{3})$. 
Then $\lambda=2+\sqrt3$,
and one verifies that
\begin{align*}
E_I^+=\{1,\lambda\}\cdot\lambda^{2\Z}=\lambda^{\Z}
\quad\text{and}\quad
E_I^-=\{\sqrt3,1+\sqrt3,\sqrt3\lambda,(1+\sqrt3)\lambda\}\cdot\lambda^{2\Z}=\{\sqrt3,1+\sqrt3\}\lambda^{\Z}
\end{align*}
with the associated values $\nu_\alpha$ 
given by
\begin{align*}
\nu_1=\nu_\lambda=1+\sqrt3,
\quad
\nu_{\sqrt3}=\nu_{\sqrt3\lambda}=\sqrt3\quad\mathrm{and}
\quad
\nu_{1+\sqrt3}=\nu_{(1+\sqrt3)\lambda}=-1+\sqrt3.
\end{align*}
Hence
\begin{align*}
E_{I,\nu}^+=\begin{cases}
\lambda^{\Z}&\text{if }\: \nu<1+\sqrt3,
\\
\emptyset&\text{if }\: \nu\geq1+\sqrt3,
\end{cases}
\end{align*}
and
\begin{align*}
E_{I,\nu}^-=\begin{cases}
\emptyset&\text{if }\: \nu\leq\sqrt3-1,
\\
\{1+\sqrt3\}\lambda^{\Z}&\text{if }\: \sqrt3-1<\nu\leq\sqrt3,
\\
\{\sqrt3,1+\sqrt3\}\lambda^{\Z}&\text{if }\: \sqrt3<\nu.
\end{cases}
\end{align*}
The corresponding partition of $\R_{>0}$ is $\R_{>0}=\sqcup_{j=1}^4S_{I,j}$
where $S_{I,j}=S_j$ 
with $S_1=(0,\sqrt{3}-1]$, $S_2=(\sqrt{3}-1,\sqrt{3}]$, $S_3=(\sqrt{3},\sqrt{3}+1)$ and $S_4=[\sqrt{3}+1,\infty)$.
The corresponding constants $A_j:=A_{I,j}$ and $B_j:=B_{I,j}$ 
are
\begin{align*}
(A_1,B_1)=(6+4\sqrt{3},0),\qquad (A_2,B_2)=(3+2\sqrt{3},2\sqrt{3})
\end{align*}
and
\begin{align*}
(A_3,B_3)=(2,3+8\sqrt{3}),\qquad (A_4,B_4)=(0,11+12\sqrt{3}).
\end{align*}
Plugging these formulas into \eqref{equ:genforap} and using also $\Delta_K=12$ and $\zeta_K(2)=\frac{\pi^4}{36\sqrt{3}}$ (see Remark \ref{rmk:zeta2}), we get
\begin{align}\notag
a_{\cP}=\frac{\sqrt{3}}{16\pi^4}\Area(\cW)\biggl(&(6+4\sqrt{3})\int_{\tilde{S}_{1}}r(\theta)^{-2}\,\text{d}\theta+\int_{\tilde{S}_{2}}r(\theta)^{-2}\left((3+2\sqrt{3})+2\sqrt{3}\nu(\theta)^{-2}\right)\,\text{d}\theta
\\ \label{GhaPexplres1al}
&+\int_{\tilde{S}_{3}}r(\theta)^{-2}\left(2+(3+8\sqrt{3})\nu(\theta)^{-2}\right)\,\text{d}\theta+(11+12\sqrt{3})\int_{\tilde{S}_{4}}r(\theta)^{-2}\nu(\theta)^{-2}\,\text{d}\theta\biggr).
\end{align}
Again by the symmetry identity \eqref{Gensplitting},
we see that \eqref{GhaPexplres1al} is equivalent to the formula 
\eqref{GHaPexplres1} in Proposition \ref{GHaPexpl}.
The last statement of Proposition \ref{GHaPexpl}
is again verified using \eqref{INTell2inv2formula}.
This completes the proof of Proposition \ref{GHaPexpl}.

Finally we consider the case when $K=\Q(\sqrt{5})$. 
Then $\lambda=\frac{1+\sqrt{5}}{2}$,
and one verifies that
\begin{align*}
E_I^+=\{1\}\cdot\lambda^{2\Z}=\lambda^{2\Z}
\quad\text{and}\quad
E_I^-=\{\lambda\}\cdot\lambda^{2\Z}=\lambda^{2\Z+1}\end{align*}
with the associated values $\nu_\alpha$
given by
\begin{align*}
\nu_1=\tfrac{1+\sqrt{5}}{2}\quad \mathrm{and}
\quad
\nu_{\lambda}=\tfrac{\sqrt5-1}{2}.
\end{align*}
Hence
\begin{align*}
E_{I,\nu}^+=\begin{cases}
\lambda^{2\Z}&\text{if }\: \nu<\frac{1+\sqrt5}{2},
\\
\emptyset&\text{if }\: \nu\geq\frac{1+\sqrt5}{2},
\end{cases}
\quad
\text{and}
\quad
E_{I,\nu}^-=\begin{cases}
\emptyset&\text{if }\: \nu\leq\frac{\sqrt5-1}{2},
\\
\lambda^{2\Z+1}&\text{if }\: \frac{\sqrt5-1}{2}<\nu.
\end{cases}
\end{align*}
The corresponding partition of $\R_{>0}$ is $\R_{>0}=\sqcup_{j=1}^3S_{I,j}$
where $S_{I,j}=S_j$ 
with $S_1=(0,\frac{\sqrt{5}-1}{2}]$, $S_2=(\frac{\sqrt{5}-1}{2},\frac{1+\sqrt{5}}{2})$ and $S_3=[\frac{1+\sqrt{5}}{2},\infty)$.
The corresponding constants $A_j:=A_{I,j}$ and $B_j:=B_{I,j}$ 
are
\begin{align*}
(A_1,B_1)=(B_3, A_3)=(\tfrac{5+3\sqrt{5}}{4},0)\quad \text{and}\quad (A_2,B_2)=(\tfrac{1+\sqrt{5}}{4},\tfrac{1+\sqrt{5}}{4}).
\end{align*}
Plugging these formulas into \eqref{equ:genforap} and using also $\Delta_K=5$ and $\zeta_K(2)=\frac{2\pi^4}{75\sqrt{5}}$ (see Remark \ref{rmk:zeta2}) and the symmetry identity \eqref{Gensplitting}, we obtain
\eqref{TTaPexplres1}.
The last statement of Proposition \ref{TTaPexpl}
is again verified using \eqref{INTell2inv2formula}.
This completes the proof of Proposition \ref{TTaPexpl}.

\section{Details on the examples and numerical computations}
\label{numcompSEC}
In Sections \ref{ABexplicitSEC}--\ref{TTTexplicitSEC} below
we show how the formulas
\eqref{PABwDEF}, \eqref{PGhwDEF} and \eqref{PTTTwDEF} 
follow from the cut-and-project presentation of these vertex sets 
in \cite[Ch.\ 7.3]{mBuG2013}.
Then in Section \ref{numcompSUBSEC}
we discuss numerical computations of 
the corresponding gap distributions.

\subsection{Ammann-Beenker tilings}
\label{ABexplicitSEC}

Recall the definition of $\scrW_{\AB}\subset\R^2$ in Section \ref{examplessec}.
Let $\xi_8=e^{\pi i/4}$ (an 8th root of unity),
and let $\star$ be the automorphism of the cyclotomic field $\Q(\xi_8)$
given by $\xi_8^\star=\xi_8^3$.
From now on we will use the standard identification of $\C$ with $\R^2$,
$z\leftrightarrow (\Re(z),\Im(z))$;
in particular we consider the ring
$\Z[\xi_8]$ as a subset of $\R^2$.
Set 
\begin{align}\label{ABgentransl}
\scrT_{\AB}:=\{\vecw\in \scrW_{\AB}\col 
\Z[\xi_8]\cap\partial(\scrW_{\AB}+\vecw)=\emptyset\}; 
\end{align}
this is the set of admissible ``generic'' translates of $\scrW_{\AB}$.
Next, for any $\vecw\in\R^2$, set
\begin{align}\label{ABPwdef}
\scrP_{\AB,\vecw}:=\{z\in\Z[\xi_8]\col z^\star\in \scrW_{\AB}+\vecw\} \qquad(\subset\R^2).
\end{align}
Then for every $\vecw\in \scrT_{\AB}$,
$\scrP_{\AB,\vecw}$ is the vertex set 
of an Ammann-Beenker tiling
(see \cite{Beenker82} and \cite[Ch.~7.3; Ex.\ 7.8]{mBuG2013});
and we note that $\bn\in\scrP_{\AB,\vecw}$. 
In fact
\textit{every} Ammann-Beenker tiling
(appropriately scaled and rotated, and translated so that $\bn$ is a vertex), 
has vertex set either equal to $\scrP_{\AB,\vecw}$ for some
$\vecw\in \scrT_{\AB}$, or equal to a \textit{limit}
of a sequence of such point sets, with respect to an appropriate metric on the 
family of locally finite point sets in $\R^2$.
(In the limit case,\label{ABlimitcasefootnote}
it follows that there exists some $\vecw\in \overline{\scrW_{\AB}}$
such that the vertex set $\scrP$ agrees with $\scrP_{\AB,\vecw}$
``up to density zero'', viz., the symmetric difference set
$\scrP\triangle\scrP_{\AB,\vecw}$ has asymptotic density zero in $\R^2$.
This implies that $\scrP$ and $\scrP_{\AB,\vecw}$ have the same limiting distribution 
of normalized gaps between directions, and therefore
the methods of the present paper apply also to these vertex sets, unless $\vecw\in\partial \scrW_{\AB}$.)

In order to rewrite \eqref{ABPwdef} as in \eqref{PABwDEF},
let $K=\Q(\sqrt2)$;
then $\scrO_K=\Z[\sqrt2]$.
Also set $g_1=\smatr101{\sqrt2}$
and $g_2:=\smatr10{1/\sqrt2\:}{1/\sqrt2}$. 
Using $\xi_8^2=\sqrt2\xi_8-1$, one verifies that
$\Z[\xi_8]=\scrO_K\oplus\scrO_K\xi_8$.
Note that for any $x_1,x_2\in\scrO_K$,
we have $x_1+x_2\xi_8=(x_1,x_2)g_2$ under our identification of $\C$ with $\R^2$,
and in particular, 
$\Z[\xi_8]=\scrO_K^2g_2$ as a subset of $\R^2$.
Furthermore, using $(\sqrt2)^\star=(\xi_8+\xi_8^{-1})^\star=\xi_8^3+\xi_8^{-3}=-\sqrt2$,
we have for any $(x_1,x_2)\in\scrO_K^2$:
\begin{align*}
(x_1+x_2\xi_8)^\star=\sigma(x_1)+\sigma(x_2)\xi_8^3
=(\sigma(x_1),\sigma(x_2))g_1^{-1},
\end{align*}
and this point lies in $\scrW_{\AB}+\vecw$ if and only if
$(\sigma(x_1),\sigma(x_2))\in(\scrW_{\AB}+\vecw)g_1$.
It follows that $\scrP_{\AB,\vecw}=\scrP((\scrW_{\AB}+\vecw)g_1,\scrL_K)g_2$,
i.e.\ the formula in \eqref{PABwDEF} holds.
(In the case $\vecw=\bn$, the formula \eqref{PABwDEF} was noted in 
\cite[{Sec.\ 3.2}]{Hammarhjelm2022}.)

\subsection{G\"ahler's shield tilings}
\label{GSexplicitSEC}
Recall the definition of $\scrW_{\Gh}\subset\R^2$ in Section~\ref{examplessec}.
Let $\xi_{12}=e^{\pi i/6}$ (a 12th root of unity),
and let $\star$ be the automorphism of the cyclotomic field $\Q(\xi_{12})$
given by $\xi_{12}^\star=\xi_{12}^5$.
As before, 
we identify $\C$ with $\R^2$ in the standard way.
For any $\vecw\in\R^2$, set 
\begin{align}\label{GHPwdef}
\scrP_{\Gh,\vecw}:=\{z\in\Z[\xi_{12}]\col z^\star\in \scrW_{\Gh}+\vecw\}
\qquad(\subset\R^2).
\end{align} 
Then for every 
$\vecw\in \scrW_{\Gh}$ satisfying $\Z[\xi_{12}]\cap\partial(\scrW_{\Gh}+\vecw)=\emptyset$,
$\scrP_{\Gh,\vecw}$ is the vertex set of a G\"ahler's shield tiling
\cite[Ch.\ 7.3; Ex.\ 7.12]{mBuG2013}.

In order to rewrite \eqref{GHPwdef} as in \eqref{PGhwDEF},
let $K=\Q(\sqrt3)$;
then $\scrO_K=\Z[\sqrt3]$.
Also set $g_1:=\smatr10{\sqrt3}2$
and $g_2:=\smatr10{\sqrt3/2}{1/2}$.
Using $\xi_{12}^2=\sqrt3\xi_{12}-1$ one verifies that
$\Z[\xi_{12}]=\scrO_K\oplus\scrO_K\xi_{12}$.
Note that for any $x_1,x_2\in\scrO_K$,
we have $x_1+x_2\xi_{12}=(x_1,x_2)g_2$ 
under our identification of $\C$ with $\R^2$,
and in particular, $\Z[\xi_{12}]=\scrO_K^2g_2$ as a subset of $\R^2$.
Furthermore, using
$(\sqrt3)^\star=(\xi_{12}+\xi_{12}^{-1})^\star
=\xi_{12}^5+\xi_{12}^{-5}=-\sqrt3$,
we have for any $(x_1,x_2)\in\scrO_K^2$:
\begin{align*}
(x_1+x_2\xi_{12})^{\star}=\sigma(x_1)+\sigma(x_2)\xi_{12}^5,
\end{align*}
and this point lies in $\scrW_{\Gh}+\vecw$
if and only if $(\sigma(x_1),\sigma(x_2))\in (\scrW_{\Gh}+\vecw)g_1$.
It follows that $\scrP_{\Gh,\vecw}=\scrP((\scrW_{\Gh}+\vecw)g_1,\scrL_K)g_2$,
i.e.\ the formula in \eqref{PGhwDEF} holds.

\subsection{T\"ubingen triangle tilings}
\label{TTTexplicitSEC}

Recall the definition of $\scrW_{\TT}\subset\R^2$ in Section~\ref{examplessec}.
Let $\xi_5=e^{2\pi i/5}$ (a 5th root of unity),
and let $\star$ be the automorphism of the cyclotomic field $\Q(\xi_5)$
given by $\xi_5^\star=\xi_5^2$.
As in Section \ref{ABexplicitSEC},
we identify $\C$ with $\R^2$ in the standard way.
For any $\vecw\in\R^2$, set 
\begin{align}\label{TTPwdef}
\scrP_{\TT,\vecw}:=\{z\in\Z[\xi_5]\col z^\star\in \scrW_{\TT}+\vecw\}
\qquad(\subset\R^2).
\end{align}
Then for every 
$\vecw\in \scrW_{\TT}$ satisfying $\frac15\Z[\xi_{5}]\cap\partial(\scrW_{\TT}+\vecw)=\emptyset$,
$\scrP_{\TT,\vecw}$ is the vertex set of a T\"ubingen triangle tiling
\cite[Ch.\ 7.3; Ex.\ 7.10]{mBuG2013},
\cite[Sec.\ 4]{BKSZ90}.

In order to rewrite \eqref{TTPwdef} as in \eqref{PTTTwDEF},
let $K=\Q(\sqrt5)$;
then $\scrO_K=\Z[\tau]$, where we recall that $\tau=\frac12(1+\sqrt5)$.
Also set $g_1:=\smatr{1}{0}0{\sqrt{(2+\tau)/5}}\smatr10{\tau}2$
and 
$g_2:=\smatr10{(\tau-1)/2\:}{\sqrt{2+\tau}/2}$.
Using 
$\xi_5^2=(\tau-1)\xi_5-1$ one verifies that
$\Z[\xi_5]=\scrO_K\oplus\scrO_K\xi_5$.
Note that for any $x_1,x_2\in\scrO_K$,
we have $x_1+x_2\xi_5=(x_1,x_2)g_2$ 
under our identification of $\C$ with $\R^2$,
and in particular, $\Z[\xi_5]=\scrO_K^2g_2$ as a subset of $\R^2$.
Furthermore, using
$(\sqrt5)^\star=(1+2\xi_5+2\xi_5^4)^\star=1+2\xi_5^2+2\xi_5^3
=-(1+2\xi_5+2\xi_5^4)=-\sqrt5$,
we have for any $(x_1,x_2)\in\scrO_K^2$:
\begin{align*}
(x_1+x_2\xi_5)^{\star}=\sigma(x_1)+\sigma(x_2)\xi_{5}^2,
\end{align*}
and this point lies in $\scrW_{\TT}+\vecw$
if and only if $(\sigma(x_1),\sigma(x_2))\in (\scrW_{\TT}+\vecw)g_1$. 
It follows that $\scrP_{\TT,\vecw}=\scrP((\scrW_{\TT}+\vecw)g_1,\scrL_K)g_2$,
i.e.\ the formula in \eqref{PTTTwDEF} holds.

\subsection{Numerical computations}
\label{numcompSUBSEC}

Next we discuss numerical computations
of the  
gaps between directions
of 
points in the Ammann-Beenker, G\"ahler's shield, and T\"ubingen triangle tilings,
$\scrP_{\AB,\vecw}$, $\scrP_{\Gh,\vecw}$ and $\scrP_{\TT,\vecw}$.
In Table \ref{asympttable1} below,
for a few examples of points $\vecw\in\scrW_{\AB}$,
we present conjectural approximate values of $s^2F(s)$,
where $F(s)$ is the limiting gap distribution function
associated to $\scrP_{\AB,\vecw}$,
and, for comparison, 
the exact values of $a_{\scrP}$.
Tables \ref{asympttable2}  and \ref{asympttable3} 
present similar types of data for
$\scrP_{\Gh,\vecw}$ and $\scrP_{\TT,\vecw}$.

\begin{table}[h!]
\textbf{Experimental 
approximate values of $s^2F(s)$ for $\scrP=\scrP_{\AB,\vecw}$}
\begin{center}
\begin{tabular}{c | c || rrrrrrr}
 \hline
$\vecw$ & $a_{\scrP}$  
& \hspace{3pt} $s=10$ & \hspace{3pt} $s=20$ & \hspace{3pt} $s=50$ 
& \hspace{3pt} $s=100$ & \hspace{3pt} $s=200$  \hspace{-5pt} &  $s=500$ \hspace{-5pt} & $s=1000$ \hspace{-10pt}
\\ \hline
$(0,0)$\rule{0pt}{16pt} & 0.24638\ldots & 0.2810 & 0.2628 & 0.2525 & 0.2497 & 0.250 & 0.246 & 0.23
\\[2pt]
$(0.8,0.1)$ & 0.21809\ldots & 0.2460 & 0.2317 & 0.2236 & 0.2212 & 0.220 & 0.219 & 0.21
\\[2pt]
$(0.9,0.3)$ & 0.20416\ldots & 0.2286 & 0.2160 & 0.2090 & 0.2067 & 0.206 & 0.205 & 0.19
\\[2pt]
$(1.2,0.5)$ & 0.17444\ldots & 0.2211 & 0.2048 & 0.1815 & 0.1767 & 0.177 & 0.178 & 0.18
\\[3pt]\hline
\multicolumn{2}{c||}{Error\rule[-6pt]{0pt}{20pt} $\lesssim$} & 0.0005 & 0.0006 & 0.0011 & 0.0023 & 0.005 & 0.015 & 0.05
\\\hline
\end{tabular}
\end{center}
\caption{Values of $a_{\scrP}$ and approximate values of $s^2F(s)$,
for $\scrP=\scrP_{\AB,\vecw}$. 
The last line gives a conjectural approximate upper bound on the absolute error
of each value in the corresponding column.
(Thus, for example, we expect that for $\vecw=\bn$ and $s=50$,
$|s^2F(s)-0.2525|\lesssim0.0011$.)}
\label{asympttable1}
\end{table}

\begin{table}[h!]
\textbf{Experimental 
approximate values of $s^2F(s)$ for $\scrP=\scrP_{\Gh,\vecw}$}
\begin{center}
\begin{tabular}{c | c || rrrrrrr}
 \hline
$\vecw$ & $a_{\scrP}$  
& \hspace{3pt} $s=10$ & \hspace{3pt} $s=20$ & \hspace{3pt} $s=50$ 
& \hspace{3pt} $s=100$ & \hspace{3pt} $s=200$  \hspace{-5pt} &  $s=500$ \hspace{-5pt} & $s=1000$ \hspace{-10pt}
\\ \hline
$(1.7,0.6)$\rule{0pt}{16pt} & 0.17790\ldots & 0.1974 & 0.1873 & 0.1817 & 0.180 & 0.180 & 0.176 & 0.19
\\[2pt]
$(0.1,0)$ & 0.21374\ldots & 0.2460 & 0.2289 & 0.2199 & 0.216 & 0.217 & 0.215 & 0.23
\\[2pt]
$(1.5,-0.3)$ & 0.19210\ldots & 0.2143 & 0.2027 & 0.1964 & 0.194 & 0.195 & 0.190 & 0.20
\\[2pt]
$(0.3,1.2)$ & 0.20870\ldots & 0.2341 & 0.2208 & 0.2136 & 0.211 & 0.211 & 0.205 & 0.23
\\[3pt]\hline
\multicolumn{2}{c||}{Error\rule[-6pt]{0pt}{20pt} $\lesssim$} & 0.0004 & 0.0007 & 0.001 & 0.004 & 0.006 & 0.019 & 0.05
\\\hline
\end{tabular}
\end{center}
\caption{Values of $a_{\scrP}$ and approximate values of $s^2F(s)$,
for $\scrP=\scrP_{\Gh,\vecw}$.}
\label{asympttable2}
\end{table}

\begin{table}[h!]
\textbf{Experimental 
approximate values of $s^2F(s)$ for $\scrP=\scrP_{\TT,\vecw}$}
\begin{center}
\begin{tabular}{c | c || rrrrrrr}
 \hline
$\vecw$ & $a_{\scrP}$  
& \hspace{3pt} $s=10$ & \hspace{3pt} $s=20$ & \hspace{3pt} $s=50$ 
& \hspace{3pt} $s=100$ & \hspace{3pt} $s=200$  \hspace{-5pt} &  $s=500$ \hspace{-5pt} & $s=1000$ \hspace{-10pt}
\\ \hline
$(0.5,0.4)$\rule{0pt}{16pt} & 
0.26135\ldots & 0.2919 & 0.2762 & 0.2673 & 0.2643 & 0.263 & 0.265 & 0.26
\\[2pt]
$(0.4,0)$ & 0.28875\ldots & 0.3283& 0.3080 & 0.2965& 0.2927& 0.291& 0.294& 0.30
\\[2pt]
$(0.15,0.3)$ &  0.29103\ldots & 0.3329 & 0.3116 & 0.2995 & 0.2954 & 0.294 & 0.297 & 0.30
\\[2pt]
$(1.3,0.4)$ & 0.18732\ldots & 0.2195 & 0.1977 & 0.1907 & 0.1892 & 0.188 & 0.190 & 0.19
\\[3pt]\hline
\multicolumn{2}{c||}{Error\rule[-6pt]{0pt}{20pt} $\lesssim$} & 
0.0004 & 0.0004 & 0.0007 & 0.0017 & 0.004 & 0.014 & 0.03
\\\hline
\end{tabular}
\end{center}
\caption{Values of $a_{\scrP}$ and approximate values of $s^2F(s)$,
for $\scrP=\scrP_{\TT,\vecw}$.}
\label{asympttable3}
\end{table}

The approximations of $F(s)$ 
in these tables were obtained by collecting all the points of 
$\scrP$ (where $\scrP=\scrP_{\AB,\vecw}$ or $\scrP_{\Gh,\vecw}$ or $\scrP_{\TT,\vecw}$) in fairly large discs $\scrB_R^2$
and evaluating the ratio in the left hand side of 
Theorem \ref{gaplimitdistrexistsTHM};
by repeating this  
for several large $R$-values we also obtained
conjectural error bounds;
see  
below for a more detailed discussion.
Recall that $\lim_{s\to\infty}s^2F(s)=a_{\scrP}$
by Theorem~\ref{MAINTHM1},
in fact with a rate $O(s^{-\frac12+\ve})$.
We certainly expect that the absolute difference
$|s^2F(s)-a_{\scrP}|$ typically \textit{decreases} as $s$ runs through the
values $10,20,50,100,200,500,1000$,
and 
it should be noted that the data in Tables~\ref{asympttable1}--\ref{asympttable3} 
is consistent with this hypothesis,
when considering also the conjectural error bounds 
displayed in the last line of the tables.

Precise numerical values of $a_{\scrP}$
such as those presented in Tables~\ref{asympttable1}--\ref{asympttable3}
(but also to much higher precision),
are quite easy to compute
using Propositions~\ref{ABaPexpl}--\ref{TTaPexpl}.
To describe how to do this, first note that 
in the case $\scrP=\scrP_{\AB,\vecw}$,
it follows from \eqref{PABwDEF},
the last statement in Theorem~\ref{gaplimitdistrexistsTHM}
and Remark~\ref{rmk1:sl2inva},
that this $a_{\scrP}$ is the same as for
$\scrP=\scrP(\scrW_{\AB}+\vecw,\scrL_K)$;
similarly, for $\scrP=\scrP_{\Gh,\vecw}$,
$a_{\scrP}$ is the same as for
$\scrP=\scrP(\scrW_{\Gh}+\vecw,\scrL_K)$ (see  \eqref{PGhwDEF}),
and for $\scrP=\scrP_{\TT,\vecw}$,
$a_{\scrP}$ is the same as for
$\scrP=\scrP(\scrW_{\TT}+\vecw,\scrL_K)$ (see  \eqref{PTTTwDEF}).
Now, somewhat more generally,
assume that $\scrW=\scrW_0+\vecw$ where
$\scrW_0$ is the open regular $N$-gon with vertices at the points
$R_0 e^{2\pi i(\frac1{2N}+\frac kN)}$ ($k\in\Z/N\Z$),
for some integer $N\geq3$ and some $R_0>0$,
and $\vecw\in\scrW_0$.
(We obtain $\scrW_0=\scrW_{\AB}$ by taking
$N=8$ and $R_0=\sqrt{1+\sqrt{1/2}}$,$\:$
$\scrW_0=\scrW_{\Gh}$ by taking 
$N=12$ and $R_0=\sqrt{2+\sqrt3}$,
and $\scrW_0=\scrW_{\TT}$ by taking $N=10$ and $R_0=\frac{\sqrt2}{20}(5+\sqrt5)^{3/2}$.)
For such a window $\scrW$, 
the function $r(\theta)$ in \eqref{equ:ellw}
is given by
\begin{align*}
r(\theta)=\Re\bigl(r_j e^{i(\theta_j+\theta)}\bigr)
=r_j\cos(\theta_j+\theta),
\qquad\forall \theta\in\bigl[(j-1)\tfrac{2\pi}N,j\tfrac{2\pi}N\bigr]+2\pi\Z,
\end{align*}
where 
for each $j\in\Z/N\Z$,
$r_j>0$ and $\theta_j\in\R/2\pi\Z$ 
are defined by $\vecw+R_0e^{2\pi i (\frac1{2N}-\frac jN)}=r_je^{i\theta_j}$.
Using also $\nu(\theta)=r(\theta+\pi)/r(\theta)$,
it is then easy to verify that 
all the sets $\tilde{S}_j$ in 
Propositions~\ref{ABaPexpl}--\ref{TTaPexpl}
are finite unions of intervals.
For given $\scrW_0$ and $\vecw$,
these intervals can be computed numerically,
and with this, also 
the integrals in \eqref{ABaPexplres1} and \eqref{GHaPexplres1} can be computed.
This is carried out in the program \cite[explicit\_aP.mpl]{HSYsupplement}.

\vspace{5pt}

We next describe more in detail how the experimental approximate values of
$F(s)$ in
Tables~\ref{asympttable1}--\ref{asympttable3} were computed.
(This is similar as in \cite{mBfGcHtJ2014} and \cite{Hammarhjelm2022}, 
but we have carried out more extensive computations.)
With notation as in 
Theorem \ref{gaplimitdistrexistsTHM}, set 
\begin{align}\label{FRsdef}
F_R(s):=\frac{\#\{1\leq j\leq N(R)\col N(R)(\xi_{R,j}-\xi_{R,j-1})\geq s\}}{N(R)}.
\end{align}
Then $F(s)=\lim_{R\to\infty}F_R(s)$,
and so for $R$ sufficiently large, $F_R(s)$ is a good approximation of $F(s)$.
Note that computing $F_R(s)$ involves collecting all the points of
$\scrP$ in the disc $\scrB_R^2$.
When carrying out the computations,
we noted that for fixed $s$,
the values $F_R(s)$ typically \textit{oscillate} as a function of $R$;
hence if $\fR$ is a finite set of several sufficiently large $R$-values
of the same order of magnitude,
it seems reasonable to consider the difference
\begin{align*}
\Delta_{\fR}F(s):=\max_{R\in\fR}F_R(s)-\min_{R\in\fR}F_R(s)
\end{align*}
as an \textit{indication} of the size of the error
$|F_R(s)-F(s)|$ for any $R\in\fR$.
To lessen the probability of getting a small difference 
by coincidence,
we took 
the experimental error bound given in the last line of
Tables \ref{asympttable1}--\ref{asympttable3} to be the
\textit{maximum} of $\ts^2\Delta_{\fR}F(\ts)$ 
over $100$ values of $\ts$
roughly equispaced in the interval $[\frac 45s,\frac65s]$,
and over all four choices $\vecw$.
Because of the factor $\ts^2$, 
it is not surprising  
that these experimental error bounds 
get worse as $s$ increases.

In order to produce Table \ref{asympttable1}
(i.e.\ the case of Ammann-Beenker tilings)
we used 
\begin{align}\label{fRDEF}
\fR=\{22000,23000,24000,25000\},
\end{align}
and for each $s$ we 
chose as approximate value of $F(s)$
the mid-point
$\frac12\bigl(\min_{R\in\fR}F_R(s)+\max_{R\in\fR}F_R(s)\bigr)$.
We remark that computing the functions $F_R$ in this case involve collecting a bit over
$2.3\cdot10^9$ points from each point set $\scrP_{\AB,\vecw}$.
For example, for $R=25000$ we found the number of points in
$\scrP_{\AB,\bn}\cap\scrB_R^2\setminus\{\bn\}$ 
to be $2370148592$
(this may be compared with $c_{\scrP}\pi R^2=\frac{1+\sqrt2}2 \pi R^2=2370148622.42\ldots$).\label{explnumbpointsR25000}

Similarly, for Table \ref{asympttable2} (G\"ahler's shield tilings)
we used
\begin{align}\label{fRDEFGh}
\fR=\{12000,13000,14000,15000\},
\end{align}
and computing the functions $F_R$ involved collecting a bit over
$2.6\cdot10^9$ points from each point set $\scrP_{\Gh,\vecw}$.
For example, for $R=15000$ and $\vecw=(1.7,0.6)$
we found the number of points in
$\scrP_{\Gh,\vecw}\cap\scrB_R^2\setminus\{\bn\}$ 
to be $2638031485$
(which may be compared with $c_{\scrP}\pi R^2=(2+\sqrt3) \pi R^2=2638031264.97 \ldots$).\label{explnumbpointsGhR25000}
\label{gHcPformula}

Similarly, for Table \ref{asympttable3} (T\"ubingen triangle tilings) we used
\begin{align}\label{fRDEFTTT}
\fR=\{17000,18000,19000,20000\},
\end{align}
and computing the functions $F_R$ involved collecting a bit over
$2.5\cdot10^9$ points from each point set $\scrP_{\Gh,\vecw}$.
For example, for $R=20000$ and $\vecw=(0.4,0)$
we found the number of points in
$\scrP_{\TT,\vecw}\cap\scrB_R^2\setminus\{\bn\}$ 
to be $2503118410$
(which may be compared with $c_{\scrP}\pi R^2=\frac25(\tau+1)\sqrt{\tau+2}\cdot \pi R^2=
2503118771.19 \ldots$).
\label{TTTcPformula}

\vspace{5pt}

Let us next discuss the graphs in Figure \ref{densityplot1}.
These were generated by
distributing the normalized gaps
$N(R)(\xi_{R,j}-\xi_{R,j-1})$, $j=1,\ldots,N(R)$
-- discarding the vanishing gaps --
in bins of width $0.02$,
and then drawing the corresponding histogram, appropriately scaled.
In other words, using the notation in \eqref{FRsdef},
each panel displays the graph of the 
function 
\begin{align}\label{ERhdef}
s\mapsto E_{R,h}(s):=\frac{\tF_R(s_h)-F_R(s_h+h)}h,
\qquad\text{with }\: s_h:=h\lfloor s/h\rfloor
\text{ and }h:=0.02,
\end{align}
where $\tF_R(0):=F_R(0+):=\lim_{s\to 0+}F_R(s)$ while $\tF_R(s):=F_R(s)$ for all $s>0$.\footnote{The fact that we
use $\tF_R$ in place of $F_R$ corresponds exactly to the fact that we
discard those gaps
$N(R)(\xi_{R,j}-\xi_{R,j-1})$ which vanish;
without this modification the graph would have a huge peak for $s\in[0,h)$;
see the discussion below.}
We point out that 
the graphs of $E_{R,h}$,
as they are scaled in 
Figure \ref{densityplot1},
look identical to the naked eye 
for each $R\in\fR$,
i.e.\ the picture looks the same if we plot all four graphs in the same coordinate system;
in fact 
the maximum over $s\in[0,6]$ of the difference
$\max_{R\in\fR}E_{R,h}(s)-\min_{R\in\fR}E_{R,h}(s)$
was in each of the six cases
found to be less than $0.0018$, and the average of that same difference was
less than $0.0004$.

These findings support the conjecture that, 
in the cases considered, the derivative $F'(s)$ exists, is continuous, 
and is essentially correctly depicted by the graphs in 
Figure \ref{densityplot1}.

It is also of interest to consider the 
\textit{areas under the curves} in 
Figure \ref{densityplot1}.
Note that it follows from
\eqref{ERhdef} that the area under the graph of $E_{R,h}(s)$ is
\begin{align*}
\int_0^\infty E_{R,h}(s)\,ds=F_R(0+).
\end{align*}
Furthermore,
as $R\to\infty$, $F_R(0+)$ tends to $\kappa_{\scrP}$,
the relative density of visible points in \label{kappaPDEFrep}
$\scrP$, 
a quantity which was defined and proved to exist for general regular cut-and-project sets
in \cite{MarklofStrombergsson2015}.
(The proof of $\lim_{R\to\infty}F_R(0+)=\kappa_{\scrP}$ is immediate
by comparing with the definition of $\kappa_{\scrP}$ in
\cite{MarklofStrombergsson2015}.)
The values which we obtained for $F_R(0+)$ for varying $R$ in $R\in\fR$ 
consistently agreed to within ca $3\cdot 10^{-6}$;
hence we expect that these values  
approximate $\kappa_{\scrP}$ to within an absolute error of similar order of magnitude.
In particular for $\scrP_{\AB,\vecw}$ with $\vecw=\bn$ and $R=25000$, we obtained 
${\displaystyle F_R(0+)=\frac{1368315872}{2370148592}=0.57731227\ldots}$.\footnote{Similar numerics,
but for considerably smaller $R$-values,
also appears in
\cite[{Table 2}]{mBfGcHtJ2014}
and
\cite[{Table 1}]{Hammarhjelm2022}.}
Thus, this value  
is the area under the curve in panel I of Figure \ref{densityplot1},
and it may be compared with the 
known exact value
\cite[{p.\ 10}]{Hammarhjelm2022}\label{density0p57731262}
\begin{align}\label{kappaPforsymmAB}
\kappa_{\scrP}=\frac{2|\sigma(\lambda)|}{\zeta_K(2)}
=\frac{2(\sqrt2-1)}{\pi^4/(48\sqrt2)}
=0.57731262\ldots. 
\end{align}
For the other panels in Figure \ref{densityplot1}
we obtained the following values of $F_R(0+)$ (viz., area under the curve):
II. $0.71023\ldots$.
III. $0.76447\ldots$.
IV. $0.62863\ldots$.
V. $0.60173\ldots$.
VI. $0.73515\ldots$.

Finally, to avoid possible confusion,
let us remind about the fact that,
as explained in Remark~\ref{visibleREM},
in the present paper we are using a different normalization of ``$F(s)$''
than what was used in \cite{mBfGcHtJ2014}, \cite{MarklofStrombergsson2015} and \cite{Hammarhjelm2022}.
This has the effect that the graph in the top left 
panel of Figure \ref{densityplot1}
is quantitatively strongly different from the graphs in
\cite[Fig.\ 9]{mBfGcHtJ2014}
and
\cite[Fig.\ 2]{Hammarhjelm2022},
despite the fact that each of these three graphs
depict the gap distribution
for directions in the point set $\scrP=\scrP_{\AB,\bn}$.
In fact, in our notation,
the plots in
\cite[Fig.\ 9]{mBfGcHtJ2014}
and
\cite[Fig.\ 2]{Hammarhjelm2022}
are experimental graphs of 
the function $s\mapsto -\kappa_{\scrP}^{-2}F'(\kappa_{\scrP}^{-1}s)$,
with $\kappa_{\scrP}$ as in \eqref{kappaPforsymmAB}.
Taking this rescaling into account, the three graphs agree fairly well\label{compGustavsFig2}
--- one naturally expects 
that the more erratic
small-scale behavior of the graphs in \cite{mBfGcHtJ2014}
and
\cite{Hammarhjelm2022}
is due to the fact that those were calculated using 
considerably smaller $R$-values than those used for 
our Figure \ref{densityplot1}.

\section*{Index of notation}

\begin{center}
\begin{footnotesize}
\begin{longtable}{llr}
$a_{\cP}$ & leading coefficient of the limiting gap distribution function $F(s)$ of a planar point set $\cP$ & \pageref{MAINTHM1}
\\
$\text{a}_{\bm{y}}$ & matrix element $(\text{a}_{y_1}, \text{a}_{y_2})$ with $\text{a}_y=\left(\begin{smallmatrix}
y & 0\\
0 & y^{-1}\end{smallmatrix}\right)$ & \pageref{p:iwdec}
\\
$\cA$ & closure of $\pi_{\intl}(\cL)$ with $\cL$ a lattice in $\R^{d+m}$ and $\pi_{\intl}$ the natural projection to the internal space & \pageref{scrADEF}
\\
$\cA^{\circ}$ & identity component of $\cA$ & \pageref{scrAcircDEF}
\\
$\fA_{\vectheta,\vecy}(s)$ & the set $(y_1^{-1}s^{1/2}\lambda^{-2r}\ell(\theta_1))\times (y_2^{-1}\lambda^{2r}\ell_{\scrW}(\theta_2))$ & \pageref{GEisbounddisc2}
\\
$\alpha_{I,J}(y)$ & the function $\min\left\{\alpha\in I\cap\R_{>0}\col y\,\sigma(\alpha)\in J\right\}$ with $I\subset \cO_K$ an integral ideal & \\
 & and $J\subset \R$ bounded measurable with non-empty interior  & \pageref{alphapWiDEF}
\\
$\tilde{\alpha}_{I,J}(y)$ & the function $\max\bigl\{\alpha\in I\cap\R_{<0}\col y\,\sigma(\alpha)\in J\bigr\}$ with $I, J$ as above & \pageref{moreexplicitPROPpf2}
\\
$b_{\alpha,\nu}$ & $b_{\alpha,\nu}:=\sigma(\alpha)$ if $\sigma(\alpha)>0$ and $b_{\alpha,\nu}:=\nu^{-1}\sigma(\alpha)$ if $\sigma(\alpha)<0$ 
 & \pageref{p:balpha1}
\\
$\scrB_R^2$ & open disc in $\R^2$ with center $\bm{0}$ and radius $R$  & \pageref{p:ball}
\\
$\fB_{\vectheta,\vecy,\vecx,\alpha}(s)$ & the set $\left(y_1\mathfrak{J}_1-\alpha x_1\right)\times \left(y_2\mathfrak{J}_2-\sigma(\alpha)x_2\right)$ & \pageref{condfinEQU}
\\
$\mathfrak{c}_i$ & the vector $(-c_i, a_i)\in \cO_K^2$  such that  $\frac{a_i}{c_i}=k_i$ is the $i$-cusp of $\G_K$ & \pageref{p:frakci}
\\
$C_K$ & ideal class group of $K$ & \pageref{p:clgrp}
\\
$c_K$ & normalizing factor of a Haar measure of $\mathrm{H}$ which equals $\Delta_K^{-3/2}\zeta_K(2)^{-1}$ & \pageref{equ:compck}
\\
$c_{\cP}$ & asymptotic density of a planar point set $\cP$ & \pageref{asymptdensityDEF}
\\
$\Delta_K$ & discriminant of $K$   & \pageref{cPformula}, \pageref{discriminantKdef}
\\
$\Delta_R$ & directions of points of length bounded by $R$ in a planar point set $\cP$ & \pageref{deltaRDEF}
\\
$E_I^{\pm}$ & set of positively or negatively extremal points with respect to an integral ideal $I$ & \pageref{p:eipm}
\\
$E_{I,\nu}$ & set of $\nu$-extremal points with respect to an integral ideal $I$ & \pageref{p:einu}
\\
$E^{\pm}_{I,\nu}$ & set of positive or negative $\nu$-extremal points with respect to $I$ & \pageref{p:einu}
\\
$F(s)$ & limiting distribution function of gaps of directions of a planar cut-and-project set & \pageref{gaplimitdistrexistsTHMres}
\\
$\mathfrak{f}(A,a)$ & the function $a^{-1}\inf\bigl\{m\bigl((x,x+a)\setminus A\bigr)\col x\in\R\bigr\}$ with $a>0$ and $A\subset \R$ Lebesgue measurable & \pageref{ffAaDEF}
\\
$\tff(A,a)$ & either $0$ or $\mathfrak{f}(A,a)$ depending on whether $A$ is empty or not & \pageref{tffAaDEF}
\\
$\cF_{\G_K}$ & a Siegel fundamental domain for $\G_K\bk \mathrm{H}$ & \pageref{equ:siegefd}
\\
$\cF_i$ & a fixed fundamental domain for $\G_i\bk H$ & \pageref{equ:cusp1}
\\
$\cF_i(t)$ & the set $\left\{\xi_i\text{n}_{\bm{x}}\text{a}_{\bm{y}}\text{k}_{\bm{\theta}}\in \cF_i\col y_1y_2\geq t\right\}$ & \pageref{equ:cusp2}
\\
$\mathfrak{F}_i$ & a fixed fundamental domain for $\R^2/\iota(I_i^{-2})$ & \pageref{p:fFi}
\\
$G(s)$ & integral of the limiting gap distribution function $F(s)$   
& \pageref{Gdef}
\\
$G_i(s)$ & contribution of the $i$-th cusp of $\G_K$ to the function $G(s)$ for sufficiently large $s$ & \pageref{GisDEF}
\\
$G_{E,i}(s)$ & error term of $G_i(s)$ & \pageref{GEisDEF}
\\
$G_{M,i}(s)$ & main term of $G_i(s)$ & \pageref{GMisDEF}
\\
$\G$ & the group $\SL_4(\Z)$ & \pageref{p:sl4z}
\\
$\G_i$ &  isotropy group of the $i$-th cusp of $\G_K$ & \pageref{equ:isotropygp}
\\
$\Gamma_K$ & the group $\iota(\SL_2(\cO_K))$ & \pageref{p:gammak}
\\
$\mathrm{H}$ & the group $\SL_2(\R)\times \SL_2(\R)$ & \pageref{p:gph}
\\
$\mathrm{H}_g$ & the group $g\mathrm{H}g^{-1}$ with $g\in \SL_4(\R)$ & \pageref{p:hg}
\\
$I_i$ & ideal generated by $(a_i, c_i)\in \cO_K^2$ with $a_i, c_i$ such that $k_i=\frac{a_i}{c_i}$ & \pageref{p:idecuspi}
\\
$\iota$ &  the natural embedding either from $K$ to $\R^2$ or from $\SL_2(K)$ to $\mathrm{H}$  & \pageref{equ:embd}, \pageref{equ:iotaemb}
\\
$J$ & the set $\bigl(0,\tfrac{\pi}{4}\bigr) \cup \bigl(\tfrac{3\pi}{4},\pi\bigr)$  & \pageref{p:intervalj}
\\
$J_{\nu}$ & open interval $(-\nu, 1)$ with $\nu>0$ 
 & \pageref{p:jnu}
\\
$\mathfrak{J}_1$  & the set $\left\{t\in \R\col (y_1\alpha,t)\in s^{1/2}\lambda^{-2r}T(1)\text{k}_{-\theta_1}\right\}$ & \pageref{J1J2def}
\\
$\mathfrak{J}_2$ & the set $\left\{t\in\R\col (y_2\sigma(\alpha),t)\in \lambda^{2r}\cW \text{k}_{-\theta_2}\right\}$ & \pageref{J1J2def}
\\
$K$ & real quadratic field $\Q(\sqrt{d})$ with $d\geq 2$ square-free & \pageref{p:qufield}, \pageref{p:qufield_rep}
\\
$\overline{K}$ & the set $K\cup\{\infty\}$ & \pageref{p:overlinek}
\\
$k_i$ &  $i$-th cusp of $\G_K$ & \pageref{p:cusps}
\\
$\text{k}_{\bm{\theta}}$ &matrix element $(\text{k}_{\theta_1}, \text{k}_{\theta_2})$ with $\text{k}_{\theta}=\left(\begin{smallmatrix}
\cos\theta & -\sin\theta\\
\sin\theta & \cos\theta\end{smallmatrix}\right)$ & \pageref{p:iwdec}
\\
$\kappa$ & number of cusps of $\G_K$ & \pageref{p:kappa}
\\
$\kappa_{\cP}$ & relative density of visible points in a planar point set $\cP$ & \pageref{kappaPDEF}, \pageref{kappaPDEFrep}
\\
$\mu_K$ & the probability $\mathrm{H}$-invariant  measure on $\G_K\bk \mathrm{H}$ & \pageref{equ:haar11}
\\
$\cL_i$ & the set $\Pi_i\cap \cL_K$ 
 & \pageref{equ:fillplane2}
\\
$\widetilde{\cL}_i$ & the set $\left\{(\alpha, t_1,\sigma(\alpha), t_2)\col \alpha\in I_i,\ t_1, t_2\in \R\right\}\xi^{-1}_i$ & \pageref{tscrLiDEF}
\\
$\cL_K$ & Minkowski embedding of $\cO_K^2$ in $\R^4$ & \pageref{scrLKdef}
\\
$\ell(\theta)$ & projection of $T(1)\text{k}_{-\theta}$ on the $x$-axis with $\theta\in \R/2\pi\Z$ & \pageref{p:elltheta}
\\
$\ell_{\cW}(\theta)$ & projection of $\cW\text{k}_{-\theta}$ on the $x$-axis  & \pageref{p:elltheta}
\\
$\lambda$ & fundamental unit of $K$  & \pageref{p:lambda}
\\
$m$ & Lebesgue measure on $\R$ & \pageref{mLEBdef}
\\
$\mathrm{N}$ & standard norm on $K$ defined by $\mathrm{N}(\alpha)=\alpha\sigma(\alpha)$ & \pageref{p:stnorm}
\\
$\text{Nr}$ & absolute norm of ideals & \pageref{equ:covolfor}
\\
$\text{n}_{\bm{x}}$ & matrix element $(\text{n}_{x_1}, \text{n}_{x_2})$ with $\text{n}_x=\left(\begin{smallmatrix}
1 & x\\
0 & 1\end{smallmatrix}\right)$ & \pageref{p:iwdec}
\\
$\cO_K$ & ring of integers of $K$ & \pageref{p:qufield}
\\
$\cO_K^{\times}$ & unit group of $\cO_K$ & \pageref{p:unitgp}
\\
$\nu_{\alpha}$ & critical value where the $\nu$-extremality of $\alpha\in E_I^{+}\sqcup E_I^-$ changes & \pageref{lem:semicont}
\\
$\nu(\theta)$ & parameter in the expression $\ell_{\cW}(\theta)=r(\theta)(-\nu(\theta), 1)$ & \pageref{equ:ellw}
\\
$\scrP_{\AB,\vecw}$ & the Ammann-Beenker tiling corresponding to the translate $\vecw$ & \pageref{PABwDEF}, \pageref{ABPwdef}
\\
$\scrP_{\Gh, \vecw}$ & the G\"ahler's shield tiling corresponding to the translate $\vecw$ & \pageref{PGhwDEF}, \pageref{GHPwdef}
\\
$\cP(\cW, \cL)$ & cut-and-project set associated to a window set $\cW$ and a lattice $\cL$ & \pageref{scrPWLdef}
\\
$\pi$ & natural projection from $\R^{d+m}$ to the physical space $\R^d$ & \pageref{p:projections}
\\
$\pi_{\intl}$ & natural projection from $\R^{d+m}$ to the internal space $\R^m$ & \pageref{p:projections}
\\
$\Pi_i$ &  the plane $(\{0\}\hspace{-3pt}\times\hspace{-3pt} \R\hspace{-3pt}\times\hspace{-3pt}\{0\}\hspace{-3pt}\times\hspace{-3pt}\R)\xi^{-1}_i$ & \pageref{equ:fillplane1}
\\
$R_{\alpha, \nu}$ & a rectangle depending on $\alpha\in I\cap \R_{>0}$ and $\nu>0$ with $I$ an integral ideal  & \pageref{p:ralphanu}
\\
$R_{\cW}(\theta,x)$ & the set $\{y\in\R\col (x,y)\text{k}_\theta\in\scrW\}$ & \pageref{RWdef}
\\
$r(\theta)$ & parameter in the expression $\ell_{\cW}(\theta)=r(\theta)(-\nu(\theta), 1)$ & \pageref{equ:ellw}
\\
$S_1^1$ & the unit circle & \pageref{p:unitcirc}
\\
$S_{I_i,j}$ & the $j$-th interval in the partition of $\R_{>0}$ corresponding to the $i$-th cusp of $\G_K$ & \pageref{thm:moreexfor}
\\
$\tilde{S}_{I_i,j}$ & the set $\left\{\theta\in [0, 2\pi)\col \nu(\theta)\in S_{I_i,j}\right\}$ & \pageref{thm:moreexfor}
\\
$\SL_2(\cO_K)$ & Hilbert modular group & \pageref{equ:slok}
\\
$\sigma$ & the unique nontrivial automorphism of the real quadratic field $K$ & \pageref{p:sigma} 
\\
$t_1$ & an absolute positive number depending only on $K$ & \pageref{p:t1}
\\
$\tau$ &  a generator of $\cO_K$ which is
$\frac{1+\sqrt{d}}{2}$ if $d\equiv 1\Mod{4}$ and $\sqrt{d}$ if $d\equiv 2, 3\Mod{4}$ 
 & \pageref{p:sigma}
\\
$T(s)$ & open triangle with vertices at $(0,0)$ and $(s/c_\scrP)^{1/2}(1,\pm 1)$ & \pageref{p:tstriangle}
\\
$\scrT(s)$ & the set $\lambda^{-2r}T(s)\times \lambda^{2r}\cW$ with $r\in \N$ such that $\lambda^{2r}\leq s^{1/4}<\lambda^{2(r+1)}$ & \pageref{p:cts}
\\
$\scrT'(s)$ & the set $T(s)\times\scrW$ & \pageref{p:scrt's}
\\
$\cW$ & window set of the cut-and-project set $\cP(\cW,\cL)$ & \pageref{p:window}
\\
$\cW_{\AB}$ & an open regular octagon with particular size and position 
& \pageref{p:waboct}
\\
$\cW_{\Gh}$ & an open regular dodecagon with particular size and position & \pageref{p:vecwgh}
\\
$\cW_{\vecw}^{(\AB)}$ & the set $(\scrW_{\AB}+\vecw)\smatr101{\sqrt2}$  
& \pageref{PABwDEF}
\\
$\cW_{\vecw}^{(\Gh)}$ & the set $(\scrW_{\Gh}+\vecw)\smatr10{\sqrt3}2$ & \pageref{GHaPexplres2}
\\
$\cW^*$ & polar set of $\cW$ & \pageref{lem:polarbody}
\\
$\xi_i$ & scaling matrix with respect to the $i$-th cusp of $\G_K$ & \pageref{def:xidef}
\\
$Y$ & the set $\left\{(y_1, y_2)\in (\R_{>0})^2\col y_2\in [1,\lambda^2)\right\}$ & \pageref{p:sety}
\\
$Y_{t_1}$ & the set $\left\{(y_1,y_2)\in(\R_{>0})^2\col 1\leq y_1/y_2<\lambda^2,\ y_1y_2>t_1\right\}$ & \pageref{Yt1DEF}
\\
$Z$ & the set $(0,\pi)\times (0, 2\pi)$ & \pageref{p:setz}
\\
$Z_J$ & the set $J\times (0, 2\pi)$ & \pageref{p:zj}
\\
$\zeta_K$ & Dedekind zeta function corresponding to $K$ & \pageref{p:dezefun}
\\
\end{longtable}
\end{footnotesize}
\end{center}

\bibliographystyle{abbrv}
\bibliography{DKbibliog}
\end{document}